\documentclass[a4paper]{article}

\usepackage{latexsym,amssymb,amsmath, mathtools, amsfonts}
\usepackage{blkarray}
\usepackage{multirow}
\usepackage{array}
\usepackage{t1enc}
\usepackage[english]{babel}
\usepackage[title]{appendix}
\usepackage{color}
\usepackage{cancel}
\usepackage[normalem]{ulem}

\newtheorem{theorem}{Theorem}[section]
\newtheorem{definition}[theorem]{Definition}
\newtheorem{proposition}[theorem]{Proposition}
\newtheorem{corollary}[theorem]{Corollary}
\newtheorem{lemma}[theorem]{Lemma}
\newtheorem{remark}[theorem]{Remark}
\newenvironment{rem}{\begin{remark} \em}{\end{remark}}

\newtheorem{example}[theorem]{Example}

\newcommand{\RR}{\mathbb R}

\newcommand{\C}{\mathcal C}
\newcommand{\G}{\mathcal G}
\newcommand{\OO}{\mathcal O}

\newcommand{\U}{\mathcal U}
\newcommand{\V}{\mathcal V}

\newcommand{\W}{\mathcal W}
\newcommand{\X}{\mathcal X}

\newcommand{\HH}{\mathcal H}

\newcommand{\PTO}{\mathcal T}

\newcommand{\TR}{\mathcal TR}

\newcommand{\ua}{\underline{a}}
\newcommand{\ub}{\underline{b}}
\newcommand{\uI}{\underline{\I}}
\newcommand{\rr}{\underline{r}}
\newcommand{\kk}{\underline{k}}
\newcommand{\mm}{\underline{m}}
\newcommand{\ww}{\underline{w}}

\newcommand{\ualpha}{\underline{\alpha}}

\newcommand{\nsc}{\ensuremath{\text{\normalsize $0$}}}

\def\Gl{\mathop{\rm Gl}\nolimits}

\def\diag{\mathop{\rm diag}\nolimits}
\def\rank{\mathop{\rm rank}\nolimits}
\def\Ima{\mathop{\rm Im}\nolimits}

\def\re{\mathop{\rm re}\nolimits}
\def\tr{\mathop{\rm tr}\nolimits}

\def\Ker{\mathop{\rm Ker}\nolimits}

\newcommand{\la}{\ensuremath{\lambda}}
\newcommand{\wt}{\widetilde}
\newcommand{\wh}{\widehat}

\newcommand{\PP}{\mathcal P}

\newcommand{\CC}{\mathbb C}
\newcommand{\wC}{\widetilde{C}}
\newcommand{\whJ}{\widehat{J}}
\newcommand{\wtJ}{\widetilde{J}}
\newcommand{\whX}{\widehat{X}}

\newcommand{\whW}{\widehat{W}}
\newcommand{\wtW}{\widetilde{W}}

\newcommand{\T}{\mathcal T}
\newcommand{\I}{\mathcal I}

\newcommand{\dZ}{Z^{\diamond}}

\newcommand{\A}{\mathcal A}
\newcommand{\R}{\mathcal R}
\newenvironment{smallarray}[1]
 {\null\,\vcenter\bgroup\scriptsize
  \renewcommand{\arraystretch}{0.7}
  \arraycolsep=.13885em
  \hbox\bgroup$\array{@{}#1@{}}}
 {\endarray$\egroup\egroup\,\null}

\makeatletter
\def\adots{\mathinner{\mkern1mu\raise\p@
\vbox{\kern7\p@\hbox{.}}\mkern2mu
\raise4\p@\hbox{.}\mkern2mu\raise7\p@\hbox{.}\mkern1mu}}
\makeatother

\title{A Local Parametrization of the State-Feedback Matrices in the Pole Assignment Problem}
\author{
\sc I. Baraga\~na
\footnotemark \hspace*{1mm}
\thanks{e-mail: itziar.baragana@ehu.eus.
Supported by the Agencia Estatal de Investigaci\'on of Spain MCIN/AEI/10.13039/501100011033/
 and by ``ERDF A way of making Europe'' of the ``European Union'' through grant PID2021-124827NB-I00.}
\and
\sc F. Puerta
\footnotemark \hspace*{1mm}
\thanks{e-mail:  ferran.puerta@upc.edu.}\
\and
\sc I. Zaballa
\footnotemark \hspace*{1mm}
\thanks{Corresponding author, e-mail: ion.zaballa@ehu.eus.
Supported by the Agencia Estatal de Investigaci\'on of Spain MCIN/AEI/10.13039/501100011033/
 and by ``ERDF A way of making Europe'' of the ``European Union'' through grant PID2021-124827NB-I00.}\\ \\
\addtocounter{footnote}{-6}
\footnotemark \hspace*{2mm}
Departamento de Ciencia de la
Computaci\'on e IA\\
University of the Basque Country UPV/EHU\\
Apdo. 649.
20080 Donostia-San Sebasti\'an. Spain
\vspace{7pt}\\
\addtocounter{footnote}{+1}
\footnotemark \hspace*{2mm}
Departament de Matem\`{a}tica Aplicada I.\\
E.T.S. Enginyeria Industrial de Barcelona. UPC\\
Institut d'Organitzaci\'o i Control(IOC) UPC\\
Diagonal 647. 08028 Barcelona. Spain
\vspace{7pt}\\
\addtocounter{footnote}{+1}
\footnotemark \hspace*{2mm}
Departamento de Matem\'aticas.\\
University of the Basque Country UPV/EHU\\
Apdo. 644.
48080 Bilbao. Spain
}

\begin{document}
\maketitle

\begin{abstract}
Given a controllable system $(F,G)$, 
a local parametrization is obtained for the set  of  the feedback gain matrices $K$ such that the state matrix, $F+GK$, of the
closed-loop system is in a prescribed similarity class.  It is shown that this set can be endowed
with the structure of a differentiable manifold whose dimension is also computed. Then a local parametrization
and a local system of coordinates are obtained using a diffeomorphism between this set of state-feedback matrices
and the orbit space of a set of truncated observability matrices via the action of a Lie group.
\end{abstract}



\section{Introduction}
\label{secintroduction}

One of the most classic problems in the literature on linear control systems is the \textit{pole assignment problem}
(also known as the \textit{pole placement problem}).
Assume that we are given a linear, time-invariant control system
\[
\dot{x}(t)=Fx(t)+Gu(t),\quad F\in\RR^{n\times n}, G\in\RR^{ n\times m}
\]
where $\RR$ is the field of real numbers and $\RR^{p\times q}$ the set of $p\times q$ matrices over $\RR$.
The system above will be identified with the pair of matrices $(F,G)$. 
Assume that we are also given a  self-conjugate 
sequence $\Lambda$ of $n$ complex numbers (i.e., such that $\overline{\Lambda}=\Lambda$), the pole placement problem
is to find a state-feedbak matrix $K$ such that the eigenvalues of $F+GK$ are the complex numbers of $\Lambda$.
If the given  system $(F,G)$ is controllable, it is well known that  such a state-feedback matrix always exists.

A more general and difficult problem,  sometimes called \textit{the general pole assignment problem}, consists of
assigning not only the eigenvalues but the complete similarity class for $F+GK$. Recall that two matrices 
$A, A' \in \RR^{n\times n}$ are said to be {\em similar} if $A'=T^{-1}AT$ with
$T\in \Gl(n)$, the general linear group of all invertible matrices in $\RR^{n\times n}$.
A complete system of invariants for matrix similarity  is given by the
invariant polynomials (see, for example, \cite[Ch. 6]{Gant59}). 
Thus,  the orbit of $A$ under the similarity action only depends on its invariant polynomials 
$\ualpha: \alpha_1(s)\mid \dots \mid \alpha_n(s)$ and will be denoted by $\OO(\ualpha)$.
Therefore, given the system $(F,G)$, the general pole assignment
problem is to find $K\in\RR^{m\times n}$ such that
$F+GK\in\OO(\ualpha)$ for a given sequence of monic polynomials
$\ualpha: \alpha_1(s)\mid \dots \mid \alpha_n(s)$ such that $\sum_{i=1}^n\deg(\alpha_i(s))=n$.

Necessary and sufficient conditions for such a matrix $K$ to exist were given by Rosenbrock  (\cite{Rosen70})
when the system $(F,G)$ is controllable   and by Zaballa (\cite{Za89}) when this controllability restriction is removed 
(see Proposition \ref{prop.necss} below). In both cases
the proofs are constructive (see also \cite{BoLoBa97,Flamm80} for alternative geometric proofs of Rosenbrock's result).
In general, if  there is a matrix $K$ such that $F+GK$ is in a prescribed  similarity class, it is not unique. We aim to 
\textit{locally parameterize} this set of feedback gain matrices. To achieve this, we first study
the geometry of the set of state-feedback matrices for the given controllable system $(F,G)$:
\begin{equation}\label{eqH}
\HH_{(F,G)}=\{K\in  \RR^{m\times n}\; : \; F+GK\in \OO(\ualpha)\}.
\end{equation}
Specifically, we aim to show that $\HH_{(F,G)}$
can be endowed with the structure of a differentiable manifold making it is an (immersed) submanifold of $\RR^{m\times n}$. 
In addition, its dimension will be computed.  It is a well-known general result that immersed submanifolds admit 
\textit{local parametrizations} and \textit{local systems of coordinates} (see, for example\cite[Lemma 8.18]{Lee03}).
In this paper, a local parametrization for $\HH_{(F,G)}$ is obtained using a diffeomorphism of this set and an orbit
space of matrices via the action of a Lie group.

Having a local parametrization defined in $\HH_{(F,G)}$ may be interesting in several respects. For example,
if $K\in\HH_{(F,G)}$ and $K'$ is a small, arbitrary
perturbation of $K$, then $K'$ may not be in $\HH_{(F,G)}$. A differentiable structure and a local parametrization
in $\HH_{(F,G)}$ can be used to determine the possible perturbations of $K$ that remain in $\HH_{(F,G)}$; i.e., those such that $F+GK$ and $F+GK'$ have the same invariant polynomials. On the other hand, a local coordinate system provides the
minimum number of parameters required to fully describe the perturbed feedback matrices $K'\in\HH_{(F,G)}$. 
A local parametrization may also be useful to determine the optimal, in some sense, feedback gain matrix
$K$.  This idea is explored  in \cite{Ma19} in relation to the \textit{optimal pole placement problem}  (see also
the references therein). When tackling this problem  the starting point is usually to parametrize the set of allowable
matrices $K$. 

The rest of the paper is organized as follows. Section \ref{sepreliminaries} presents the notation and the necessary preliminary 
results. In particular, it reviews the complex and real Jordan and Weyr canonical forms and their centralizers, recalls the feedback 
equivalence of linear control systems and the statement of Rosenbrock's theorem on pole assignment, and introduces a new 
Brunovsky canonical form. In Section \ref{secgeomstrh} we study the geometry of  $\HH_{(F,G)}$ in \eqref{eqH} for a given 
controllable system $(F,G)$. It is shown that this set is a differentiable submanifold  of $\RR^{m\times n}$ and we also
compute its dimension. The final three sections are dedicated to obtaining a local parametrization and a local system of 
coordinates for $\HH_{(F,G)}$.  First, we will show in Section \ref{secmanifpc} that $\HH_{(F,G)}$ is diffeomorphic to an orbit space 
of truncated observability matrices via the action of a Lie group. This Lie group consists of the invertible matrices that commute 
with the state matrix of the observability system from which the truncated observability matrix is derived.
Using this action a unique local reduced form is found in
Section \ref{secreducedform} for each orbit. This reduced form provides a local parametrization of the orbit space of
truncated observability matrices that is translated in Section \ref{secparameterization} to a local parametrization and
a local system of coordinates of $\HH_{(F,G)}$ by means of the diffeomorphism found in Section \ref{secmanifpc}.
An example that illustrates the whole process is provided in Section \ref{secparameterization}.

\section{Notation and Preliminary results}
\label{sepreliminaries}


\subsection{Partitions}
\label{sec.partitions}

If $s$ and $p$ are positive integers ($0<s\leq p$),
$
Q_{s,p}:=\{(i_1, \dots, i_s)\;:\; 1\leq i_1<\dots <i_s\leq p\}
$
and $Q_{0, p}:=\{\emptyset\}$. 
If $A\in \RR^{n \times m}$, $I\in Q_{s, n}$ and
$J\in Q_{r, m}$ then $A(I, J)$ will denote the $s\times r$ submatrix of  $A$
formed by the rows in $I$ and the columns in $J$; that is, if $I=(i_1,\ldots, i_s)$,
$J=(j_1,\ldots, j_r)$ and $B=A(I,J)\in\RR^{s\times r}$ then $b_{k\ell}=a_{i_kj_\ell}$,
$k=1,\ldots, s$ and $\ell=1,\ldots, r$. Similarly, $A(I, :)\in \RR^{s\times m}$ and
$A(:, J)\in \RR^{n\times r}$ are the submatrices of $A$ formed by the rows in $I$
(and all columns) and the columns in $J$ (and all rows), respectively.
If $p\leq q$ are integers the symbol $p:q$
denotes the sequence $(p, p+1,\ldots, q)$.

It is well known that a \textit{partition} is a finite or infinite sequence $\ua=(a_1, a_2,\ldots, )$
of nonnegative integers almost all zero.  In this manuscript we will only use partitions whose components
are arranged in nonincreasing order.  The sequence $\ua$ is said to be a partition of $n$ if
$n=\sum_{i\geq 1} a_i$.
If $\ua$ and $\ub$ are partitions then $\ua + \ub=(a_1+b_1, a_2+b_2,\ldots)$ and
$\ua\cup\ub$ is the partition whose elements are those of $\ua$ and those of $\ub$
(with possible repetitions) reordered so that they do not increase. If $\ua$ is a partition
of $n$ and we define $b_i=\#\{j: a_j\geq i\}$, where $\#$ stands for cardinality, then
$\ub=(b_1,b_2,\ldots)$ is said to be the \textit{conjugate} or \textit{dual} partition of $\ua$
and $\sum_{i\geq 1}b_i=n$.
It is well-known (see, for example, \cite[Section 7.B]{MOA11}), that the conjugation 
of non-increasingly ordered partitions is an involution; i.e., if $\ub$ is the conjugate 
partition of $\ua$ then the latter is the conjugate partition of the former. On the other hand,
if $\ua=(a_1,a_2,\ldots, a_n)$ and $\ub=(b_1,b_2,\ldots, b_n)$ are finite partitions with
$a_1\geq a_2\geq \cdots\geq a_n\geq 0$ and $b_1\geq b_2\geq\cdots\geq b_n\geq 0$,
following \cite{HaLiPo67} (see also \cite[Chapter 1]{MOA11}),
we say that $\ua$ is majorized by $\ub$, and we write $\ua \prec \ub$, if 
\begin{equation}\label{eq.maj1}
\sum_{j=1}^{k}a_{j} \leq \sum_{j=1}^{k} b_{j}\quad 1 \leq k \leq m,\quad \text{ and }\quad
\sum_{i=1}^n a_i = \sum_{i=1}^n b_i.
\end{equation}
It should be noted that we can always assume that any two finite partitions have the same number of components
by adding or deleting zeros.
The following result is well-known (see, for example, \cite[p. 6,7]{Mac95}):
\begin{proposition}\label{prop.uncmajc}
If $\ua$ and $\ub$ are partitions then
\begin{itemize}
\item[(i)] $(\ua+\ub)^\ast= \ua^\ast\cup\ub^\ast$
\item[(ii)] $\ua \prec \ub\quad\Leftrightarrow \quad \ub^\ast \prec  \ua^\ast$
\end{itemize}
where $\ua^\ast$, $\ub^\ast$ and  $(\ua+\ub)^\ast$ are the conjugate partitions of  $\ua^\ast$, $\ub^\ast$
and  $(\ua+\ub)^\ast$, respectively.
\end{proposition}
A proof for item (ii) above can also be found in\cite[Section 7.B]{MOA11}.

\subsection{Real Jordan and Weyr Canonical Forms}\label{sec.realjordanweyr}

%

A canonical form for the similarity of real square matrices is the {\em real Jordan canonical form}
(see, for instance, \cite[Theorem 6.7.1]{LaTi85} or \cite[Theorem[12.2.2]{GoLaRo86}).
Let $A\in \RR^{n\times n}$ be a matrix with invariant polynomials 
$\ualpha: \alpha_1(s)\mid \dots \mid \alpha_n(s)$, let
$\lambda_1, \dots, \lambda_t\in\CC$ be the distinct eigenvalues of $A$
and split $\alpha_i(s)$ as a product of powers of irreducible polynomials
\[
\alpha_{n-i+1}(s)=(s-\la_1)^{m_{1i}}(s-\la_2)^{m_{2i}}\cdots (s-\la_t)^{m_{t\,i}},\quad 1\leq i\leq n.
\]
Then $m_{i1}\geq m_{i2}\geq\cdots\geq m_{i\,w_i}>0=m_{i\,w_{i+1}}=\cdots=m_{i\,n}$  for $i=1,\ldots, t$ and the sequence
$\big((m_{11}, \dots, m_{1w_1}), \dots, (m_{t1}, \dots, m_{t\,w_t})\big)$ is known as \textit{the Segre characteristic} of $A$.
If $n_i$ is the algebraic multiplicity of $\la_i$ then the Segre
characteristic of $A$, $(m_{i1}, \dots, m_{i\,w_i})$,  is a
partition of $n_i$. Its conjugate partition will be denoted by $(w_{i1},\ldots, w_{i\,m_i})$  (observe
that $w_{i1}=w_i$ and $m_{i1}=m_i$).The sequence of partitions
$\big((w_{11}, \dots, w_{1m_1}), \dots, (w_{t1}, \dots, w_{t\,m_t})\big)$ is called \textit{the Weyr characteristic} of $A$
(see, for example, \cite{Shapi99}). Each of these characteristics has an associated canonical form: 
the Jordan canonical form for the Segre characteristic and the Weyr canonical form for the Weyr characteristic.
If the eigenvalues are all real; i. e., $\la_1,\ldots, \la_t\in\RR$, then the well-known Jordan canonical form of $A$ is
$J= \diag\big(J(\la_1), \dots, J(\la_t)\big)$ where
\begin{equation}\label{eq.defJk}
 \begin{array}{c}
 J(\la_i)=\diag\big(J_{1}(\lambda_i), \dots, J_{w_i} (\lambda_i)\big)\in\RR^{n_i\times n_i}, \quad 1\leq i\leq t\\
J_{k}(\lambda_i)=\begin{bmatrix}
\lambda_i&1&\cdots &0&0\\
0&\lambda_i&\cdots &0&0\\
\vdots&\vdots&\ddots &\ddots&\vdots\\
0&0&\dots &\lambda_i&1\\
0&0&\dots &0&\lambda_i\\
\end{bmatrix}\in \RR^{m_{ik}\times m_{ik}},\quad 1\leq k\leq w_i.
\end{array}
\end{equation}
The Weyr canonical form is, perhaps, less known but it is more convenient for our developments
(see \cite[Chapter 2]{OmClVi11} or \cite{Shapi99}; we willl use the notation of the latter):
$W= \diag\big(W(\la_1), \dots, W(\la_t)\big)$ where, for $1\leq i\leq t$,
\begin{equation}\label{eq.defWk}
W(\la_i)=\begin{bmatrix}\la_i I_{w_{i1}}& I_{w_{i1},w_{i2}} &\cdots&0&0\\
0 & \la_i I_{w_{i2}}& \cdots&0&0\\
\vdots&\vdots&\ddots&\ddots&\vdots\\
 0& 0& \cdots&\la_i I_{w_{i\,m_i-1}}&I_{w_{i\,m_i-1},w_{i\,m_i}} \\
0&0&\cdots&0&\la_i I_{w_{i\,m_i}}\end{bmatrix}\in\RR^{n_i\times n_i}, 
\end{equation}
and for $p\geq q$
\begin{equation}\label{eq.defIpq}
I_{p,q}=\begin{bmatrix} I_q\\0\end{bmatrix}\in\RR^{p\times q}.
\end{equation}

The Jordan and Weyr canonical forms of $A$ are closely related. The following lemma shows that they can be obtained
from each other by a permutation similarity. This is a well-known result
(see, for example, \cite[Chapter 2]{OmClVi11}
or \cite{MiMoRo20} for a generalization to arbitrary fields). We offer a simple proof that will be used to define
the real Weyr canonical form.
 
\begin{lemma}(\cite[Chapter 2]{OmClVi11})\label{lem.FCJW}
Let $H\in\RR^{n\times n}$ be a matrix with $\la_0\in\RR$ as its only eigenvalue and let $(m_1,\ldots, m_w)$
and $(w_1,\ldots, w_m)$  be their Segre and Weyr characteristics, respectively, where $m=m_1$ and $w=w_1$.
Let $J$ and $W$ be the Jordan and Weyr canonical forms of $H$ and define
\begin{equation}\label{eq.defsi}
 s_i=w_1+\ldots+w_i,\quad i=1,\ldots, m.
\end{equation}
Let $e_k$ be the $k$-th column of $I_{n}$, $k=1,\ldots, n$, and
\begin{equation}\label{eq.defQ}
\begin{array}{c}
Q=\begin{bmatrix}Q_1^T&Q_2^T&\cdots& Q_w^T\end{bmatrix}^T,\\
Q_i=\begin{bmatrix}
e_i& e_{s_1+i}& \cdots & e_{s_{m_i-1}+i}\end{bmatrix}^T\in\RR^{m_i\times n},\;\;1\leq i\leq w
\end{array}
\end{equation}
where  ${}^T$ stands for transpose. Then $W=Q^TJQ$ .
\end{lemma}

\textbf{Proof.} Since $J$ and $W$ are the Jordan and Weyr canonical forms of $H$, they are similar.
Put $\wtJ=\la_0I_n-J$ and $\wtW=\la_0I_n-W$.
The matrices $T\in\Gl(n)$ such that $\wtJ T=T\wtW $ form an open set of the following linear subspace of
dimension $wn$:
\begin{equation}\label{eq.defPTO}
\PTO(\wtW,\mm)=\left\{T=\begin{bsmallmatrix}T_1\\T_2\\\vdots\\\\T_w\end{bsmallmatrix},
T_i=\begin{bsmallmatrix} t_i\\t_i \wtW\\\vdots\\\\ t_i \wtW^{m_i-1}\end{bsmallmatrix}, t_i\in\RR^{1\times n},
1\leq i\leq w\right\}
\end{equation}
It follows from  the definition of $\wtW$ that, for $1\leq i \leq w$, $1\leq j \leq m_i-1$, 
$ e^T_{s_{j-1}+i}\wtW=e^T_{s_{j}+i}$ ($s_0=0$) and so $e_{s_{j}+i}^T=e_i^T\wtW^j$. Hence,
$Q_i=\begin{bsmallmatrix}
  e_i^T\\ e_{i}^T\wtW\\  \vdots \\ e_{i}^T\wtW^{m_i-1}\end{bsmallmatrix}$, $1\leq i\leq w$ and 
  $Q\in\PTO(\wtW,\mm)$. Therefore $Q^TJQ=W$ as desired.\hfill $\Box$




\medskip
Assume now that $A\in\RR^{n\times n}$ has real and nonreal eigenvalues:
$\lambda_1, \dots, \lambda_p\in \RR$ and $\lambda_{p+1}, \dots, \lambda_{p+q}, \overline{\lambda}_{p+1}, 
\dots, \overline{\lambda}_{p+q} \in \CC\setminus \RR$ where
$\overline{\lambda}_i$ stands for the complex conjugate of $\lambda_i$ and put $t=p+2q$.  
In this case the prime factorization of $\alpha_{n-i+1}(s)$, $i=1,\ldots, n$, would be
\begin{equation}\label{eq.primefacalfa}
\alpha_{n-i+1}(s)=(s-\la_1)^{m_{1i}}\cdots (s-\la_p)^{m_{pi}} (s^2+c_1s+c_2)^{m_{p+1\,i}}
\cdots(s^2+c_{2q-1}s+c_{2q})^{m_{p+q\,i}},
\end{equation}
where we can assume without loss of generality that
$s^2+c_{2k-1}s+c_{2k}= (s-\la_{p+k})(s-\overline{\la}_{p+k})$, $k=1,\ldots, q$. Then
$m_{i1}\geq m_{i2}\geq\cdots\geq m_{i\,w_i}>0=m_{i\,w_{i+1}}=\cdots=m_{in}$  for $i=1,\ldots, t=p+2q$ and 
$\big((m_{11}, \dots, m_{1w_1}), \dots, (m_{t1}, \dots, m_{t\,w_t})\big)$ is the  Segre characteristic of $A$. Note
that for $k=p+1,\ldots, t=p+2q$, the Segre characteristics of $A$ for the eigenvalues
$\la_k$ and $\overline{\la}_k$ coincide; i.e., $m_{kj}=m_{k+q\,j}$ for $j=1,\ldots, w_k$. Then the
Real Jordan canonical form of $A$ is  (see \cite[Theorem 6.7.1]{LaTi85} or \cite[Theorem 12.2.2]{GoLaRo86}):

\begin{equation}\label{eqAJ}
  J_R=\diag\left(J(\la_1), \dots, J(\la_p), \whJ(\la_{p+1},\overline{\la}_{p+1}), \dots, 
  \whJ(\la_{p+q},\overline{\la}_{p+q})\right),
\end{equation}
where $J(\la_i)$ are the matrices of \eqref{eq.defJk},
$$
\whJ(\la_j,\overline{\la}_j))=\diag\left(\whJ_{1} (\lambda_j,\overline{\lambda}_{j}), 
\dots, \whJ_{w_j} (\lambda_j,\overline{\lambda}_{j})\right), 
\quad p+1\leq j\leq p+q,
$$
and, if $\lambda_j=a_j+b_ji \in  \CC\setminus \RR$ then
\begin{equation}\label{eq.defJklajlajc}
\whJ_{k}(\lambda_j,\overline{\lambda}_j)=
\begin{bmatrix}
B_j&I_2&0&\dots &0&0\\
0&B_j&I_2&\dots &0&0\\
0&0&B_j&\dots &0&0\\
\vdots&\ddots&\ddots &\ddots &\vdots&\vdots \\
0&0&0&\dots &B_j&I_2\\
0&0&0&\dots &0&B_j\\
\end{bmatrix}\in \RR^{2m_{jk}\times 2m_{jk}}, \quad 1\leq k\leq w_i,
\end{equation}
with $B_j=\begin{bsmallmatrix}a_j&b_j\\-b_j&a_j\end{bsmallmatrix}$.

Now, for $i=1,\ldots, t=p+2q$, let $(w_{i1},\ldots, w_{i\,m_i})$ be the conjugate partition of
$(m_{i1}, \dots, m_{i\,w_i})$ ($m_i=m_{i1}$ and $w_i=w_{i1}$). Then, as in the case when all eigenvalues are real,
$\big((w_{11}, \dots, w_{1m_1}), \dots, (w_{t1}, \dots, w_{t\,m_t})\big)$ is the Weyr characteristic
of $A$. Weyr canonical forms for matrices over arbitrary fields where studied in \cite{MiMoRo20}.
However an explicit definition of a real Weyr canonical form is not provided. It can be obtained, after some
manipulations,  by applying the technique of  Section 3 in that paper to the \textit{Generalized Jordan canonical 
form of the first kind} in \cite[Theorem 2.7]{MiMoRo20} for matrices with real entries . We take a more direct approach generalizing
Lemma \ref{lem.FCJW} to the case of nonreal eigenvalues.

Recall that if $X\in\RR^{m\times n}$ the Kronecker product
$I_p\otimes X=\diag(\overbrace{X,\ldots, X}^{p})$. For
notational simplicity we will use the notation
\begin{equation}\label{eq.defX(n)}
X^{(p)}=I_p\otimes X=\diag(\overbrace{X,\ldots, X}^{p}.
\end{equation}

\begin{lemma}\label{lem.FCJWC}
Let $H\in\RR^{2n\times 2n}$ be a matrix with $\la_0,\overline{\la}_0\in\CC\setminus\RR$ as its only eigenvalues 
and let $(m_1,\ldots, m_w)$ and $(w_1,\ldots, w_m)$ be the common Segre and Weyr characteristics of $H$
for both $\la_0$ and $\overline{\la}_0$. 
Assume that $\la_0=a_0+b_0i$ and let $B_0=\begin{bsmallmatrix}a_0 &b_0\\-b_0&a_0\end{bsmallmatrix}$.
Let the real Jordan canonical form of $H$ be
$\wh J(\la_0,\overline{\la}_0)=\diag\left( \whJ_{1} (\lambda_0,\overline{\lambda}_{0}), 
\dots, \whJ_{w} (\lambda_0,\overline{\lambda}_{0})\right)$ where $\whJ_k(\lambda_0,\overline{\lambda}_{0})$
is the matrix of \eqref{eq.defJklajlajc} with $j=0$ and of size $2m_k\times 2m_k$, $k=1,\ldots,w$. Define
\begin{equation}\label{eq.blockWeyrcomplex}
 \whW(\la_0,\overline{\la}_0)=
  \begin{bmatrix}
    B_0^{(w_{1})}&I_{2w_{1},2w_{2}}&\dots &0&0\\
    0&B_0^{(w_{2})}&\dots &0&0\\
    \vdots&\vdots&\ddots &\ddots&\vdots\\
    0&0&\dots &B_0^{(w_{m-1})}&I_{2w_{m-1},2w_{m}}\\
    0&0&\dots&0&B_0^{(w_m)}
\end{bmatrix},
\end{equation}
For $i=1,\ldots,m$, let $s_i$ be the
positive integer of \eqref{eq.defsi}.
Let $E_k=\begin{bmatrix}e_{2k-1}&e_{2k}\end{bmatrix}\in\RR^{2n\times 2}$ where $e_k$ is the $k$-th column
of $I_{2n}$, $k=1,\ldots, n$, and
\begin{equation}\label{eq.defQC}
\begin{array}{c}
Q=\begin{bmatrix}Q_1^T&Q_2^T&\cdots& Q_w^T\end{bmatrix}^T,\\
Q_i=\begin{bmatrix}
E_i& E_{s_1+i}& \cdots & E_{s_{m_i-1}+i}
\end{bmatrix}^T\in\RR^{2m_i\times 2n},\;\;1\leq i\leq w.
\end{array}
\end{equation}
Then $\whW(\la_0,\overline{\la}_0)=Q^T\wh J(\la_0,\overline{\la}_0)Q$.
\end{lemma}

\textbf{Proof.} Let $\wtJ=\wh J(\la_0,\overline{\la}_0)-B_0^{(n)}$ and
 $\wtW=\whW(\la_0,\overline{\la}_0)-B_0^{(n)}$. Then
the only eigenvalue of $\wtJ$ and $\wtW$ is $0$, the Segre characteristic of $\wtJ$ is $\mm\cup\mm=
(m_1,m_1,m_2,m_2,\ldots, m_w,m_w)$ and the Weyr characteristic of $\wtW$ is $\ww+\ww=
(2w_1,2w_2,\ldots,2w_m)$. Since $\mm\cup\mm$ and $\ww+\ww$ are conjugate partitions, $\wtJ$ and $\wtW$
are similar matrices. As in Lemma \ref{lem.FCJW}, the set of matrices $T\in\Gl(2n)$ such that $\wtJ T=T \wtW$
form an open set of a linear subspace of dimension $4wn$:
\[
\TR(\wtW,\mm+\mm)=\left\{T=\begin{bsmallmatrix}T_1\\T_2\\\vdots\\\\T_w\end{bsmallmatrix},
T_i=\begin{bsmallmatrix} X_i\\X_i \wtW\\\vdots\\\\ X_i\wtW^{m_i-1}\end{bsmallmatrix}, X_i\in\RR^{2\times 2n},
1\leq i\leq w\right\}.
\]
Now, as in the proof of Lemma \ref{lem.FCJW}, for $1\leq i \leq w$, $1\leq j \leq m_i-1$, 
$ E^T_{s_{j-1}+i}\wtW=E^T_{s_{j}+i}$ ($s_0=0$), $E_{s_{j}+i}^T=E_i^T\wtW^j$,
$Q_i=\begin{bsmallmatrix}
E_i^T\\ E_{i}^T\wtW\\  \vdots \\ E_{i}^T\wtW^{m_i-1}\end{bsmallmatrix}$, $1\leq i\leq w$, 
$Q\in\PTO(\wtW,\mm)$ and so $Q^T\wtJ Q=\wtW$. Since
$B_0^{(n)}$ is a block diagonal matrix with repeated $2\times 2$ diagonal blocks,
it is not difficult to see that $Q^TB_0^{(n)}Q=
B_0^{(n)}$ and the Lemma follows.\hfill $\Box$

\medskip
In general, if $A\in\RR^{n\times n}$, there is permutation matrix $Q\in\RR^{n\times n}$
such that
\begin{equation}\label{eq.defWR}
W_R=Q^TJ_RQ= \diag(W(\lambda_1), \dots, W(\lambda_p), 
\whW(\lambda_{p+1}, \overline{\lambda}_{p+1}), \dots, \whW(\lambda_{p+q}, \overline{\lambda}_{p+q})),
\end{equation}
\[
\whW(\lambda_j, \overline{\lambda}_j)=
\begin{bmatrix}
    B_j^{(w_{j1})}&I_{2w_{j1},2w_{j2}}&\dots &0&0\\
    0&B_j^{(w_{j2})}&\dots &0&0\\
    \vdots&\vdots&\ddots &\ddots&\vdots\\
    0&0&\dots &B_j^{(w_{j\,m_j-1})}&I_{2w_{j\,m_j-1},2w_{j\,m_j}}\\
    0&0&\dots&0&B_j^{(w_{j\,m_j})}
\end{bmatrix}.
\]
The matrix $W_R$ will be called the \textit{real Weyr canonical form} of $A$.

\subsection{The centralizer}\label{sec.centralizer}
Given a matrix $A\in \RR^{n\times n}$, we denote by $C_A$ the {\em centralizer} of $A$, i.e.
$$
C_A=\{X\in \RR^{n\times n}\, : \, XA=AX\}.
$$
It is well-known (see, for example, \cite[Ch. 8]{Gant59}) that if $\alpha_1(s)\mid \dots \mid \alpha_n(s)$
are the invariant polynomials of $A$ then  $C_A$ is a  subspace of $\RR^{n\times n}$ of  $\dim C_A=N$,
where
\begin{equation}\label{eqqN}
N=\deg(\alpha_n)+3\deg(\alpha_{n-1})+5\deg(\alpha_{n-2})\dots =
\sum_{k=1}^n(2k-1)\deg(\alpha_{n-k+1}).
\end{equation}

Let  $\wC_A$ be the  subgroup of $\Gl(n)$ formed by the invertible matrices of $C_A$, 
$\wC_A=\{X\in \Gl(n) \, : \, X\in C_A\}.$
Then $\wC_A$ is an open subset of a linear manifold and $\dim \wC_A=N$.
It turns out that  (see, for instance, \cite[Proposition 3.2]{BaPu19}, \cite[Theorem 2.1]{DeEd95} or
the proof of Theorem 9.16 in \cite{Lee03}), also $\OO(\ualpha)$ is a differentiable manifold of codimension $N$.

\medskip
Let $A=J_R$ be the matrix in (\ref{eqAJ}). Then it is easily seen that
$X\in C_A$ if and only if $X=\diag(X_1, \dots, X_{p}, \whX_{p+1}, \dots, \whX_{p+q})$ where
$X_i \in C_{J(\la_i)}$, $1\leq i \leq p$ and $\whX_i \in C_{\whJ(\la_i,\overline{\la}_i)}$, $p+1\leq i \leq p+q$.
Similarly, if $A=W_R$ is the matrix of \eqref{eq.defWR} then $X\in C_A$ if and only if
$X=\diag(X_1, \dots, X_{p}, \whX_{p+1}, \dots, \whX_{p+q})$ where
$X_i \in C_{W(\la_i)}$, $1\leq i \leq p$ and $\whX_i \in C_{\whW(\la_i,\overline{\la}_i)}$, $p+1\leq i \leq p+q$.

\medskip
 Let $\la_0 \in \RR$  be an eigenvalue of $A$ and let $(m_1,\ldots, m_w)$ be its Segre Characteristic. If
 $J(\la_0)=\diag \left(J_1(\la_0),\dots ,  J_w(\la_0)\right)$  is the block associated to $\la_0$ in the Jordan canonical
 form  of $A$, then the characterization of
the centralizer of $J(\la_0)$,  $\C_{J(\la_0)}$, can be found in many books (for example in
\cite[Ch. 8]{Gant59}, \cite[Theorem 12.4.2]{GoLaRo86} or \cite[Ch. 12]{LaTi85}). We are interested in the less known
characterization of the centralizer of $W(\la_0)$,
the block associated to $\la_0$ in the Weyr canonical form of $A$. Taking into account that
$Q^TJ(\la_0) Q=W(\la_0)$ where $Q$ is the matrix of \eqref{eq.defQ}, $X\in C_{J(\la_0)}$ if and only if
$Q^TXQ\in C_{W( \la_0)}$. Using this property we get

\begin{lemma}\label{lemmazreal}
 Let $\la_0 \in \RR$ and $W(\la_0)$ be the matrix of \eqref{eq.defWk} with $i=0$ and Weyr characteristic
 $(w_1,w_2,\ldots, w_m)$. Then $Y\in C_{W(\la_0)}$ if and only if
 \begin{equation}\label{eq.QXQr}
 Y=
\begin{bmatrix}
  Y_{11}&Y_{12}& \dots & Y_{1,m}\\
   0&Y_{22}& \cdots & Y_{2 m}\\
  \vdots&\vdots&\ddots&\vdots&\\
  0& 0&\cdots & Y_{mm}\\
\end{bmatrix}\
 \end{equation}
where
\begin{itemize}
\item[(i)]
\begin{equation}\label{eq.QYQm}
Y_{1j}=\begin{bmatrix}
D_{11}^{(j)} & D_{12}^{(j)} & \cdots & D_{1\,m-j+1}^{(j)} \\
\vdots & \vdots &  & \vdots\\
D_{j1}^{(j)} &D_{j2}^{(j)} & \cdots & D_{j\,m-j+1}^{(j)}\\
0 &D_{j+1\,2}^{(j)} & \cdots & D_{j+1\,m-j+1}^{(j)}\\
\vdots &\vdots &&\vdots\\
0 & 0 & \cdots & D_{m\,m-j+1}^{(j)}
\end{bmatrix},
\end{equation}
and
\begin{equation}\label{eq.qxqdijr}
D^{(j)}_{i,k}\in \RR^{(\tau_i-\tau_{i-1})\times (\tau_k-\tau_{k-1})}, \quad
1\leq i, j\leq m,\, \max\{i-j+1, 1\}\leq k \leq m-j+1,
\end{equation}
with $\tau_i=w_{m-i+1}$,  $0\leq i \leq m$ ($w_{m+1}=0$).
\item[(ii)] For $1\leq i\leq  j \leq m-1$
\begin{equation}\label{eq.qxqyijr}
Y_{i+1\, j+1}=I_{w_{i}, w_{i+1}}^T Y_{i\,j}I_{w_{j}, w_{j+1}}.
 \end{equation}
\end{itemize}
\end{lemma}
\bigskip

\begin{rem}\label{rem.comweyr}{\rm
Condition \eqref{eq.qxqyijr} means that $Y_{ij}$ has de form $Y_{ij}=\begin{bsmallmatrix} Y_{i+1\,j+1} & \ast
\\0 & \ast\end{bsmallmatrix}$. So, all distinct parameters of $Y$ are concentrated in $Y_{1j}$, $1\leq j\leq m$. The
number of parameters in $Y_{1j}$ is $(w_j-w_{j+1})w_1+(w_{j+1}-w_{j+2})w_2+\ldots+ (w_{m-1}-w_m)w_{m-j}+w_m w_{m-j+1}$.
Thus the number of distinct parameters in $Y$ is 
\[
\sum_{j=1}^m(w_j-w_{j+1})w_1+\sum_{j=1}^m(w_{j+1}-w_{j+2})w_2+\ldots+\sum_{j=1}^m w_m w_{m-j+1}=
w_1^2+w_2^2+\ldots+w_m^2.
\]
This is, actually, the value of $N$ in \eqref{eqqN} when $A$ has only one eigenvalue. In fact, in that case, if
$\mm=(m_1,\ldots, m_w)$ is the Segre characteristic of $A$ then $N=\sum_{j=1}^n (2j-1)m_j$. Now, 
$w_i-w_{i+1}=\#\left\{j: m_j=i\right\}$; that is,  there are $w_m$ numbers in $\mm$
equal to $m_1=m$, $w_{m-1}-w_m$ equal to $m-1$, $w_{m-2}-w_{m-1}$ equal to $m-2$ , \ldots, $w_2-w_1$ equal
to $1$ (of course $w_j-w_{j+1}$ can be $0$ for some $j$). Hence, with the agreement $\sum_{i=p+1}^{p}:=0$ ($p\geq 0$), 
we get
\[
\begin{array}{rcl}
N&=&\sum\limits_{j=1}^n (2j-1)m_j=\sum\limits_{j=1}^{w_m}(2j-1)m+\sum\limits_{j=w_m+1}^{w_{m-1}}(2j-1)(m-1)\\
&&+\sum\limits_{j=w_{m-1}+1}^{w_{m-2}}(2j-1)(m-2)+
\ldots  +\sum\limits_{j=w_2+1}^{w_1}(2j-1)1\\
&=& w_m^2m+(w_{m-1}^2-w_m^2)(m-1)+ (w_{m-2}^2-w_{m-1}^2)(m-2)\\
&&+\ldots+(w_1^2-w_2^2)1
=w_m^2+w_{m-1}^2+w_{m-2}^2+\ldots+ w_1^2.
\end{array}
\]
}\hfill$\Box$
\end{rem}

%
\begin{example}\label{exJCreal0}{\rm
Assume that $A\in\RR^{12\times 12}$ has $\la_0\in\RR$ as its only eigenvalue and let
$\mm=(4,2,2,2,1,1)$ and $\ww=(6,4,1,1)$ be its Segre and Weyr characteristics, respectively.
Then $N=4+3\cdot 2+5\cdot 2+7\cdot 2+9\cdot 1+11\cdot 1=6^2+4^2+1+1=54$ 
and the matrices of $C_{W(\la_0)}$ have the following form (recall that $\tau_i=w_{m-i+1}$ and so $\tau_1=\tau_2=1$, 
$\tau_3=4$ and $\tau_4=6$ and note that $\tau_2-\tau_1=0$):
\begin{equation}\label{eqexXr}
\renewcommand\arraystretch{1.3}
Y=
\begin{blockarray}{cccccccc}
 \begin{smallarray}{c}\tau_1\\1\end{smallarray}& \begin{smallarray}{c}\tau_3-\tau_2\\3\end{smallarray}& 
 \begin{smallarray}{c}\tau_4-\tau_3\\2\end{smallarray}&
\begin{smallarray}{c}\tau_1\\1\end{smallarray}&\begin{smallarray}{c}\tau_3-\tau_2\\3\end{smallarray} &
\begin{smallarray}{c}\tau_1\\1\end{smallarray}& 
 \begin{smallarray}{c}\tau_1\\1\end{smallarray}\\
\begin{block}{[ccc|cc|c|c]c}
d_{11}^{(1)} & D_{13}^{(1)} &  D_{14}^{(1)} & d_{11}^{(2)} & D_{13}^{(2)} & d_{11}^{(3)}  & d_{11}^{(4)} &  
\mbox{\scriptsize $\tau_1=1$}\\
0                      & D_{33}^{(1)} & D_{34}^{(1)} &  0                      & D_{33}^{(2)}   & D_{31}^{(3)}& D_{31}^{(4)} & 
\mbox{\scriptsize $\tau_3-\tau_2=3$}\\
0                     &           0             & D_{44}^{(1)} & 0                       & D_{43}^{(2)} & 0                      &  D_{41}^{(4)}&
\mbox{\scriptsize $\tau_4-\tau_3=2$}\\
\BAhhline{- - - - - - - }
  0                    &          0              &         0              &d_{11}^{(1)} & D_{13}^{(1)}&d_{11}^{(2)} &  d_{11}^{(3)} &
  \mbox{\scriptsize $\tau_1=1$} \\
   0                   & 0                       &   0                     &     0                  &D_{33}^{(1)} &0     &D_{31}^{(3)}& 
   \mbox{\scriptsize $\tau_3-\tau_2=3$}\\
 \BAhhline{- - - - - - - }
     0                   & 0                       &   0                     &     0               &0                        &d_{11}^{(1)} & d_{11}^{(2)} & 
     \mbox{\scriptsize $\tau_1=1$} \\
    \BAhhline{- - - - - - - }  
   0                   & 0                       &   0                     &     0               &0                        &  0                        &d_{11}^{(1)}  &
   \mbox{\scriptsize $\tau_1=1$}. \\
\end{block}
\end{blockarray}
\end{equation}
}\hfill $\Box$
\end{example}

Let $\la_0,\overline{\la}_0\in  \CC\setminus \RR$ be eigenvalues of $A\in\RR^{n\times n}$. Let
$(m_1,\ldots, m_w)$  and $(w_1,\ldots, w_m)$ be the Segre and Weyr characteristics for both $\la_0$ and
$\overline{\la}_0$. Let $\whJ(\la_0,\overline{\la}_0)=\diag\left(\whJ_{1} (\la_0,\overline{\la}_0), \dots, 
\whJ_{w} (\la_0,\overline{\la}_0)\right)$ be the block associated to $\la_0$ and $\overline{\la}_0$ in the real
Jordan canonical form of $A$. A characterization of $C_{\whJ(\la_0,\overline{\la}_0)}$ can be found in
several publications (see, for example, \cite[Theorem 12.4.2]{GoLaRo86}, \cite[Theorem 5.6]{MiMoRo20}
or \cite[Section 3]{LaZa12}). As in the case when all eigenvalues are real,  we are interested in the centralizer of 
$\whW(\la_0,\overline{\la}_0)$, the block associated to the pair of eigenvalues $\la_0$ and
$\overline{\la}_0$ in the real Weyr canonical form of $A$. Since $Q^T\whJ(\la_0,\overline{\la}_0) Q=
\whW(\la_0,\overline{\la}_0)$ where $Q$ is the matrix
of \eqref{eq.defQC}, $X\in C_{\whJ(\la_0,\overline{\la}_0)}$ if and only if
$Q^TXQ\in C_{\whW(\la_0,\overline{\la}_0)}$. Using this property we get

\begin{lemma}\label{lemmazcomplex}
 Let $\la_0,\overline{\la}_0 \in \CC\setminus\RR$ and $\whW(\la_0,\overline{\la}_0)$ be the matrix of 
 \eqref{eq.blockWeyrcomplex} with  Weyr characteristic
 $(w_1,w_2,\ldots, w_m)$ for each eigenvalue $\la_0$ and $\overline{\la}_0 $. Then 
 $Y\in C_{\whW(\la_0,\overline{\la}_0)}$ if and only if $Y$ has the structure of \eqref{eq.QXQr} satisfying the
 properties \eqref{eq.QYQm},  for $1\leq i\leq  j \leq m-1$,
 \[
Y_{i+1\, j+1}=I_{2w_{i}, 2w_{i+1}}^T Y_{i\,j}I_{2w_{j}, 2w_{j+1}},
\]
 and for $1\leq i, j\leq m$ and $\max\{i-j+1, 1\}\leq k \leq m-j+1$
 \begin{equation}\label{eq.qxqdijc}
 D^{(j)}_{i,k}=\begin{bmatrix}T_{\alpha\,\beta}^{(j)}\end{bmatrix}_{\begin{smallarray}{l}
 \tau_{i-1}+1\leq \alpha\leq \tau_i\\\tau_{k-1}+1\leq \beta\leq\tau_k\end{smallarray}}
 \in \RR^{2(\tau_i-\tau_{i-1})\times 2(\tau_k-\tau_{k-1})}, 
\end{equation}
where $\tau_0=0$, $\tau_i=w_{m-i+1}$, $1\leq i\leq m$, and $T_{\alpha, \beta}^{(j)}=
\begin{bsmallmatrix}x_{\alpha\, \beta}^{(j)}&y_{\alpha\, \beta}^{(j)}\\-y_{\alpha\, \beta}^{(j)}&x_{\alpha\, \beta}^{(j)}
\end{bsmallmatrix}\in \RR^{2\times 2}$.
\end{lemma}

Note that if $B=\begin{bsmallmatrix}a&b\\-b& a\end{bsmallmatrix}\in\RR^{2\times 2}$ with $b\neq0$ then
$X\in\RR^{2\times 2}$ commutes with $B$ if and only if 
$X=\begin{bsmallmatrix}x_1&x_2\\-x_2& x_1\end{bsmallmatrix}$ for any $x_1,x_2\in\RR$.

On the other hand, because of \eqref{eqqN}, the dimension of $C_A$ is the same whether it is computed
on $\CC$ or $\RR$. Hence,  since $\la_0$ and $\overline{\la}_0$ have the same associated Segre and Weyr
characteristics, the number of parameters  in $Y$ when $\la_0$ and $\overline{\la}_0$ are the
only eigenvalues of $A$ is
\[
N=2\sum\limits_{j=1}^n (2j-1)m_j=2(w_m^2+w_{m-1}^2+\ldots+ w_1^2)=
2(\tau_1^2+\tau_2^2+\ldots+ \tau_m^2).
\]
\begin{example}\label{exJCinvcomplex0}{\rm
As in Example \ref{exJCreal0}, assume that $A\in\RR^{12\times 12}$ has $\la_0,\overline{\la}_0\in\CC\setminus\RR$ 
as its only eigenvalues and let
$\mm=(4,2,2,2,1,1)$ and $\ww=(6,4,1,1)$ be the Segre and Weyr characteristics, respectively, for both $\la_0$ and 
$\overline{\la}_0$.
Then $N=2(4+3\cdot 2+5\cdot 2+7\cdot 2+9\cdot 1+11\cdot 1)=2(6^2+4^2+1+1)=108$ 
and the matrices of $C_{W(\la_0)}$ have the following form :
\[
  Y=
\begin{blockarray}{cccccccc}
 \mbox{$\begin{smallarray}{c}2\tau_1\\2\end{smallarray}$} & 
 \mbox{$\begin{smallarray}{c}2(\tau_3-\tau_2)\\6\end{smallarray}$}& 
 \mbox{$\begin{smallarray}{c}2(\tau_4-\tau_3)\\4\end{smallarray}$}&
 \mbox{$\begin{smallarray}{c}2\tau_1\\2\end{smallarray}$}&
 \mbox{$\begin{smallarray}{c}2(\tau_3-\tau_2)\\6\end{smallarray}$}& 
  \mbox{$\begin{smallarray}{c}2\tau_1\\2\end{smallarray}$}&
   \mbox{$\begin{smallarray}{c}2\tau_1\\2\end{smallarray}$}\\
\begin{block}{[ccc|cc|c|c]c}
D_{11}^{(1)} & D_{13}^{(1)} &  D_{14}^{(1)} & D_{11}^{(2)} & D_{13}^{(2)} & D_{11}^{(3)}  & D_{11}^{(4)} &  
\mbox{\scriptsize $2\tau_1=2$}\\
0                      & D_{33}^{(1)} & D_{34}^{(1)} &  0                      & D_{33}^{(2)}   & D_{31}^{(3)}& D_{31}^{(4)} & 
\mbox{\scriptsize $2(\tau_3-\tau_2)=6$}\\
0                     &           0             & D_{44}^{(1)} & 0                       & D_{43}^{(2)} & 0                      &  D_{41}^{(4)}&
\mbox{\scriptsize $2(\tau_4-\tau_3)=4$}\\
\BAhhline{- - - - - - - }
  0                    &          0              &         0              &D_{11}^{(1)} & D_{13}^{(1)}&d_{11}^{(2)} &  D_{11}^{(3)} &
  \mbox{\scriptsize $2\tau_1=2$} \\
   0                   & 0                       &   0                     &     0                  &D_{33}^{(1)} &0     &D_{31}^{(3)}& 
   \mbox{\scriptsize $2(\tau_3-\tau_2)=6$}\\
 \BAhhline{- - - - - - - }
     0                   & 0                       &   0                     &     0               &0                        &D_{11}^{(1)} & D_{11}^{(2)} & 
     \mbox{\scriptsize $2\tau_1=2$} \\
    \BAhhline{- - - - - - - }  
   0                   & 0                       &   0                     &     0               &0                        &  0                        &D_{11}^{(1)}  &
   \mbox{\scriptsize $2\tau_1=2$},\\
\end{block}
\end{blockarray}
\]
with
\[
\begin{array}{l}
D_{11}^{(i)}=\begin{matrix}T_{11}^{(i)}\end{matrix}, 1\leq i\leq 4,\;\;
D_{13}^{(i)}=\begin{bmatrix}T_{12}^{(i)} &T_{13}^{(i)}&T_{14}^{(i)}\end{bmatrix}, i=1,3,\;\;
D_{14}^{(1)}=\begin{bmatrix}T_{15}^{(1)} &T_{16}^{(1)}\end{bmatrix} \\
D_{33}^{(i)}=\begin{bmatrix}T_{22}^{(i)} &T_{23}^{(i)}&T_{24}^{(i)}\\
T_{32}^{(i)} &T_{33}^{(i)}&T_{34}^{(i)}\\
T_{42}^{(i)} &T_{43}^{(i)}&T_{44}^{(i)}\end{bmatrix}, i=1,2,\;\;
D_{34}^{(1)}=\begin{bmatrix}T_{25}^{(1)} &T_{26}^{(1)}\\
T_{35}^{(1)} &T_{36}^{(1)}\\T_{45}^{(1)} &T_{66}^{(1)}\end{bmatrix},\;\;
D_{31}^{(i)}=\begin{bmatrix}T_{21}^{(i)}\\ T_{31}^{(i)}\\T_{41}^{(i)}\end{bmatrix}, i=3,4\\
D_{44}^{(1)}=\begin{bmatrix}T_{55}^{(1)} &T_{56}^{(1)}\\T_{65}^{(1)} &T_{66}^{(1)}\end{bmatrix},
D_{43}^{(2)}=\begin{bmatrix}T_{52}^{(2)} &T_{53}^{(2)}&T_{54}^{(2)}\\
T_{62}^{(2)} &T_{63}^{(2)}&T_{64}^{(2)}\end{bmatrix},
D_{41}^{(4)}=\begin{bmatrix}T_{51}^{(4)}\\ T_{61}^{(i4)}\end{bmatrix}
\end{array}
\]
and
\[
T_{\alpha, \beta}^{(i)}=\begin{bsmallmatrix}x_{\alpha, \beta}^{(i)}&y_{\alpha, \beta}^{(i)}\\
-y_{\alpha, \beta}^{(i)}&x_{\alpha, \beta}^{(i)}
\end{bsmallmatrix}\in \RR^{2\times 2}, \quad 1\leq \alpha, \beta \le 6,\quad 1\leq i \leq 4.
\]}\hfill$\Box$
\end{example}


\medskip

\subsection{Feedback Equivalence and the Pole Assignment Problem by State-Feedback}

Let $\Sigma^c$ be the  open subset of $\RR^{n\times n}\times \RR^{n\times m}$
formed by the controllable pairs of matrices. That is to say,
$$
\Sigma^c = \{(F, G)\in \RR^{n\times n}\times \RR^{n\times m}\;:\;
\rank\begin{bmatrix}G &FG\dots &F^{n-1}G\end{bmatrix}=n\}.$$
Let
$\G_c = \left\{
\begin{bsmallmatrix}
P & 0\\
R & Q
\end{bsmallmatrix} \;:\;  P\in \Gl(n), Q\in \Gl(m), R\in \RR^{m\times n} \right\}$ denote the feedback group.
Two pairs  of matrices $(F, G), (F', G')\in \Sigma^c$ are said to be {\em feedback equivalent} if
$
\begin{bsmallmatrix}
  F' &  G'
\end{bsmallmatrix}= P^{-1}\begin{bsmallmatrix}
F&G
\end{bsmallmatrix} \begin{bsmallmatrix}
P & 0\\
R & Q
\end{bsmallmatrix} 
$ with $\begin{bsmallmatrix}
  P & 0\\
R & Q
\end{bsmallmatrix}\in \G_c$.

It is well-known that the controllability indices form a
complete system of invariants for the feedback equivalence relation. However also \textit{the 
Brunovsky indices} form a complete system of invariants.  Both
indices are closely related. Let us briefly  recall their definitions. 

Assume that we are given a controllable system $(F,G)\in\Sigma^c$ and $\rank G=r$. For $i=1,\ldots n$ let 
$r_1+\ldots+r_i=\rank\begin{bmatrix} G & FG & \cdots & F^{i-1}G\end{bmatrix}$ (see \cite{Bru70}).
Then there is a positive integer $k$ such that $r=r_1\geq r_2\geq\cdots\geq r_k>0=r_{k+1}=\ldots=r_n$. 
The nonnegative integers $r_1\geq r_2\geq\cdots \geq r_n\geq 0$ are called the Brunovsky indices of $(F,G)$
(they were called $r$-numbers in \cite{ BaZa90}).  
As $r_1+r_2+\cdots+ r_n=n$, $\rr=(r_1, r_2, \ldots,  r_n)$ is a partition of $n$.
The controllability indices of $(F,G)$ are the components of its \textit{conjugate partition}.
That is to say, for $i=1,\ldots, m$, $k_i$ is the number of elements of
$\rr$ that are not smaller than $i$: $k_i=\#\{j: r_j\geq i\}$. Hence, bearing in mind that $r=r_1$,
 $k=k_1\geq k_2\geq \cdots\geq k_{r}>0
= k_{r+1}=\cdots=k_m$.


There are canonical representatives in each feedback equivalence class associated to either the controllability
indices or the Brunovsky indices. The so-called and well-known \textit{Brunovsky canonical form} is associated
to the controllability indices (see, for instance, \cite[Theorem 6.2.5]{GoLaRo86}):
Let $(F,G)\in \Sigma^c$ be a controllable pair with controllability indices 
$\kk: k_1\geq \dots \geq k_r>0=k_{r+1}=\dots =k_m$. 
Then $(F, G)$ is feedback equivalent to $(F_c, G_c)$, where
\[
  \label{eqFcGc}
\begin{array}{c}
F_c=\diag(J_{1}(0), \dots, J_{r}(0))\in \RR^{n\times n}, \quad G_c=\begin{bmatrix}
  G_1 & \nsc
\end{bmatrix}\in \RR^{n\times (r+(m-r))}, \\\\
J_i(0)=\begin{bsmallmatrix}
0 & I_{k_i-1}\\0 & 0\end{bsmallmatrix}\in \RR^{k_i\times k_i},\;
G_1=\begin{bsmallmatrix}
   E_1\\ \vdots \\E_r  
\end{bsmallmatrix}\in \RR^{n\times r}, \;E_i=\begin{bsmallmatrix}\nsc\\e_i^T 
\end{bsmallmatrix}\in \RR^{((k_i-1)+1)\times r},
\end{array}
\]
and $e_i$ is the $i$-th column of the identity matrix $I_r$, $1\leq i \leq r$.

Using the permutation matrix of \eqref{eq.defQ} with $m_i=k_i$, we get (compare with
\cite[Theorem 3.3]{ BaZa90} that it is sometimes called the \textit{dual
Brunovsky canonical form}):
\begin{equation}
\begin{array}{c}
  \label{eq.FpGp}
F_p=Q^TF_cQ=
\begin{bmatrix} 0 & I_{r_1,r_2} & 0&\cdots&0&0\\
0 & 0& I_{r_2,r_3} &\cdots&0&0\\
\vdots&\vdots&\vdots&\vdots&\ddots&\vdots\\
0 & 0& 0& \cdots&0&I_{r_{k-1},r_k} \\
0&0&0&\cdots&0&0
\end{bmatrix}, \\
G_{p}=Q^TG_c=\begin{bmatrix}
0 & 0 & \cdots & 0 &E_{r_1-r_{2}} & 0\\
0&0&\cdots & E_{r_{2}-r_{3}} & 0&0\\
\vdots&\vdots &\vdots& \vdots & \vdots&\vdots\\
0 &E_{r_{k-1}-r_k} & \cdots & 0 & 0 & 0\\
I_{r_k} & 0 & \cdots & 0 & 0 & 0
\end{bmatrix}
\end{array}
\end{equation}
where $I_{r_i,r_{i+1}}$ is defined in \eqref{eq.defIpq}
and
\[
E_{r_i-r_{i+1}}=\begin{bmatrix} 0\\I_{r_i-r_{i+1}}\end{bmatrix}\in\RR^{r_i\times (r_i-r_{i+1})}, i=1,2,\ldots, k-1.
\]
Note that  $F_p$ is the Weyr canonical form of $F_c$. The pair $(F_p,G_p)$  will be called the
\textit{permuted dual Brunovsky canonical form} or,  for short, the \textit{$p$-Brunovsky canonical form} of $(F,G)$.

\medskip
Recall that for a given controllable system $(F,G)\in \Sigma^c$ and a given sequence
of monic polynomials $\ualpha: \alpha_1(s)\mid \dots \mid \alpha_n(s)$ with $\sum_{i=1}^{n}\deg (\alpha_i(s))=n$,
we aim to parametrize the set of feedback matrices $K\in\RR^{m\times n}$ such that $F+GK$ has the polynomials
in $\ualpha$ as invariant polynomials; i.e., such that $F+GK\in \OO(\ualpha)$. Necessary and sufficient conditions for such a set not to be empty
were obtained in  \cite{Rosen70}  when $(F,G)$ is controllable and in \cite{Za89} in the general case.
We state the result for the controllable case.

\begin{proposition}\em{\cite[Ch. 5, Sec. 4]{Rosen70}, \cite[Theorem 2.6]{Za89}}\label{prop.necss}
  Let   $(F,G)\in \RR^{n \times n}\times  \RR^{n \times m}$ be a controllable pair  
and let $k_1\geq  \dots \geq k_m$ be its controllability indices. Let $\ualpha: \alpha_1(s)\mid \dots \mid \alpha_n(s)$ 
be monic polynomials. 
There exists $K\in \RR^{m\times n}$ such that
 $ F+GK \in \OO(\ualpha)$ if and only if (see \eqref{eq.maj1})
\begin{equation}\label{eqnecss}
 (k_1,k_2,\ldots, k_m)\prec (\deg(\alpha_n(s)),\deg(\alpha_{n-1}(s)),\ldots,\deg(\alpha_1(s))).
\end{equation}
\end{proposition}


\begin{rem}\label{rem.necssweyr}
Assume that the prime factorization of $\alpha_{n-i+1}(s)$ is given by \eqref{eq.primefacalfa} and 
$\alpha_1(s)=\cdots =\alpha_{h}(s)=1\neq \alpha_{h+1}(s)$.
\begin{itemize}
\item If $\mm_i=(m_{i1}, \dots, m_{iw_i})$ is the Segre characteristic for $\la_i$, $1\leq i\leq t$, then
\[
(\deg(\alpha_n(s)),\deg(\alpha_{n-1}(s)),\ldots,\deg(\alpha_{h+1}(s)))=\mm_1+\mm_2+\cdots +\mm_{p+2q}.
\]
Therefore, \eqref{eqnecss} is equivalent to
\begin{equation}\label{eq.necsssegre}
 (k_1,k_2,\ldots, k_m)\prec (m_{11},\ldots, m_{1w_1})+\cdots+(m_{p+2q\, 1},\ldots, m_{p+2q\,w_{p+2q}})
\end{equation}
\item For $i=1,\ldots, t$, let $\ww_i=(w_{i1},\ldots, m_{iw_i})$ be the Weyr characteristic of $\la_i$ and
let $\rr=(r_1,\ldots, r_n)$ be
the Brunovsky indices of $(F,G)$. Then,  by definition, $\ww_i$ and $\rr$ are the conjugate partitions of $\mm_i$,
$1\leq i\leq p+2q$, and $\kk$, respectively. It follows from \eqref {eq.necsssegre} and Proposition \ref{prop.uncmajc}
that condition \eqref{eqnecss} is equivalent to
\begin{equation}\label{eq.necssweyr}
(w_{11},\ldots,w_{1m_1})\cup\cdots\cup (w_{p+2q\, 1},\ldots, w_{p+2q\,m_{p+2q}})\prec (r_1, r_2,\ldots, r_n).
\end{equation}
Note that $\ww_1\cup\ww_2\cup\cdots\cup\ww_t=(\deg(\alpha_n(s)),\deg(\alpha_{n-1}(s)),\ldots,\deg(\alpha_1(s)))^\ast$,
the conjugate partition of $(\deg(\alpha_n(s)),\deg(\alpha_{n-1}(s)),\ldots,\deg(\alpha_1(s)))$. 
\end{itemize}\hfill $\Box$
\end{rem}

\section{Geometric structure of $\HH_{(F,G)}$}
\label{secgeomstrh}

Let $(F,G)\in\Sigma^c$ be a controllable system with $k_1\geq k_2\geq\cdots\geq k_r>0=k_{r+1}=\cdots= k_m$
as controllability indices and let $\ualpha: \alpha_1(s)\mid \dots \mid \alpha_n(s)$ be monic polynomials.
Let $\HH_{(F,G)}$ is the set of \eqref{eqH} and assume that it is not empty; i.e., condition  \eqref{eqnecss} 
holds true. In this section we will prove that $\HH_{(F,G)}$ is a submanifold of $\RR^{m\times n}$ whose
dimension is $\dim \HH_{(F,G)}=nm-N$, where $N$ is given in \eqref{eqqN}.


\begin{lemma}\label{lemmatang}
Let  $\ualpha: \alpha_1(s)\mid \dots \mid \alpha_n(s)$ be monic polynomials such that
$\sum_{i=1}^{n}\deg(\alpha_{n})=n$ and let $A\in \OO(\ualpha)$. Then the tangent space of $\OO(\ualpha)$
at $A$ is
$$
T_A\OO(\ualpha)=\{[A,X]\; :\; X\in  \RR^{n\times n}\},
$$
where $[A,X]=AX-XA$ is the commutator of $A$ and $X$.
\end{lemma}
{\bf Proof.}  
Let $\gamma_A: \Gl(n)\longrightarrow \RR^{n\times n}$ be the map defined by
$\gamma_A(P)=P^{-1}AP$. It follows from the proof of Theorem 9.16 in \cite{Lee03} 
(see also \cite[Proposition 3.2]{BaPu19}) that
$T_A\OO(\ualpha)=\Ima d\gamma_{A,I_n}$. For $X \in \RR^{n\times n}$,
$$ \gamma_A(I_n+ \epsilon X)=(I_n+ \epsilon X)^{-1}A(I_n+ \epsilon X)=(I_n- \epsilon X+ \epsilon^2 X^2- \dots )A(I_n+ \epsilon X).$$
Then
$$\gamma_A(I_n+ \epsilon X)-\gamma_A(I_n)=\epsilon(A X -XA)+\epsilon^2P(\epsilon)$$
where $P(\epsilon)$ is a polynomial matrix whose coefficients depend on $A$ and $X$. Therefore
$d\gamma_{A,I_n}(X)=[A,X]$.
\hfill $\Box$

\medskip
In the proof of the following theorem we will use the Frobenius  inner product in $\RR^{n\times n}$:
if $A,B\in\RR^{n\times n}$, $<A,B>=\tr(A^TB)$, where $\tr$ stands for trace.

\begin{theorem}\label{theovd}
Let $(F,G)\in \RR^{n \times n}\times  \RR^{n \times m}$ be a controllabe pair   with controllability indices
 $k_1\geq  \dots \geq k_m$
 and  let $\ualpha: \alpha_1(s)\mid \dots \mid \alpha_n(s)$  be
monic polynomials satisfying (\ref{eqnecss}).  
Then  the set $\HH_{(F,G)}$ defined in (\ref{eqH}) is a submanifold of $\RR^{m\times n}$ and
$\dim \HH_{(F,G)}=nm-N$, where $N$ is given in (\ref{eqqN}).
\end{theorem}
{\bf Proof.}
Let $\varphi: \RR^{m\times n}\longrightarrow \RR^{n\times n}$ be the differentiable map defined by
$\varphi(K)=F+GK$. Then
$\varphi^{-1}(\OO(\ualpha))=\HH_{(F,G)}$ and
$d\varphi_K(V)=GV$ for all $V\in\RR^{m\times n}$.

If we prove that $\varphi$ is transversal to $\OO(\ualpha)$ then  (\cite[p. 28]{GuPo74}) 
$\varphi^{-1}(\OO(\ualpha))=\HH_{(F,G)}$ would be a submanifold of $\RR^{m\times n}$
of dimension $\dim \HH_{(F,G)}=mn-N$, as desired.

We take $K\in \varphi^{-1}(\OO(\ualpha))$ and we are to prove that
$
\Ima d\varphi_k+T_{\varphi(K)}\OO(\ualpha)=T_{\varphi(K)}\RR^{n\times n}.
$
Equivalently, bearing in mind that $T_{\varphi(K)}\RR^{n\times n}=\RR^{n\times n}$,
$$
(\Ima d\varphi_k)^\perp \cap (T_{\varphi(K)}\OO(\ualpha))^\perp =\{0\}.
$$
Let $U\in \RR^{n\times n}$.   On one hand, by Lemma \ref{lemmatang},
$U\in (T_{\varphi(K)}\OO(\ualpha))^\perp$ if and only if 
$0=<U,[F+GK,X]>
=\tr(U^T(F+GK)X-U^TX(F+GK))=\tr(U^T(F+GK)X-(F+GK)U^TX)=
\tr((U^T(F+GK)-(F+GK)U^T)X)$ for all $X\in \RR^{n\times n}$, i.e., if and only if 
$[F+GK,U^T]=0$. On the other hand,  $U\in (\Ima d\varphi_k)^\perp$ if and only if
$\tr(U^TGV)=0$ for all $V\in \RR^{m\times n}$, i.e., if and only if $U^TG=0$.

Therefore, if  $U\in (\Ima d\varphi_k)^\perp \cap (T_{\varphi(K)}\OO(\ualpha))^\perp$ then
$$
\begin{array}{rl}
0=&\begin{bmatrix}U^TG& (F+GK)U^TG& \dots &(F+GK)^{n-1}U^TG\end{bmatrix}\\=&
\begin{bmatrix}U^TG& U^T(F+GK)G& \dots &U^T(F+GK)^{n-1}G\end{bmatrix}\\=&
U^T\begin{bmatrix}G& (F+GK)G& \dots &(F+GK)^{n-1}G\end{bmatrix}.\end{array}
$$
Taking into account that $(F+GK, G)$ is controllable, the matrix
$$\begin{bmatrix}G& (F+GK)G& \dots &(F+GK)^{n-1}G\end{bmatrix}$$
is right invertible and so $U = 0$ as desired.
\hfill $\Box$

\medskip
The next proposition shows that in order to study the geometry of the set $\HH_{(F, G)},$ the pair $(F,G)$  can
be replaced by any other pair in its orbit of feedback equivalence.

\begin{proposition}\label{propfeedback}
Let $(F,G), (F',G') \in \RR^{n \times n}\times  \RR^{n \times m}$ be controllabe pairs.
If $(F,G), (F', G') $ are feedback equivalent then
$\HH_{(F,G)}$ and  $\HH_{(F',G')}$ are diffeomorphic.
\end{proposition}
{\bf Proof.}
There exist $P\in \Gl(n)$, $Q\in \Gl(m)$ and    $R\in \RR^{m \times n}$ such that $P^{-1}FP+P^{-1}GR=F'$ and $P^{-1}GQ=G'$.

Let $K \in \RR^{m\times n}$. Then  $F'+G'Q^{-1}(KP-R)=P^{-1}(F+GK)P$. Therefore, the map
$\psi: \HH_{(F,G)}\longrightarrow \HH_{(F',G')}$ defined by $\psi(K)=Q^{-1}(KP-R)$ is well
defined and bijective. It is easily seen that it is a diffeomorphism.
\hfill $\Box$

\section{The manifold $\PP_{(A; \rr)}/\wC_A$}
\label{secmanifpc}

In order to obtain a parameterization of the manifold $\HH_{(F,G)}$ defined in Section \ref{secgeomstrh},
we will prove that it is diffeomorphic  to an orbit space by the action of a Lie group.
%
We are led by the following idea taken from \cite[Section 2.2]{Antou83} (see also 
\cite[Section 2.3]{BaPuPuPuMCSS10}): Assume that we are given a controllable system
$(F_c,G_c)\in \RR^{n \times n}\times  \RR^{n \times m}$  in Brunovsky canonical form 
with $\kk: k_1\geq \dots \geq k_{r}>0=k_{r+1}=\dots =k_m$ as controllability indices. Then, by 
\cite[Theorem 2.15]{Antou83}, a subspace $\V\subset\RR^n$ of dimension $d$ is $(F_c,G_c)-invariant\, 
(i.e.; F_c\V\subset \V+\Ima\;G_c$)
if and only if there exists a pair of matrices $(\wh H,\wh F)\in \RR^{r \times d}\times  \RR^{d \times d}$ such that
 $\V=\Ima O_\pi(\wh H,\wh F)$ where $\wh H=\begin{bmatrix} h_1^T &\cdots& h_r^T\end{bmatrix}^T$ and
 \[ 
 O_\pi(\wh H, \wh F)=\begin{bmatrix}O_1^T &\cdots & O_r^T\end{bmatrix}^T,\;
 O_i=\begin{bmatrix}h_i^T& \wh F^Th_i^T\, \cdots\, (\wh F^T)^{k_i-1}h_i^T\end{bmatrix}^T,\, 1\leq i \leq r.
\]

Using Antoulas' notation,  $O_\pi(\wh H,\wh F)$ is a \textit{permuted and truncated observability matrix}
of $(\wh H,\wh F)$.  Note that $O_\pi(\wh H,\wh F)\in\PTO(\wh F,\kk)$ (cf.  \eqref{eq.defPTO}). It is then shown
(\cite[Corollary 2.18]{Antou83}) that $\wh F=(F_c+G_cK)_{\mid \V}$, i.e., $\wh F$ is the restriction of $F_c+G_cK$
to $\V$ for some state-feedback matrix $K$. In other words; if
$R=\begin{bmatrix}O_\pi(\wh H,\wh F) & X\end{bmatrix}\in \Gl(n)$ then
$(F_c+G_cK)R=O_\pi(\wh H,\wh F)\wh F$ for some feedback transformation $K$. In particular, if $\V=\RR^n$
then there is a pair  of matrices $(\wh H,\wh F)\in \RR^{r \times n}\times  \RR^{n \times n}$ such that
$O_\pi(\wh H,\wh F)$ is invertible and
 \begin{equation}\label{eq.antoulas}
 O_\pi(\wh H,\wh F)\wh F=(F_c+G_cK)O_\pi(\wh H,\wh F).
 \end{equation}
This result establishes a close relationship between the Antoulas' permuted 
and truncated observability matrices with fixed state matrix $A\in\OO(\ualpha)$ and the set $\HH_{(F_c,G_c)}$
and, by Proposition \ref{propfeedback}, with the set $\HH_{(F,G)}$ provided that $(F,G)$ and $(F_c,G_c)$ are
feedback equivalent. As Antoulas himself remarks if $(F,G)=(T^{-1}F_c T, T^{-1}G_c)$ for some invertible matrix
$T$, then a subspace is $(F,G)$-invariant if and only if it is spanned by $T^{-1}O_\pi(\wh H, \wh F)$.
In order to simplify the computations in Section \ref{secreducedform},
it is most convenient for us to work with some matrices whose rows are obtained by permuting in a precise way
the rows of Antoulas' permuted  and truncated observability matrices $O_\pi(\wh H, \wh F)$.  Specifically
if $Q$ is the permutation matrix of \eqref{eq.defQ} then $(F_p,G_p)=(Q^TF_cQ,Q^T G_c)$ is the
$p$-Brunovsky canonical form of \eqref{eq.FpGp} and so
 a subspace is $(F_p,G_p)$-invariant if and only if it is spanned by $P=Q^TO_\pi(\wh H, \wh F)$.
 A direct computation shows that $P$
has the following form: If  $\rr=(r_1, r_2,\ldots, r_k)$  ($r_{k+1}=0$) is the
conjugate partition of $\kk=(k_1, k_2,\ldots, k_r)$ then
\[
\begin{array}{l}
P=\begin{bmatrix}P_1\\P_2\\\vdots\\\\P_{k}\end{bmatrix}, \quad
P_1=\begin{bmatrix}P_{11}\\P_{12}\\\vdots\\P_{1k}\end{bmatrix}\in\RR^{r_1\times d}, 
P_{1i}=\begin{bmatrix}h_{r_{k-i+2}+1}\\h_{r_{k-i+2}+2}\\\vdots\\h_{r_{k-i+1}}\end{bmatrix}\in
\begin{array}[t]{l}\RR^{(r_{k-i+1}-r_{k-i+2})\times d},\\
1\leq i\leq k,\end{array}
\end{array}
\]
and for $i=1,\ldots, k-1$
\[
P_{i+1}=I_{r_{i}, r_{i+1}}^TP_{i}\wh F  =I_{r_{1}, r_{i+1}}^TP_{1}\wh F^{i}=\begin{bmatrix}
P_{11}\wh F^{i}\\P_{12}\wh F^{i}\\\vdots\\\\P_{1k-i}\wh F^{i}\end{bmatrix}\in\RR^{r_{i+1}\times d},
\]
where $I_{p,q}$ is the matrix of \eqref{eq.defIpq}.

Note that $P$ is a truncated observability matrix of $(\wh H, \wh F)$ but it is obtained from that matrix without
permuting its rows.

%
%
%
%
 
We define formally the set of matrices $P$ introduced above.
Given an arbitrary matrix $A\in \RR^{d \times d}$ and arbitrary
nonnegative integers $\rr: r=r_1\geq r_{2}\geq \dots \geq r_{k}>0=r_{k+1}=\cdots=r_d$ 
such that $n=\sum_{i=1}^d r_i\geq d$, we define 
\[
\PP_{(A;\rr)}:=
\left\{
P=\begin{bsmallmatrix}P_1\\P_2\\\vdots\\\\P_{k}\end{bsmallmatrix},\; P_i=\begin{bsmallmatrix}
P_{11}A^{i-1}\\P_{12}A^{i-1}\\\vdots\\\\P_{1k-i+1}A^{i-1}\end{bsmallmatrix}\in\RR^{r_i\times d},
i=1,\ldots, k, \rank P=d\right\}.
\]
Note that if $\PP_{(A;\rr)}\neq \emptyset$ then it is an open set of a linear subspace or dimension $rd$. Hence
$\PP_{(A;\rr)}$ is a linear manifold of dimension $rd$.

\begin{rem}\label{remPArnoempty}
\begin{itemize}
\item[(i)] We have seen that if $P=\begin{bsmallmatrix}P_1\\P_2\\\vdots\\\\P_k\end{bsmallmatrix}\in\PP(A,\rr)$ 
with $P_1=\begin{bsmallmatrix}p_1\\p_2\\\vdots\\\\p_r\end{bsmallmatrix}$ then $P_{1i}=\begin{bsmallmatrix}
p_{r_{k-i+2}+1}\\p_{r_{k-i+2}+2}\\\vdots\\\\p_{r_{k-i+1}}\end{bsmallmatrix}$, $i=1,\ldots, k$. It follows from this and
$r_{i}-r_{i+1}=\#\{j: k_j= i\}$ that
\begin{equation}\label{eq.PApA}
\begin{bsmallmatrix}P_{11}A^{k_1-1}\\P_{12}A^{k_1-2}\\\vdots\\\\P_{1\,k-1}A\\P_{1k}\end{bsmallmatrix}=
\begin{bsmallmatrix}p_1A^{k_1-1}\\p_2 A^{k_2-1}\\\vdots\\\\p_{r-1}A^{k_{r-1}-1}\\p_rA^{k_r-1}\end{bsmallmatrix}.
\end{equation}

\item[(ii)]It is worth-noticing that $\PP_{(A;\rr)}$ can be empty for some matrices $A$ and some sequences $\rr$.
 For example, if $A=I_4$ and $\rr=(2,2)$ then
 $\rank \begin{bsmallmatrix} p_1\\p_2\\p_1A\\p_2A \end{bsmallmatrix}<4$ for all vectors
 $p_1,p_2\in\RR^{1\times 4}$.

 \end{itemize}
\hfill $\Box$ 
\end{rem}
The following proposition provides a necessary and sufficient condition for $\PP_{(A;\rr)}\neq \emptyset$. As it is lengthy
and does not significantly contribute to the aim of the paper, its proof is deferred to the Appendix.

\begin{proposition}\label{prop.parnemp}
With the above notation,  if $\alpha_1(s)\mid\cdots\mid\alpha_d(s)$ are the invariant polynomials of $A$ and
$(w_1, \dots, w_d)$ is the conjugate partition of $(\deg(\alpha_d(s)), $ $\dots, \deg(\alpha_1(s)))$ then
$\PP_{(A;\rr)}\neq \emptyset$ if and only if each of the following equivalent conditions holds:
\begin{equation}\label{eq.kalpha}
\sum_{j=i+1}^r k_j\geq \sum_{j=1}^{d-i}\deg(\alpha_j),\quad i\geq 1
\end{equation}
\begin{equation}\label{eq.wr}
  \sum_{j=1}^i w_j\leq \sum_{j=1}^ir_j,\quad 1\leq i\leq d.
 \end{equation}
\end{proposition}

\begin{rem}\label{rem.deqn}
It should be noted that when $n=\sum_{i=1}^d r_i=d$ then, since $\sum_{i=1}^d \deg (\alpha_i(s))=d$,
 $\sum_{j=1}^d w_j=\sum_{j=1}^dr_j$. This and \eqref{eq.wr} implies $(w_1, \dots, w_d)\prec
(r_1,\ldots, r_d)$. Taking into account the second item of Remark \ref {rem.necssweyr}, if
$\sum_{i=1}^d r_i=d$ then $\HH_{(F,G)}\neq\emptyset$ if and only if $\PP_{(A;\rr)}\neq\emptyset$. In addition, in this case,
all matrices in $\PP_{(A;\rr)}$ are square and invertible.
\end{rem}

Our goal in this section is to show that, for any $A\in\RR^{d\times d}$, $\PP_{(A; \rr)}/\wC_A$ is a differentiable manifold 
and  that  there is a diffeomorphism between $\HH_{(F,G)}$ and $\PP_{(A; \rr)}/\wC_A$ when $d=n$ and
$A\in\OO(\ualpha)$.

\begin{proposition}\label{propfreeproper}
The action  $\sigma: \wC_A\times \PP_{(A;\rr)}\longrightarrow \PP_{(A;\rr)}$  of $\wC_A$ on 
$\PP_{(A;\rr)}$ defined by $\sigma(X, P)=PX$ is free and proper.
\end{proposition}
{\bf Proof.}
If $P\in \PP_{(A;\rr)}$ then $P$ is left invertible and thus $\sigma$ is free. 

Let $\{P_i\}$ be a convergent sequence in $\PP_{(A;\rr)}$ and $\{X_i\}$ a sequence in $\wC_A$ such that
$\{P_iX_i\}$ converges. Then  $\{(P_i^TP_i)^{-1}P_i^TP_iX_i\}=\{X_i\}$ converges. By \cite[Proposition 9.13]{Lee03},
the action $\sigma$ is proper.
\hfill $\Box$

As a consequence we can apply the quotient manifold theorem (see, for example, \cite[Theorem 9.16]{Lee03}).
\begin{corollary}\label{corpcman}
The space of orbits $\PP_{(A; \rr)}/\wC_A$ is a differentiable manifold,
 the natural projection $\pi:\PP_{(A;\rr)}\longrightarrow \PP_{(A;\rr)}/\wC_A$  is a submersion
 and $\dim \PP_{(A;\rr)}/\wC_A=rd-\dim \wC_A$.
\end{corollary}

In what follows $\ualpha: \alpha_1(s)\mid\cdots\mid\alpha_n(s)$ will be assumed to be  monic polynomials
satisfying $\sum_{i=1}^n\deg(\alpha_i(s))=n$ and  $(F, G) \in\Sigma^c$ a given controllable pair with
controllability indices $\kk: k_1\geq \dots \geq k_{r}>0=k_{r+1}=\dots =k_m$ satisfying \eqref{eqnecss}
and Brunovsky indices $r_1\geq \cdots r_k>0=r_{k+1}=\cdots=r_n$ ($k=k_1$ and $r=r_1$).
We  aim to obtain a parameterization of $\HH_{(F, G)}$. This will be achieved in
Section \ref{secparameterization} throughout a diffeomorphism between $\HH_{(F, G)}$ and 
$\PP_{(A;\rr)}/\wC_A$ where $A\in \RR^{n\times n}$ is a matrix with
$\alpha_1(s)\mid\cdots\mid\alpha_n(s)$ as invariant polynomials. That $\HH_{(F, G)}$ and
$\PP_{(A;\rr)}/\wC_A$ are diffeomorphic is proved in Theorem \ref{teodifeo} below.

Since $\ualpha$ and $\kk$ satisfy \eqref{eqnecss}, $\HH_{(F, G)}\neq \emptyset$ and, by Remark \ref{rem.deqn},
$\PP_{(A;\rr)}\neq\emptyset$ and $\PP_{(A;\rr)}\subset \Gl(n)$. Also, it follows from $k_r>0=k_{r+1}$ that $\rank G=r$.

\begin{rem}\label{remrankr}
By Proposition \ref{propfeedback}, we can assume that $(F, G)=(F_p, G_p)$ where $(F_p, G_p)$
is the $p$-Brunovsky canonical form given in (\ref{eq.FpGp}). 
Let 
$G=\begin{bmatrix}G_1&O\end{bmatrix},$ with $ G_1\in \RR^{n\times r}$, $\rank G_1 =r$. 
If  $K= \begin{bmatrix}K_1\\K_2\end{bmatrix}\in \RR^{(r+(m-r))\times n}$, then
$F+GK= F+G_1K_1$  and therefore
$$
\HH_{(F,G)}=\HH_{(F, G_1)}\times \RR^{(m-r)\times n}.
$$
Thus, it is enough to obtain a parameterization of $\HH_{(F, G_1)}$.
\end{rem} 

The following lemma gives the counterpart of \eqref{eq.antoulas} when the matrices of $\PP_{(A;\rr)}$ are used.
\begin{lemma}\label{lem.isformula}
Let $A\in\RR^{n\times n}$ and let $(F,G)$ be in $p$-Brunovsky canonical form with
$G=\begin{bmatrix}G_1&\nsc\end{bmatrix}$, $G_1\in \RR^{n\times r}$, $\rank G_1 =r$.  Let
$k_1\geq k_2\geq\cdots\geq k_r>0$ and $r_1\geq r_2\geq\cdots\geq r_k>0$ be the nonzero controllability 
and Brunovsky indices of $(F,G)$. Then:
\begin{itemize}
\item[(i)] For each $P\in\PP_{(A;\rr)}$ the matrix $K=\begin{bsmallmatrix} p_1A^{k_1}\\\vdots\\p_rA^{k_r}
\end{bsmallmatrix}P^{-1}$, where $p_1$,\ldots, $p_r$ are the first $r$ rows of $P$, is in $\HH_{(F,G_1)}$.

\item[(ii)] For each $K\in\HH_{(F,G_1)}$ there is $P\in\PP_{(A;\rr)}$ such that $PA=(F+G_1K)P$ and
$KP=\begin{bsmallmatrix} p_1A^{k_1}\\\vdots\\p_rA^{k_r}
\end{bsmallmatrix}$ where $p_1$,\ldots, $p_r$ are the first $r$ rows of $P$.
\end{itemize}
\end{lemma}

{\bf Proof.}
Assume that $P\in\PP_{(A;\rr)}$. Then $P=\begin{bsmallmatrix}P_1\\P_2\\\vdots\\P_{k}\end{bsmallmatrix}$
with $P_i=\begin{bsmallmatrix}P_{i1}\\P_{i2}\\\vdots\\P_{i \,k-i+1}\end{bsmallmatrix}\in\RR^{r_i\times d}$
and $P_{ij}=P_{1j}A^{i-1}\in\RR^{(r_{k-j+1}-r_{k-j+2})\times d},\;\; i=1,\ldots, k,\; j=1,\ldots, k-i+1$.
Let $K=\begin{bsmallmatrix} p_1A^{k_1}\\\vdots\\p_rA^{k_r}
\end{bsmallmatrix}P^{-1}$ where $p_1$,\ldots, $p_r$ are the rows of $P_1$ (recall that $r_1=r$).  
We aim to show that $PA=(F+G_1K)P$. As $\rank P=n$, this implies that $F+G_1K=PAP^{-1}\in\OO(\ualpha)$ and so
$K\in\HH_{(F,G_1)}$. In fact, define $Z=KP$. Then 
$Z=\begin{bsmallmatrix}P_{11}A^{k_1}\\P_{12}A^{k_1-1}\\\vdots\\\\P_{1\,k-1}A^2\\P_{1k}A\end{bsmallmatrix}$.
Put $Z_i=P_{1i}A^{k_1-i+1}$, $i=1,\ldots, k=k_1$.
Now, $PA=\begin{bsmallmatrix}P_1A\\P_2A\\\vdots\\P_{k}A\end{bsmallmatrix}$
and it follows from \eqref{eq.FpGp} that
\begin{equation}\label{eq.FPGPZ}
\renewcommand\arraystretch{0.85}
FP=
\begin{blockarray}{cc}
\begin{block}{[c]c}
\text{\scriptsize $P_2$} & \text{\scriptsize $r_2$}\\
\text{\scriptsize $0$} & \text{\scriptsize $r_1-r_2$}\\\cline{1-1}
\text{\scriptsize $P_3$} &\text{\scriptsize $r_3$}\\
\text{\scriptsize $0$}&\text{\scriptsize $r_2-r_3$}\\\cline{1-1}
\vdots&\vdots\\\cline{1-1}
\text{\scriptsize $P_k$}&\text{\scriptsize $r_k$}\\
\text{\scriptsize $0$}& \text{\scriptsize $r_{k-1}-r_k$}\\\cline{1-1}
\text{\scriptsize $0$} &\text{\scriptsize $r_k$}\\
\end{block}
\end{blockarray}\quad\text{and}\quad G_1Z=\begin{blockarray}{cc}
\begin{block}{[c]c}
\text{\scriptsize $0$} & \text{\scriptsize $r_2$}\\
\text{\scriptsize $Z_k$}&\text{\scriptsize $r_1-r_2$}\\\cline{1-1}
\text{\scriptsize $0$} &\text{\scriptsize $r_3$}\\
\text{\scriptsize $Z_{k-1}$}&\text{\scriptsize $r_2-r_3$}\\\cline{1-1}
\vdots&\vdots\\\cline{1-1}
\text{\scriptsize $0$}& \text{\scriptsize $r_k$}\\
\text{\scriptsize $Z_2$}&\text{\scriptsize $ r_{k-1}-r_k$}\\\cline{1-1}
\text{\scriptsize $Z_1$} &\text{\scriptsize $r_k$}\\
\end{block}
\end{blockarray}.
\end{equation}
Bearing in mind that
\[
P_2=\begin{bmatrix}P_{11}A\\P_{12}A\\\vdots\\P_{1\,k-1}A\end{bmatrix},\;
P_3= \begin{bmatrix}P_{11}A^2\\P_{12}A^2\\\vdots\\P_{1\,k-2}A^2\end{bmatrix},\;\ldots,\;
P_k= P_{11} A^{k_1-1},
\]
and
\[
Z_k=P_{1k}A, \; Z_{k-1}=P_{1\,k-1}A^2,\;\ldots,\; Z_2=P_{12}A^{k_1-1},\;
Z_1=P_{11} A^{k_1},
\]
we get $PA=FP+G_1Z=(F+G_1K)P$, as claimed. Finally, by \eqref{eq.PApA},
$Z=\begin{bsmallmatrix}P_{11}A^{k_1}\\P_{12}A^{k_1-1}\\\vdots\\\\P_{1\,k-1}A^2\\P_{1k}A\end{bsmallmatrix}=
\begin{bsmallmatrix} p_1A^{k_1}\\p_2A^{k_2}\\\vdots\\p_{r-1}A^{k_{r-1}}\\p_rA^{k_r}\end{bsmallmatrix}$.

Conversely, if $K\in\HH_{(F,G_1)}$ then there is $P\in\Gl(n)$ such that $PA=(F+G_1K)P$. Split
$P=\begin{bsmallmatrix}P_1\\P_2\\\vdots\\P_{k}\end{bsmallmatrix}$
with $P_i=\begin{bsmallmatrix}P_{i1}\\P_{i2}\\\vdots\\P_{i \,k-i+1}\end{bsmallmatrix}\in\RR^{r_i\times d}$
and $P_{ij}\in\RR^{(r_{k-j+1}-r_{k-j+2})\times d},\;\; i=1,\ldots, k,\; j=1,\ldots, k-i+1$. Put $Z=KP=
\begin{bsmallmatrix}Z_1\\\vdots\\Z_k\end{bsmallmatrix}$ with $Z_i\in\RR^{(r_{k-i+1}-r_{k-i+2})\times n}$,
$i=1,\ldots, k$. By using that $PA=\begin{bsmallmatrix}P_1A\\P_2A\\\vdots\\P_{k}A\end{bsmallmatrix}$ and
\eqref{eq.FPGPZ} we get $P_{ij}=P_{1j}A^{i-1}$ for $ i=1,\ldots, k$, $ j=1,\ldots, k-i+1$, and
$Z_i=P_{1i}A^{k_1-i+1}$, $i=1,\ldots, k$. Therefore $P\in\PP_{(A;\rr)}$ and by \eqref{eq.PApA},
$KP=Z=\begin{bsmallmatrix} p_1A^{k_1}\\\vdots\\p_rA^{k_r}
\end{bsmallmatrix}$ where $p_1$,\ldots, $p_r$ are the rows of $P_1$.
\hfill$\Box$

We are ready to prove that $\HH_{(F,G_1)}$ and $\PP_{(A;\rr)}/\wC_A$ are diffeomorphic manifolds.

\begin{theorem}\label{teodifeo}
Let $(F, G)\in\Sigma^c$ be in $p$-Brunovsky canonical form, with $G=\begin{bmatrix}G_1&\nsc\end{bmatrix}$,
$G_1\in \RR^{n\times r}$, $\rank G_1 =r$ and $k_1\geq\cdots\geq k_r>0=k_{r+1}=\cdots=k_m$ as controllability
indices, and let $A\in\RR^{n\times n}$ be a matrix in $\OO(\ualpha)$.  Then the map
\begin{equation}\label{eq.defphi}
\begin{array}{rccc}
\phi: &
\PP_{(A;\rr)}/\wC_A   &\longrightarrow &   \HH_{(F,G_1)}\\
&\wt P& \mapsto & \begin{bsmallmatrix}p_1A^{k_1}\\p_{2}A^{k_{2}}\\\vdots \\p_{r}A^{k_r}\\
\end{bsmallmatrix}P^{-1},
\end{array}
\end{equation}
where $P\in\PP_{(A;\rr)}$ is any matrix in the orbit $\wt P$ and $p_1,\ldots, p_r$ are its first $r$ rows,
is a diffeomorphism. 
\end{theorem}
  {\bf Proof.} Let us see first that $\phi$ is well-defined. If $\wt P_1=\wt P_2$ then for any $P_1\in \wt P_1$ and
 $P_2\in\wt P_2$,  $P_1=P_2X$ for some $X\in\wC_A$. Then if $\phi(\wt P_i)=
\begin{bsmallmatrix}p_{i1}A^{k_1}\\p_{i2}A^{k_{2}}\\\vdots \\p_{ir}A^{r}\\
\end{bsmallmatrix}P_i^{-1}$, $i=1,2$, then
\[
\phi(\wt P_1)=\begin{bsmallmatrix}p_{21}XA^{k_1}\\\vdots\\p_{2r}XA^{k_r}\end{bsmallmatrix}X^{-1}P_2^{-1}=
\begin{bsmallmatrix}p_{21}A^{k_1}\\\vdots\\p_{2r}A^{k_r}\end{bsmallmatrix}XX^{-1}P_2^{-1}=\phi(\wt P_2).
\]
Next, let $P\in\wt P$. By item (i) of Lemma \ref{lem.isformula}, the matrix
$K=\begin{bsmallmatrix} p_1A^{k_1}\\\vdots\\p_rA^{k_r} \end{bsmallmatrix}P^{-1}$ is in $\HH_{(F,G_1)}$.
Thus $\phi$ is well-defined.

Conversely, if  $K\in\HH_{F,G_1)}$, by item (ii) of Lemma \ref{lem.isformula}, there is $P\in\PP(A,\kk)$ such that
$KP=\begin{bsmallmatrix} p_1A^{k_1}\\\vdots\\p_rA^{k_r} \end{bsmallmatrix}$.
Therefore $\phi$ is surjective.

Assume now  that $K_1,K_2\in\HH_{(F,G_1)}$ and $K_1=K_2$. By Lemma \ref{lem.isformula} there are
$P_1,P_2\in\PP(A,\kk)$ such that $P_1AP_1^{-1}=F+G_1K_1=F+G_1K_2=P_2AP_2^{-1}$. Then
$X=P_1^{-1}P_2\in\wt C_A$ and $P_2=P_1X$. Hence $\wt P_1=\wt P_2$ and $\phi$ is injective.

In order to prove that $\phi$ is a diffeomorphism, we introduce the set
\[
\wh \HH_{(F, G_1)}=\{F+G_1K\; : \; K\in \HH_{(F, G_1)} \}
\]
and the map $\wh \theta :\RR^{r\times n}\longrightarrow \RR^{n\times n}$ defined by $\wh\theta(K)=F+G_1K$.
This is a linear map whose differential has constant rank. Then it is an embedding
and $\wh\theta(\HH_{(F,G_1)})=\wh\HH_{(F, G_1)}$.
If $\wh\theta_{\HH_{(F,G_1)}}$ is the restriction of $\wh\theta$ to $\HH_{(F,G_1)}$, since $\HH_{(F,G_1)}$ is
a smooth manifold (Theorem \ref{theovd}), we can provide $\wh\HH_{(F,G_1)}$ with a smooth structure
for which $\wh\theta_{\HH_{(F,G_1)}}$ is a diffeomorphism.

Let $\psi:\PP_{(A;\rr)}/\wt C_A\longrightarrow  \wh \HH_{(F, G_1)}$ be the map defined by $\psi(\wt P)=PAP^{-1}$,
where $P\in\PP_{(A,\rr)}$ is any matrix in $\wt P$.
This is a well-defined and bijective map because $\psi=\wh\theta_{\HH_{(F,G_1)}}\circ\phi$. In fact,
\[
(\wh\theta_{\HH_{(F,G_1)}}\circ\phi)(\wt P)=\wh\theta_{\HH_{(F,G_1)}}\left(\begin{bsmallmatrix}p_1A^{k_1}\\
\vdots\\p_rA^{k_r}\end{bsmallmatrix}P^{-1}\right)=F+G_1K
\]
where  $P\in\PP_{(A,\rr)}$ is any representative of $\wt P$ and $K=\begin{bsmallmatrix}p_1A^{k_1}\\
\vdots\\p_rA^{k_r}\end{bsmallmatrix}P^{-1}$.  But $PA=(F+G_1K)P$ so that 
$(\wh\theta_{\HH_{(F,G_1)}}\circ\phi)(\wt P)=PAP^{-1}$. Taking into account that $\wh\theta_{\HH_{(F,G_1)}}$
is a diffeomorphism, we are going to prove that $\psi$ and $\psi^{-1}$ are differentiable. This proves that
$\phi$ is a diffeomorphism. 

We prove first that $\psi$ is differentiable. Let $\pi:\PP_{(A;\rr)}\longrightarrow\PP_{(A;\rr)}/\wC_A$ be the
natural projection, then $f:=\psi\circ \pi$ is the restriction to $\PP_{(A;\rr)}$ of  the differentiable map
$\wh f:\Gl(n)\longrightarrow\RR^{n\times n}$ defined by $\wh f(P)=PAP^{-1}$. That $f$ is differentiable follows
from the fact that $\wh f$ is differentiable. Since $f$ is differentiable and $\pi$ is,  by Corollary
\ref{corpcman}, a submersion, using \cite[Proposition 7.17]{Lee03}, we can conclude that
$\psi$ is differentiable.

Let us see now that $\psi^{-1}$ is also differentiable. First $f=\psi\circ \pi$ and so $f$ is surjective
because $\pi$ is surjective and $\psi$ is bijective. Now, for $P\in\PP_{(A,\rr)}$ and $U\in \RR^{n\times n}$,
a direct computation shows that $d\wh f_{P}(U)=UAP^{-1}-PAP^{-1}UP^{-1}$. Let $\wt\PP_{(A,\rr)}=T_P \PP_{(A,\rr)}$
be the tangent space of $\PP_{(A,\rr)}$ at $P$. Then
\[
 \wt\PP_{(A,\rr)}=\left\{\begin{bsmallmatrix}P_1\\P_2\\\vdots \\P_k\end{bsmallmatrix}\in \RR^{n\times n}\; : \; 
P_i=\begin{bsmallmatrix}P_{11}A^{i-1}\\\vdots\\P_{1\, k-i+1}A^{i-1}\end{bsmallmatrix}\in\RR^{r_i\times n},
1\leq i\leq k\right\}
\]
and $\dim  \wt\PP_{(A,\rr)}=rn$. In addition, $df_{P}=d\wh f_{P}\mid_{\wt\PP_{(A,\rr)}}$ and so
$\Ker df_{P}=\{U\in \wt\PP_{(A,\rr)} \; : \; P^{-1}UA=AP^{-1}U\}= \wt \PP_{(A,\rr)} \cap PC_{A}$. We claim that
$\dim \Ker df_{P}=\dim C_A$. In fact, the map
$\alpha :C_A\longrightarrow \Ker df_P$, defined by $\alpha(X)=PX$ is well-defined because $PX\in PC_{A}$
and $P_iX=\begin{bsmallmatrix}P_{11}A^{i-1}\\\vdots\\P_{1\, k-i+1}A^{i-1}\end{bsmallmatrix}X=
\begin{bsmallmatrix}P_{11}XA^{i-1}\\\vdots\\P_{1\, k-i+1}A^{i-1}X\end{bsmallmatrix}$ so that
$PX\in\wt\PP_{(A,\rr)}$. It is easy to see that $\alpha$ is bijective. Thus, $\alpha$ is an isomorphism of linear
spaces. As a conclusion we get $\dim \Ima df_{P}=rn-N=\dim \wh \HH_{(F,G_1)}$. Therefore $f$ is a
surjective submersion. 
Using again \cite[Proposition 7.17]{Lee03} with $\psi^{-1}\circ f=\pi$ we conclude that $\psi^{-1}$ is
differentiable, as claimed. 
\hfill $\Box$

\section{Parameterization of $\PP_{(A;\rr)}/\wC_A$}
\label{secreducedform}

Let $\ualpha: \alpha_1(s)\mid \dots \mid \alpha_n(s)$ be a sequence of monic polynomials such that
$\sum_{i=1}^n\deg (\alpha_i(s))=n$ and let $(F, G) \in  \RR^{n\times n}\times\RR^{n\times m}$ be a controllable
pair with controllability indices  $\kk: k=k_1\geq\dots \geq k_{r}>0=k_{r+1}=\dots =k_m$ satisfying (\ref{eqnecss}),  and 
Brunovsky indices $\rr: r=r_1\geq\dots \geq r_{k}>0=r_{k+1}=\dots =r_n$.
By  Remark \ref{remrankr} and Theorem \ref{teodifeo}, obtaining 
a parameterization of $\HH_{(F, G)}$ is equivalent to obtaining a parameterization of $\PP_{(A;\rr)}/\wC_A$ for any matrix
$A\in\OO(\ualpha)$. So, we can assume that $\alpha_{n-i+1}(s)$ factorizes as in \eqref{eq.primefacalfa}
and that $A$ is the associated real Weyr canonical form:


\begin{equation}\label{eq.AWR}
A=\diag(W_1, \dots, W_p, \whW_{p+1}, \dots, \whW_{p+q}),
\end{equation}
where $W_i=W(\lambda_i)$, $1\leq i \leq p$ and 
$\whW_{p+i}=\whW(\lambda_{p+i}, \overline{\lambda_{p+i}})$, $1\leq i\leq q$ are the matrices of \eqref{eq.defWk} and
\eqref{eq.defWR}, respectively.
Let  $s_i=\sum_{j=1}^{m_i}w_{i,j}$, $1\leq i \leq p+q$ and 
$P=\begin{bsmallmatrix}P_1\\\vdots \\P_k\end{bsmallmatrix}\in \PP_{(A;\rr)}$. Put 
$$
P_i=\begin{bmatrix}P_{i}^{(1)}&\cdots &P_{i}^{(p)}&P_{i}^{(p+1)}&\cdots &P_{i}^{(p+q)}\end{bmatrix}, \quad 1\leq i \leq k,
$$
and
$$
P^{(j)}=\begin{bmatrix}P_1^{(j)}\\\vdots \\P_k^{(j)}\end{bmatrix}, \quad 1\leq j \leq p+q.
$$
where  for $1\leq i \leq k$,
$$
P_i^{(j)}\in \RR^{r_i\times s_j}\, 1\leq i \leq p,\quad
P_i^{(j)}\in \RR^{r_i\times 2s_j}\,  p+1\leq j \leq p+q.
$$
For $1\leq i \leq k-1$,
\begin{equation}\label{eq.Pi1PiA}
\begin{array}{rcl}
P_{i+1}&=&I_{r_{i}, r_{i+1}}^TP_{i}A\\
&=&I_{r_{i}, r_{i+1}}^T
\begin{bmatrix}P_{i}^{(1)}W_1& \cdots &P_{i}^{(p)}W_p&P_{i}^{(p+1)}\whW_{p+1}&\cdots &P_{i}^{(p+q)}\whW_{p+q}
\end{bmatrix}.
\end{array}
\end{equation}
Hence
$$
P_{i+1}^{(j)}=
I_{r_{i}, r_{i+1}}^T
P_i^{(j)}W_i, \;
1\leq i \leq p\mbox{ and }
P_{i+1}^{(j)}=
I_{r_{i}, r_{i+1}}^T
P_i^{(j)}\whW_i,
 \; p+1\leq j \leq p+q.
$$
As $P^{(j)}$ are full column rank matrices, $1\leq j \leq p+q$, $P_i\in \PP_{(W_i;\rr)}$, $1\leq i \leq p$, and
$P_i\in \PP_{(\whW_i;\rr)}$, $p+1\leq i \leq p+q$. 
Let
$$
\PP=\PP_{(W_1;\rr)}\times \dots \times  \PP_{(W_p;\rr)}\times
\PP_{(\whW_{p+1};\rr)}\times \dots \times  \PP_{(\whW_{p+q};\rr)}.
$$
and note that $\PP_{(A;\rr)}$ can be identified with the subset of $\PP$ formed by their invertible matrices .
Thus we can think of $\PP_{(A;\rr)}$ as an open subset of $\PP$.

Recall that $X\in C_A$ if and only if
$X=\diag(X_1, \dots, X_{p}, \whX_{p+1}, \dots, \whX_{p+q})$ with
$X_i \in C_{W_i}$, $1\leq i \leq p$ and $\whX_i \in C_{\whW_i}$, $p+1\leq i \leq p+q$. Then,
$PX=\begin{bsmallmatrix}P^{(1)}X_1 &\dots &P^{(p)}X_p&P^{(p+1)}X_{p+1}&\dots &P^{(p+q)}X_{p+q}\end{bsmallmatrix}$.
As a consequence (see Corollary \ref{corpcman})  we can identify
$\PP_{(A;\rr)}/\widetilde C_A$ with an open subset of
$$\PP_{(W_1;\rr)}/\wC_{W_1}\times \dots \times  \PP_{(W_p;\rr)}/\wC_{W_p}\times
\PP_{(\whW_{p+1};\rr)}/\wC_{\whW_{p+1}}\times \dots \times  \PP_{(\whW_{p+q};\rr)}/\wC_{\whW_{p+q}},
$$
and we can parametrize $\PP_{(A;\rr)}/\wC_A$ from a parametrization of $\PP_{(W_i;\rr)}/\wC_{W_i}$, $1\leq i \leq p$
and  $\PP_{(\whW_i;\rr)}/\wC_{\whW_i}$, $p+1\leq i \leq p+q$.

Then, we aim to parameterize the manifold $\PP_{(A;\rr)}/\wC_A$
in two cases, when $A$ has only one real eigenvalue, $A=W(\lambda)$,
and when $A$ has two conjugate complex  eigenvalues,
$A=W(\lambda, \overline{\lambda})$.
To do this, we will obtain in both cases  a local reduced form for the 
equivalence relation associated to the action of $\wC_A$ on  $\PP_{(A;\rr)}$.
This equivalence relation  will be denoted by $\stackrel{\wC_A}{\sim}$.  
That is to say, given $P, \wh P\in \PP_{(A;\rr)}$, we will write  $P \stackrel{\wC_A}{\sim}\wh P$
if there exists $X\in  \wC_A$ such that 
$\wh P=PX$. If there is no risk of confusion we will  write $P {\sim} \wh P$ instead of $P \stackrel{\wC_A}{\sim} \wh P$.

\subsection{Reduced form when there is only a real eigenvalue}
\label{subsecredreal}
Let
$W=W(\lambda)$, $\lambda \in \RR$, with Weyr characteristic 
$(w_1, \dots , w_m)$ and assume that   $\PP_{(W;\rr)}\neq \emptyset$.
The procedure to bring a matrix $P\in \PP_{(W;\rr)}$ to a reduced form is based
on a sequence of elementary transformations defined by some subgroups
of $\wC_W$. It is worth-recalling at this point the structure of the matrices in $C_W$ (Lemma \ref{lemmazreal})
and that  $\tau_i=w_{m-i+1}$, $1\leq i\leq m$ and $\tau_0=0$ (see \eqref{eq.qxqdijr}).

\begin{definition}\label{defelementary}\,
\begin{enumerate}
\item 
Let $T_i \in \Gl(\tau_i-\tau_{i-1})$ , $1\leq i \leq m$ and
$Y_I=\diag(Y_{11}, \dots,  Y_{mm})$ with $Y_{ii}=\diag(T_1, \dots,  T_{m-i+1})$, $1\leq i \leq m$.
The matrices of this type will be called {\em elementary matrices of type I} and they form a
subgroup of $\wC_W$.
\item 
For $j=1$, $1\leq i <k\leq m$, and for
$2\leq j \leq m$, $1\leq k \leq m-j+1$ $1\leq i \leq k+j-1$ let
$Y_{II,i,k}^{(j)}$ be the matrix of \eqref{eq.QXQr}, with, perhaps, $D_{ik}^{(j)}\neq 0$,
$$
D_{ii}^{(1)}=I_{\tau_i-\tau_{i-1}}, \quad 1\leq i \leq m,
$$
and all the other blocks zero. This type of matrices will be called
{\em elementary matrices of type II} and they form a subgroup of $\wC_W$.
\end{enumerate}
\end{definition}

In addition to these elementary matrices we will use some auxiliary results.

\begin{proposition}\label{propauxiliarrw1}
Let $P=\begin{bsmallmatrix}P_1\\\vdots \\P_k\end{bsmallmatrix}\in\PP_{(W, \rr)}$ and partition
$P_i=\begin{bmatrix} P_{i1} & P_{i2}& \cdots & P_{im}\end{bmatrix}$ with $P_{ij}\in\RR^{r_i\times w_j}$,
$1\leq i\leq k$, $1\leq j\leq m$.
Then $\rank P_{11}=w_1$.
\end{proposition}
{\bf Proof.}
Since $\rank P=\sum_{j=1}^m w_j$,
$\rank \begin{bsmallmatrix}P_{1 j}\\\vdots \\P_{k j}\end{bsmallmatrix}=w_j$, $1 \leq j \leq m$. On the other
hand (see \eqref{eq.Pi1PiA}),  $P_{i+1}=I_{r_{i}, r_{i+1}}^TP_{i}W(\lambda)=I_{r_{1}, r_{i+1}}^TP_{1}W(\lambda)^{i}$, 
$1\leq i \leq k-1$. Thus for $i=1,\ldots, k-1$,
\[
P_{i+1\,1}=I_{r_{1}, r_{i+1}}^TP_{1}W(\la)^i\begin{bmatrix}I_{w_1}\\0\end {bmatrix}=
I_{r_{1}, r_{i+1}}^TP_{1}\begin{bmatrix}\la^iI_{w_1}\\0\end{bmatrix}=I_{r_{1}, r_{i+1}}^TP_{11}\lambda^{i}.
\] 
Then
$$\begin{bsmallmatrix}P_{11}\\\vdots \\P_{k 1}\end{bsmallmatrix}=
\begin{bsmallmatrix}P_{11}\\I_{r_{1}, r_{2}}^T\lambda P_{11}\\\vdots \\
I_{r_{1}, r_{k}}^T\lambda^{k-1} P_{11}
\end{bsmallmatrix}=\diag(I_{r_1}, I_{r_{1}, r_{2}}^T\lambda,\dots  
I_{r_{1}, r_{k}}^T\lambda^{k-1} )P_{11}.
$$
Therefore
$$w_1=\rank  \begin{bsmallmatrix}P_{11}\\\vdots \\P_{k1}\end{bsmallmatrix}=\rank P_{11}.$$
\hfill $\Box$

\bigskip

Recall (see Section  \ref{sepreliminaries}) that if $s$ and $p$ are positive integers ($0<s\leq p$) then
$Q_{s,p}:=\{(i_1, \dots, i_s)\;:\; 1\leq i_1<\dots <i_s\leq p\}$ and $Q_{0, p}:=\{\emptyset\}$. 

\begin{corollary}\label{corindicesrw}
  Let  $P\in \PP_{(W;\rr)}$. Then, for each $j=1,\ldots, m$, there is a sequence of $\tau_j$ indices
$\I_j\subseteq \{1, \dots, r\}$ such that
\begin{equation}\label{eqIsubsetr}
  \I_j\subseteq \I_{j+1}, \quad 1\leq j \leq m-1,
\end{equation}
\begin{equation}\label{eqIsetminus}
\I_{j}\setminus \I_{j-1}\in  Q_{\tau_j-\tau_{j-1},r}, \quad 1\leq j \leq m, \quad  (\I_{0}=\emptyset),
\end{equation}
\begin{equation}\label{eqIrankr}
   P(\I_j, 1:\tau_j)\in \Gl(\tau_j), \quad 1\leq j \leq m,
\end{equation}
\end{corollary}

\textbf{Proof}.  With the same notation as in Proposition \ref{propauxiliarrw1},  partition 
$P_{11}=\begin{bmatrix} P_{11}^{(1)} & P_{11}^{(2)}&\cdots &P_{11}^{(m)}\end{bmatrix}$
with $P_{11}^{(j)}\in\RR^{r_1\times (\tau_{j}-\tau_{j-1})}$, $1\leq j\leq m$.
By Proposition \ref{propauxiliarrw1}, $\rank P_{11}=w_1=\tau_m$.
Thus, $\rank\begin{bmatrix}P_{11}^{(1)} & \cdots &  P_{11}^{(j)}\end{bmatrix}=\tau_j$, $1\leq j \leq m$.
 
Since $\rank P_{11}^{(1)}=\tau_1$, in $P_{11}^{(1)}$  there must be $\tau_1$ linearly independent rows
$i_1<\dots< i_{\tau_1}$. Then $\I_1=(i_1,\ldots, i_{\tau_1})\in Q_{\tau_1, r}= Q_{\tau_1-\tau_0, r}$ and 
$P(\I_1, 1:\tau_1)=P_{11}^{(1)}(\I_1, :)\in \Gl(\tau_1)$. Now, $\rank\begin{bmatrix}P_{11}^{(1)} & P_{11}^{(2)}
\end{bmatrix}=\tau_2$. Thus, in $P_{11}^{(2)}$  there must be $\tau_2-\tau_1$
rows $i_{\tau_{1}+1}<i_{\tau_{1}+2}<\dots <i_{\tau_{2}}$ such that
the rows $i_1<\dots< i_{\tau_1}, \, i_{\tau_{1}+1}< \dots < i_{\tau_{2}}$ of
$\begin{bmatrix}P_{11}^{(1)} & P_{11}^{(2)}\end{bmatrix}$ are linearly independet. 
Put $\I_2=(i_1,\ldots, i_{\tau_2})$. Then $I_1\subseteq I_2$, $I_2\setminus I_1=
(i_{\tau_{1}+1}, \dots, i_{\tau_{2}})\in  Q_{\tau_2-\tau_{1},r}$, and
$P(\I_2, 1:\tau_2)= \begin{bmatrix}P_{11}^{(1)} & P_{11}^{(2)}\end{bmatrix}(\I_2, :)\in \Gl(\tau_2)$.
Continuing the process, we can obtain $m$ sequences of indices,
$\I_1, \dots, \I_m$ satisfying \eqref{eqIsubsetr}--
 \eqref{eqIrankr}. \hfill $\Box$


\begin{definition}\label{defadmissible}
Given $P\in \PP_{(W;\rr)}$,  let $\I_i$, $1\leq i \leq m$, be  sequences of indices satisfying
(\ref{eqIsubsetr})-(\ref{eqIrankr}). Then $\uI=(\I_1, \dots, \I_{m})$ will be called an
{\em admissible sequence of indices} for $P$.
  \end{definition}

\begin{proposition}\label{propsameadmissiblerw}
Let $P, \wh P\in \PP_{(W;\rr)}$ be matrices  such that $\wh P{\sim} P$ and let $\uI=(\I_1, \dots, \I_{m})$
be an admissible sequence of indices for $P$. Then $\uI$ is also an admissible sequence of indices for
$\wh P$.
\end{proposition}
{\bf Proof.}
First of all, since $\uI=(\I_1, \dots, \I_{m})$ is an admissible sequence of indices for $P$, it
satisfies \eqref{eqIsubsetr}  and \eqref{eqIsetminus}. So, it only remains to prove that 
$\wh P(\I_j,1:\tau_j)\in \Gl(\tau_j)$ for $1\leq j\leq m$.

Since $\wh P{\sim} P$, there exists  $Y\in \wC_{W}$ such that  $\wh P= PY$ and so $\wh P(\I_j, 1:\tau_j)=
 P(\I_j, :)Y(:, 1:\tau_j)$. By \eqref{eq.QYQm},
 \[
 Y(:, 1:\tau_j)=
\begin{bsmallmatrix}
         D_{11}^{(1)}&D_{12}^{(1)}&\dots& D_{1j}^{(1)}\\
 \nsc&D_{22}^{(1)}&\dots &D_{2j}^{(1)}\\
    \vdots &\vdots&\ddots & \vdots
    \\
    \nsc&\nsc&\dots &D_{jj}^{(1)}\\\hline
    \nsc&\nsc&\dots &\nsc\\
     \vdots&\vdots &\ddots&\vdots \\
     \nsc&\nsc&\dots &\nsc
        \end{bsmallmatrix}
 \]
 where $D_{ii}^{(1)}\in \Gl(\tau_i-\tau_{i-1})$, $1\leq i \leq j$. On the other hand,
 it follows from  \eqref{eqIrankr} that for $1\leq j \leq m$, $P(\I_j, 1:\tau_j)\in\Gl(\tau_j)$.
 Henceforth
$$
\wh P(\I_j, 1:\tau_j)=
 P(\I_j, :)Y(:, 1:\tau_j)=
P(\I_j, 1:\tau_j)
 \begin{bsmallmatrix}
         D_{11}^{(1)}&D_{12}^{(1)}&\dots& D_{1j}^{(1)}\\
  \nsc&D_{22}^{(1)}&\dots &D_{2,j}^{(1)}\\
    \vdots &\vdots&\ddots & \vdots
    \\
    \nsc&\nsc&\dots &D_{jj}^{(1)}
        \end{bsmallmatrix}\in \Gl(\tau_j)
$$
and the Proposition follows.
\hfill $\Box$

\bigskip
Let
\begin{equation}\label{eq.defAsW}
\A_{W}=\{\uI=(\I_1, \dots, \I_{m})\; : \;\I_j  \mbox{ satisfies \eqref{eqIsubsetr} and \eqref{eqIsetminus}, }
1\leq j \leq m\}.
\end{equation}
Given $\uI \in \A_W$,
we  denote by $\U_{\uI}$
 the open subset of 
 $\PP_{(W; \rr)}$ formed by their matrices 
with $\uI$ as an admissible sequence of indices.

We are ready to show a procedure to bring any $P\in\PP_{(W,\rr)}$ to a reduced form. First we illustrate this procedure
with an example.


\begin{example}\label{exWCreal}{\rm
Consider the invariant polynomials of the matrix  $A\in \RR^{12\times 12}$ of Example \ref{exJCreal0}. 
Let $W\in \RR^{12\times 12}$ be its Weyr canonical form. Then its
Weyr characteristic  is $(w_1, \dots, w_4)=(6, 4, 1, 1)$, $s=\sum_{i=1}^m w_i=12$, 
$\tau_1- \tau_0=1$ (recall that $\tau_0=0$), $\tau_2- \tau_1=0,$ $\tau_3- \tau_2=3,$ $\tau_4- \tau_3=2$
and $Y\in\wC_W$ is the matrix of \eqref{eqexXr} with $d_{11}^{(1)}\in \Gl(1)$, $D_{33}^{(1)}\in \Gl(3)$ and
$D_{44}^{(1)}\in\Gl(2)$. Let $\rr=(7, 4, 2, 1)$ and note that $\rr$ and $\ww$ satisfies the conditions of
Proposition \ref{prop.parnemp} so that $\PP_{(W,\rr)}\neq\emptyset$. Let 
\[
P=\begin{bsmallmatrix}
P_1\\P_2\\P_3\\P_4
\end{bsmallmatrix} \in  \PP_{(W;\rr)}
\]
Put $P_1=\begin{bsmallmatrix} p_1\\\vdots \\p_7\end{bsmallmatrix}$ and recall that $P_j=I_{r_1,r_j}^TP_1W^{j-1}$, 
$2\leq j\leq 4$. Thus
\[
P_2=\begin{bsmallmatrix}p_1W\\\vdots \\p_4W\end{bsmallmatrix},\quad P_3=\begin{bsmallmatrix}p_1W^2\\p_2W^2
\end{bsmallmatrix},\quad P_4=p_1W^3.
\]
Since $P_i$, $2\leq i\leq 4$, can be obtained from $P_1$, we only need to reduce $P_1$. Put
$P_1=\begin{bmatrix}P_{11}& P_{12}& P_{13}& P_{14}\end{bmatrix}$,
 $P_{1j}\in \RR^{r_1\times w_j}$, $1\leq j\leq 4$, and partition $P_{1j}$ as follows:
\[
P_{1j}=\begin{bmatrix}P_{11}^{(j)} & P_{12}^{(j)}&\cdots& P_{1\;5-j}^{(j)}\end{bmatrix}, \;
P_{1k}^{(j)}\in\RR^{r_1\times (\tau_k-\tau_{k-1})}, \; 1\leq j\leq 4,\,1\leq k\leq 5-j.
\]
Specifically,
\[
P_1=
\left[\begin{smallarray}{c|ccc|cc||c|ccc||c||c}
    p_{11}& p_{12}&p_{13}& p_{14}&p_{15}& p_{16}
& p_{17}& p_{18}&p_{19}& p_{1\,10}&p_{1\,11}&p_{1\,12}\\
   p_{21}& p_{22}&p_{23}& p_{24}&p_{25}& p_{26}
& p_{27}& p_{28}&p_{29}& p_{2\,10}&p_{2\,11}&p_{2\,12}\\
    p_{31}& p_{32}&p_{33}& p_{34}&p_{35}& p_{36}
& p_{37}& p_{38}&p_{39}& p_{3\,10}&p_{3\,11}&p_{3\,12}\\
 p_{41}& p_{42}&p_{43}& p_{44}&p_{45}& p_{46}
& p_{47}& p_{48}&p_{49}& p_{4\,10}&p_{4\,11} &p_{4\,12}\\
p_{51}& p_{52}&p_{53}& p_{54}&p_{55}& p_{56}
& p_{57}& p_{58}&p_{59}& p_{5\,10}&p_{5\,11} &p_{5\,12}\\
p_{61}& p_{62}&p_{63}& p_{64}&p_{65}& p_{66}
& p_{67}& p_{68}&p_{69}& p_{6\,10}&p_{6\,11}&p_{6\,12}\\
p_{71}& p_{72}&p_{73}& p_{74}&p_{75}& p_{76}
& p_{77}& p_{78}&p_{79}& p_{7\,10}&p_{7\,11} &p_{7\,12}\\
\end{smallarray}\right].
\]
Recall now that, by Proposition \ref{propauxiliarrw1}, $\rank P_{11}=w_1=6$. This means that
$\rank P_{11}^{(1)}=\tau_1=1$, $\rank \begin{bsmallmatrix}P_{11}^{(1)}&P_{12}^{(1)}\end{bsmallmatrix}=\tau_2=1$, 
$\rank  \begin{bsmallmatrix}P_{11}^{(1)}&P_{12}^{(1)}&P_{13}^{(1)}\end{bsmallmatrix}=\tau_3=4$
and $\rank \begin{bsmallmatrix}P_{11}^{(1)}&P_{12}^{(1)}&P_{13}^{(1)}&P_{14}^{(1)}\end{bsmallmatrix}=\tau_4=6$.
Then $P_{12}^{(1)}$ is an empty matrix. Let us assume that, for example,
\[
\begin{array}{l}
P_{11}^{(1)}(3,:)=p_{31}\neq 0, \qquad \det \begin{bmatrix}P_{11}^{(1)}&P_{13}^{(1)}\end{bmatrix}((3,1,4,7),:)
\neq 0,\\\\
\det\begin{bmatrix}P_{11}^{(1)}&P_{12}^{(1)}&P_{13}^{(1)}&P_{14}^{(1)}\end{bmatrix}((3,1,4,7,5,6),:)\neq 0. 
\end{array}
\]
Then $\uI=((3), (3,1,4,7), (3,1,4,7,5,6))$ is an admissible sequence of indices for $P$. For this admissible sequence of
indices we define
$Y_{11}^{a}=\diag(\begin{bmatrix}p_{31}^{-1}\end{bmatrix}, I_3, I_2)$, 
$Y_{22}^{a}=\diag(\begin{bmatrix}p_{31}^{-1}\end{bmatrix}, I_3)$, 
$Y_{33}^{a}=Y_{44}^{a}=\begin{bmatrix}p_{31}^{-1}\end{bmatrix}$ and
$Y_I^{a}=\diag(Y_{11}^{a}, Y_{22}^{a}, Y_{33}^{a}, Y_{44}^{a})$.
Then $Y_I^{a}$ is an elementary matrix of type I  (see Definition \ref{defelementary}) and
$$
P_1Y_I^{a}=
\left[\begin{smallarray}{cccccc|cccc|c|c}
    p_{11}^a& p_{12}&p_{13}& p_{14}&p_{15}& p_{16}
& p_{17}^a& p_{18}&p_{19}& p_{1\, 10}&p_{1\, 11}^a&p_{1\,12}^a\\
   p_{21}^a& p_{22}&p_{23}& p_{24}&p_{25}& p_{26}
& p_{27}^a& p_{28}&p_{29}& p_{2\,10}&p_{2\,11}^a &p_{2\,12}^a\\
    1& p_{32}&p_{33}& p_{34}&p_{35}& p_{36}
& p_{37}^a& p_{38}&p_{39}& p_{3\,10}&p_{3\,11}^a &p_{3\,12}^a\\
 p_{41}^a& p_{42}&p_{43}& p_{44}&p_{45}& p_{46}
& p_{47}^a& p_{48}&p_{49}& p_{4\,10}&p_{4\,11}^a &p_{4\,12}^a\\
p_{51}^a& p_{52}&p_{53}& p_{54}&p_{55}& p_{56}
& p_{57}^a& p_{58}&p_{59}& p_{5\,10}&p_{5\,11}^a &p_{5\,12}^a\\
p_{61}^a& p_{62}&p_{63}& p_{64}&p_{65}& p_{66}
& p_{67}^a& p_{68}&p_{69}& p_{6\,10}&p_{6\,11}^a &p_{6\,12}^a\\
p_{71}^a& p_{72}&p_{73}& p_{74}&p_{75}& p_{76}
& p_{77}^a& p_{78}&p_{79}& p_{7\,10}&p_{7\,11}^a &p_{7\,12}^a\\
\end{smallarray}\right]$$
where $p_{ij}^a=p_{ij}p_{31}^{-1}$ for $1\leq i\leq 7$ and $j=1,7,11,12$.

Now, take the elementary matrix 
$Y_{II, 1, 3}^{(1)}$ of Definition \ref{defelementary} with $D_{13}^{(1)}=-\begin{bsmallmatrix}
p_{32}&p_{33}& p_{34}
\end{bsmallmatrix}$. This is an elementary matrix of type II 
and 
$$P_1Y_I^{a}Y_{II,1,3}^{(1)}=
\left[\begin{smallarray}{cccccc|cccc|c|c}
    p_{11}^a& p_{12}^a&p_{13}^{a}& p_{14}^{a}&p_{15}& p_{16}
& p_{17}^a& p_{18}^{a}&p_{19}^{a}& p_{1\,10}^{a}&p_{1\,11}^a&p_{1\,12}^a\\
   p_{21}^a& p_{22}^a&p_{23}^{a}& p_{24}^{a}&p_{25}& p_{26}
& p_{27}^a& p_{28}^{a}&p_{29}^{a}& p_{2\,10}^{a}&p_{2\,11}^a &p_{2\,12}^a\\
    1&0&0& 0&p_{35}& p_{36}
& p_{37}^a& p_{38}^{a}&p_{39}^{a}& p_{3\,10}^{a}&p_{3\,11}^a &p_{3\,12}^a\\
 p_{41}^a& p_{42}^a&p_{43}^{a}& p_{44}^{a}&p_{45}& p_{46}
& p_{47}^a& p_{48}^{a}&p_{49}^{a}& p_{4\,10}^{a}&p_{4\,11}^a &p_{4\,12}^a\\
p_{51}^a& p_{52}^a&p_{53}^{a}& p_{54}^{a}&p_{55}& p_{56}
& p_{57}^a& p_{58}^{a}&p_{59}^{a}& p_{5\,10}^{a}&p_{5\,11}^a &p_{5\,12}^a\\
p_{61}^a& p_{62}^a&p_{63}^{a}& p_{64}^{a}&p_{65}& p_{66}
& p_{67}^a& p_{68}^{a}&p_{69}^{a}& p_{6\,10}^{a}&p_{6\,11}^a &p_{6\,12}^a\\
p_{71}^a& p_{72}^a&p_{73}^{a}& p_{74}^{a}&p_{75}& p_{76}
& p_{77}^a& p_{78}^{a}&p_{79}^{a}& p_{7\,10}^{a}&p_{7\,11}^a &p_{7\,12}^a\\
    \end{smallarray}\right].$$
Defining elementary matrices 
$Y_{II,1,4}^{(1)}$, $Y_{II,1,1}^{(2)}$, $Y_{II,1,3}^{(2)}$, $Y_{II,1,1}^{(3)}$ and
$Y_{II,1,1}^{(4)}$ of type II in a similar way we can zero out the remaining elements
of the third row:
$$
P_1^{a}
=\left[\begin{smallarray}{cccccc|cccc|c|c}
    p_{11}^a& p_{12}^a&p_{13}^a& p_{14}^a&p_{15}^a& p_{16}^a
& p_{17}^{a'}& p_{18}^{a'}&p_{19}^{a'}& p_{1\,10}^{a'}&p_{1\,11}^{a'}&p_{1\,12}^{a'}\\
   p_{21}^{a}& p_{22}^a&p_{23}^a& p_{24}^a&p_{25}^a& p_{26}^a
& p_{27}^{a'}& p_{28}^{a'}&p_{29}^{a'}& p_{2\,10}^{a'}&p_{2\,11}^{a'} &p_{2\,12}^{a'}\\
    1& 0&0& 0&0& 0&0& 0&0& 0&0&0\\
 p_{41}^{a}& p_{42}^a&p_{43}^a& p_{44}^a&p_{4,}^a& p_{46}^a
& p_{47}^{a'}& p_{48}^{a'}&p_{4,9}^{a'}& p_{4\,10}^{a'}&p_{4\,11}^{a'} &p_{4\,12}^{a'}\\
p_{51}^{a}& p_{52}^a&p_{53}^a& p_{54}^a&p_{55}^a& p_{56}^a
& p_{57}^{a'}& p_{58}^{a'}&p_{5,9}^{a'}& p_{5\,10}^{a'}&p_{5\,11}^{a'} &p_{5\,12}^{a'}\\
p_{61}^{a}& p_{62}^a&p_{63}^a& p_{64}^a&p_{65}^a& p_{66}^a
& p_{67}^{a'}& p_{68}^{a'}&p_{6,9}^{a'}& p_{6\,10}^{a'}&p_{6\,11}^{a'} &p_{6\,12}^{a'}\\
p_{71}^{a}& p_{72}^a&p_{73}^a& p_{74}^a&p_{75}^a& p_{76}^a
& p_{77}^{a'}& p_{78}^{a'}&p_{7,9}^{a'}& p_{7\,10}^{a'}&p_{7\,11}^{a'} &p_{7\,12}^a\\
  \end{smallarray}\right].$$
 
 Let  $P^a=\begin{bsmallmatrix}P^a_1\\P^a_2\\P^a_3\\P^a_4\end{bsmallmatrix} \in  \PP_{(W;\rr)}$.
 Then $P^a\sim P$. By Proposition \ref{propsameadmissiblerw}, $\uI$ is an admissible sequence of indices
 for $P^a$. Therefore $\det P_1^a((3,1,4,7),1:4)\neq 0$ and 
$T_3=
\begin{bsmallmatrix}
p_{12}^a&p_{13}^a& p_{14}^a\\
p_{42}^a&p_{43}^a& p_{44}^a\\
p_{72}^a&p_{73}^a& p_{74}^a\\
\end{bsmallmatrix}\in \Gl(3)$.
Put
$Y_{11}^{b}=\diag(1, T_3^{-1}, I_2)$, 
$Y_{22}^{b}=\diag(1, T_3^{-1})$, 
$Y_{33}^{b}=Y_{44}^{b}=1$ and
$Y_I^{b}=\diag(Y_{11}^{b}, Y_{22}^{b}, Y_{33}^{b}, Y_{44}^{b})$.
Then 
$$P_1^aY_I^{b}=
\left[\begin{smallarray}{cccccc|cccc|c|c}
     p_{11}^a& 1&0&0&p_{15}^a& p_{16}^a
& p_{17}^{a'}& p_{18}^b&p_{19}^b& p_{1\,10}^b&p_{1\,11}^{a'}&p_{1\,12}^{a'}\\
   p_{21}^{a}& p_{22}^b&p_{23}^b& p_{24}^b&p_{25}^a& p_{26}^a
& p_{7}^{a'}& p_{28}^b&p_{29}^b& p_{2\,10}^b&p_{2\,11}^{a'} &p_{2\,12}^{a'}\\
    1& 0&0& 0&0& 0&0& 0&0& 0&0&0\\
 p_{41}^{a}& 0&1&0&p_{45}^a& p_{46}^a
& p_{47}^{a'}& p_{48}^b&p_{49}^b& p_{4\,10}^b&p_{4\,11}^{a'} &p_{4\,12}^{a'}\\
p_{51}^{a}& p_{52}^b&p_{53}^b& p_{54}^b&p_{55}^a& p_{56}^a
& p_{57}^{a'}& p_{58}^b&p_{59}^b& p_{5\,10}^b&p_{5\,11}^{a'} &p_{5\,12}^{a'}\\
p_{61}^{a}& p_{62}^b&p_{63}^b& p_{64}^b &p_{65}^a& p_{66}^a
& p_{67}^{a'}& p_{68}^b&p_{69}^b& p_{6\,10}^b&p_{6\,11}^{a'} &p_{6\,12}^{a'}\\
p_{7,1}^{a}&0&0&1&p_{7,5}^a& p_{76}^a
& p_{77}^{a'}& p_{78}^b&p_{79}^b& p_{7\,10}^b&p_{7\,11}^{a'} &p_{7\,12}^a\\
    \end{smallarray}\right].$$
Take the elementary matrix
$Y_{II,3,4}^{(1)}$ with
  $
D_{34}^{(1)}=-\begin{bsmallmatrix}
    p_{15}^a& p_{16}^a\\
p_{45}^a& p_{46}^a\\
p_{75}^a& p_{76}^a
\end{bsmallmatrix}$.
Then
$$P_1^aY_I^{b}Y_{II,3,4}^{(1)}=
\left[\begin{smallarray}{cccccc|cccc|c|c}
     p_{11}^a& 1&0&0&0&0
& p_{17}^{a'}& p_{18}^b&p_{19}^b& p_{1\,10}^b&p_{1\,11}^{a'}&p_{1\,12}^{a'}\\
   p_{21}^{a}& p_{22}^b&p_{23}^b& p_{24}^b&p_{25}^{a'}& p_{26}^{a'}
& p_{27}^{a'}& p_{,8}^b&p_{29}^b& p_{2\,10}^b&p_{2\,11}^{a'} &p_{2\,12}^{a'}\\
    1& 0&0& 0&0& 0&0& 0&0& 0&0&0\\
 p_{41}^{a}& 0&1&0&0& 0
& p_{47}^{a'}& p_{,8}^b&p_{49}^b& p_{4\,10}^b&p_{4\,11}^{a'} &p_{4\,12}^{a'}\\
p_{51}^{a}& p_{52}^b&p_{53}^b& p_{54}^b&p_{55}^{a'}& p_{56}^{a'}
& p_{57}^{a'}& p_{58}^b&p_{59}^b& p_{5\,10}^b&p_{5\,11}^{a'} &p_{5\,12}^{a'}\\
p_{61}^{a}& p_{62}^b&p_{63}^b& p_{64}^b &p_{65}^{a'}& p_{66}^{a'}
& p_{67}^{a'}& p_{68}^b&p_{69}^b& p_{6\,10}^b&p_{6\,11}^{a'} &p_{6\,12}^{a'}\\
p_{71}^{a}&0&0&1&0&0
& p_{77}^{a'}& p_{78}^b&p_{79}^b& p_{7\,10}^b&p_{7\,11}^{a'} &p_{7\,12}^a\\
    \end{smallarray}\right].$$
Repeating the process with appropriate elementary matices
$Y_{II,3,3}^{(2)}$, $Y_{II,3,1}^{(3)}$ and $Y_{II,3,1}^{(4)}$ of type II we obtain
$$P_1^{b}=
\left[\begin{smallarray}{cccccc|cccc|c|c}
     p_{11}^a& 1&0&0&0&0
& p_{17}^{a'}& 0&0&0&0&0\\
   p_{21}^{a}& p_{22}^b&p_{23}^b& p_{24}^b&p_{25}^{a'}& p_{26}^{a'}
& p_{27}^{a'}& p_{28}^{b'}&p_{29}^{b'}& p_{2\,10}^{b'}&p_{2\,11}^{a'} &p_{2\,12}^{b}\\
    1& 0&0& 0&0& 0&0& 0&0& 0&0&0\\
 p_{41}^{a}& 0&1&0&0&0
& p_{47}^{a'}&0&0&0&0&0\\
p_{51}^{a}& p_{52}^b&p_{53}^b& p_{54}^b&p_{55}^{a'}& p_{56}^{a'}
& p_{57}^{a'}& p_{58}^{b'}&p_{59}^{b'}& p_{5\,10}^{b'}&p_{5\,11}^{a'} &p_{5\,12}^{b}\\
p_{61}^{a}& p_{62}^b&p_{63}^b& p_{64}^b &p_{65}^{a'}& p_{66}^{a'}
& p_{67}^{a'}& p_{68}^{b'}&p_{69}^{b'}& p_{6\,10}^{b'}&p_{6\,11}^{a'} &p_{6\,12}^{b}\\
p_{71}^{a}&0&0&1&0&0
& p_{77}^{a'}& 0&0&0&0&0\\
  \end{smallarray}\right].$$

Finally,
let  $ P^b=\begin{bsmallmatrix}P^b_1\\P^b_2\\P^b_3\\P^b_4\end{bsmallmatrix} \in  \PP_{(W;\rr)}$.
Then $P^b\sim P$. As before, by Proposition \ref{propsameadmissiblerw}, $\uI$ is an admissible sequence of indices for $P^b$ 
and $T_4=\begin{bsmallmatrix}p_{55}^{a'}& p_{56}^{a'}\\p_{65}^{a'}& p_{66}^{a'}\end{bsmallmatrix}\in \Gl(2)$.
Putting
$Y_{11}^{c}=\diag(1, I_3, T_4^{-1})$, 
$Y_{22}^{c}=\diag(1, I_3)$, 
$Y_{33}^{c}=Y_{44}^{c}=1$ and
$Y_I^{c}=\diag(Y_{11}^{c}, Y_{22}^{c}, Y_{33}^{c}, Y_{44}^{c})$,
$$P_1^{b}Y_I^{c}=
\left[\begin{smallarray}{cccccc|cccc|c|c}
     p_{11}^a& 1&0&0&0&0
& p_{17}^{a'}& 0&0&0&0&0\\
   p_{21}^{a}& p_{22}^b&p_{23}^b& p_{24}^b&p_{2,5}^{c}& p_{26}^{c}
& p_{2,7}^{a'}& p_{2,8}^{b'}&p_{29}^{b'}& p_{2\,10}^{b'}&p_{2\,11}^{b} &p_{2\,12}^{b}\\
    1& 0&0& 0&0& 0&0& 0&0& 0&0&0\\
 p_{41}^{a}& 0&1&0&0&0
& p_{47}^{a'}&0&0&0&0&0\\
p_{51}^{a}& p_{52}^b&p_{53}^b& p_{54}^b&1&0
& p_{57}^{a'}& p_{58}^{b'}&p_{59}^{b'}& p_{5\,10}^{b'}&p_{5\,11}^{b} &p_{5\,12}^{b}\\
p_{61}^{a}& p_{62}^b&p_{63}^b& p_{64}^b &0&1
& p_{67}^{a'}& p_{68}^{b'}&p_{69}^{b'}& p_{6\,10}^{b'}&p_{6\,11}^{b} &p_{6\,12}^{b}\\
p_{71}^{a}&0&0&0&1&0
& p_{77}^{a'}& 0&0&0&0&0\\
    \end{smallarray}\right].$$
Using appropriate  elementary matrices
$Y_{II,4,3}^{(2)}$ and $Y_{II,4,1}^{(4)}$, we get
   $$ P_1^{(\re)}=
\left[\begin{smallarray}{cccccc|cccc|c|c}
     p_{11}^{r}& 1&0&0&0&0& p_{17}^{r}& 0&0&0&0&0\\
   p_{21}^{r}& p_{22}^{r}&p_{23}^{r}& p_{24}^{r}&p_{25}^{r}& p_{26}^{r}& p_{27}^{r}& p_{28}^{r}&
   p_{29}^{r}& p_{2\,10}^{r}&p_{2\,11}^{r} &p_{2\,12}^{r}\\
    1& 0&0& 0&0& 0&0& 0&0& 0&0&0\\
 p_{41}^{r}& 0&1&0&0&0& p_{47}^{r}&0&0&0&0&0\\
p_{51}^{r}& p_{52}^{r}&p_{53}^{r}& p_{54}^{r}&1&0& p_{57}^{r}& 0&0&0&p_{5\,11}^{r} &0\\
p_{61}^{r}& p_{62}^{r}&p_{63}^{r}& p_{64}^r &0&1& p_{67}^{r}& 0&0&0&p_{6\,11}^{r} &0\\
p_{71}^{r}&0&0&1&0&0& p_{77}^{r}& 0&0&0&0&0\\
    \end{smallarray}\right].$$
    Let
     $$P^{(\re)}=\begin{bsmallmatrix}
P_1^{(r)}\\P_2^{(r)}\\P_3^{(r)}\\P_4^{(r)}
\end{bsmallmatrix}=\begin{bsmallmatrix}
P_1^{(r)}\\
I_{r_2, r_1}^TP_1^{(r)}\\
I_{r_3, r_1}^T{P_1^{(r)}}^2\\
I_{r_4, r_1}^T{P_1^{(r)}}^3
\end{bsmallmatrix}.$$
We will say that $P^{(\re)}$ is a matrix in {\em reduced form}. Observe that 
$$P_1^{(\re)}((3, 1, 4, 7, 5, 6),:)=
\begin{bsmallmatrix}
 1& 0&0& 0&0& 0&0& 0&0& 0&0&0\\
p_{11}^{r}& 1&0&0&0&0& p_{17}^{r}& 0&0&0&0&0\\
p_{41}^{r}& 0&1&0&0&0& p_{47}^{r}&0&0&0&0&0\\
p_{71}^{r}&0&0&1&0&0& p_{77}^{r}& 0&0&0&0&0\\
p_{51}^{r}& p_{52}^{r}&p_{53}^{r}& p_{54}^{r}&1&0& p_{5,7}^{r}& 0&0&0&p_{5\,11}^{r} &0\\
p_{61}^{r}& p_{62}^{r}&p_{63}^{r}& p_{64}^r &0&1& p_{67}^{r}& 0&0&0&p_{6\,11}&0
\end{bsmallmatrix}$$
$$
=
\left[
\begin{smallarray}{c|ccc|cc||c|ccc||c||c}
 1& 0&0& 0&0& 0&0& 0&0& 0&0&0\\\hline
p_{11}^{r}& 1&0&0&0&0& p_{17}^{r}& 0&0&0&0&0\\
p_{41}^{r}& 0&1&0&0&0& p_{47}^{r}&0&0&0&0&0\\
p_{71}^{r}&0&0&1&0&0& p_{77}^{r}& 0&0&0&0&0\\\hline
p_{51}^{r}& p_{52}^{r}&p_{53}^{r}& p_{54}^{r}&1&0& p_{57}^{r}& 0&0&0&p_{5\,11}^{r} &0\\
p_{61}^{r}& p_{62}^{r}&p_{63}^{r}& p_{64}^r &0&1& p_{67}^{r}& 0&0&0&p_{6\,11}&0
\end{smallarray}
\right]=
\left[
\begin{smallarray}{ccc|cc|c|c}
 1 &\nsc&\nsc&\nsc&\nsc&\nsc&\nsc\\
P_{3,1}^{(r,1)}& I_3&\nsc &P_{3,1}^{(r,2)}&\nsc&\nsc&\nsc\\
P_{4,1}^{(r,1)}& P_{4,3}^{(r,1)}&I_2 &P_{4,1}^{(r,2)}&\nsc&P_{4,1}^{(r,3)}&\nsc\\
\end{smallarray}
\right],
$$
with $P_{i,k}^{(r,j)}\in \RR^{(\tau_i-\tau_{i-1})\times (\tau_k-\tau_{k-1})}$, $1\leq i \leq m=4$, $1\leq k \leq m-1=3$, $1\leq k\leq m-j$, and the number of parameters in
$P^{(\re)}$ is $30=sr-N=84-54$.
}\hfill$\Box$
\end{example}

\bigskip 
With this example in mind we define the notion of \textit{reduced form} of a matrix in $\PP_{(W; \rr)}$ and show
that any matrix in this open set is $\wC_W$-equivalent to a matrix in reduced form. Recall that $r=r_1$.

Let $P=\begin{bsmallmatrix} P_1\\P_2\\\vdots\\P_k\end{bsmallmatrix}\in\PP_{(W; \rr)}$ with $P_1=\begin{bmatrix}
P_{11}& P_{12}&\cdots &P_{1m}\end{bmatrix}$ and $P_{1j}\in\RR^{r\times w_j}$, $1\leq j\leq m$. Let
$\uI=(\I_1, \dots, \I_{m})$ be an admissible sequence of indices for $P$, with $\I_j=(i_1, \dots, i_{\tau_j})$, $1\leq j\leq m$; 
in particular, $\I_m=(i_1,\ldots, i_{w_1})$. A matrix
$R=\begin{bsmallmatrix} R_1\\R_2\\\vdots\\R_k\end{bsmallmatrix}\in\PP_{(W; \rr)}$ with $R_1=\begin{bmatrix}
R_{11}& R_{12}&\cdots &R_{1m}\end{bmatrix}$ and $R_{1j}\in\RR^{r\times w_j}$, $1\leq j\leq m$, is said to be
a \textit{$\wC_W$-reduced form of $P$ with respect to $\uI$} if
\begin{equation}\label{eq.defR11}
R_{11}(\I_m,:)=\begin{bmatrix} I_{\tau_1} & 0 & 0 & \cdots & 0\\R_{21}^{(1)}& I_{(\tau_2-\tau_1)} & 0 &\cdots & 0\\
R_{31}^{(1)} & R_{32}^{(1)} & I_{(\tau_3-\tau_2)} &\cdots & 0\\\vdots&\vdots& \vdots &\ddots& \vdots\\
R_{m1}^{(1)} & R_{m2}^{(1)}& R_{m3}^{(1)} & \cdots & I{(\tau_m-\tau_{m-1})}\end{bmatrix}
\end{equation}
and for $j=2,3,\ldots, m$,
\begin{equation}\label{eq.defR1j}
\begin{array}{l}
R_{1j}(\I_m,:)=\\
\begin{blockarray}{ccccccc}
 \mbox{\scriptsize $\tau_1$}&
 \mbox{\scriptsize $\tau_2-\tau_1$}& \mbox{\scriptsize $\tau_3-\tau_2$}&\cdots &\mbox{\scriptsize $\tau_{m-j}-\tau_{m-j-1}$}&
\mbox{\scriptsize $\tau_{m-j+1}-\tau_{m-j}$}\\
\begin{block}{[cccccc]c}
0 & 0 & 0&\cdots &0&0&\mbox{\scriptsize $\tau_1$}\\
\vdots &\vdots&\vdots& & \vdots & \vdots &\vdots\\
0& 0 & 0&\cdots &0& 0&\mbox{\scriptsize $\tau_j-\tau_{j-1}$}\\
R_{j+1\, 1}^{(j)}&0&0& \cdots &0&0&\mbox{\scriptsize $\tau_{j+1}-\tau_j$}\\
R_{j+2\, 1}^{(j)} & R_{j+2\, 2}^{(j)} &0& \cdots & 0 & 0&\mbox{\scriptsize $\tau_{j+2}-\tau_{j+1}$}\\
\vdots&\vdots&\vdots&&\vdots &\vdots&\vdots\\ 
R_{m1}^{(j)} & R_{m2}^{(j)}&R_{m3}^{(j)}&\cdots &R_{m\,m-j}^{(j)} &0&\mbox{\scriptsize $\tau_{m}-\tau_{m-1}$}\\
\end{block}
\end{blockarray}
\end{array}
\end{equation}
\begin{rem}\label{rem.redform1}{\rm
\begin{itemize}
\item[(i)]
Note that $\I_m=(i_1,\ldots, i_{w_1})$ and so $R_{1j}(\I_m,:)\in\RR^{w_1\times w_j}$, $1\leq j\leq m$; in particular,
$R_{1m}(\I_m,:)=0\in\RR^{w_1\times w_m}$. 
\item[(ii)] Since $R$ is assumed to be in $\PP_{(W,\rr)}$, $R_j=I_{r_1,r_j}^TR_1W^{j-1}$ for $1\leq j\leq k$. Therefore
the $\wC_W$-reduced forms of $P$  with respect to $\uI$ are completely determined by $R_1$.
\item [(iii)] A detailed analysis of the zero-nonzero block pattern of $R_1$ yields the following 
characterization of the $\wC_W$-reduced forms of $P$ with respect to $\uI$:
\begin{equation}\label{eqfrm}
\begin{array}{ll}
R_{ii}^{(1) }=I_{(\tau_i-\tau_{i-1})}, &1\leq i \leq m,\\
R_{ik}^{(1) }=0,& 1\leq i \leq m-1,  \, i<k \leq m,\\
R_{ik}^{(j) }=0, & 1\leq i \leq m,\, 2\leq j \leq m, \\
& \max\{i-j+1, 1\}\leq k \leq m-j+1.
\end{array}
\end{equation}
The two first conditions mean that $R_{11}(\uI_m,:)$ has the form of \eqref{eq.defR11} and the third condition
means that $R_{1j}(\uI_m,:)$ has the form of \eqref{eq.defR1j}.
\end{itemize}\hfill$\Box$}
\end{rem}

\begin{theorem}\label{theofrreal}
Let $P\in \PP_{(W; \rr)}$ and let $\uI=(\I_1, \dots, \I_{m})$ be an admissible sequence of indices for $P$. Then $P$ is
$\wC_W$-equivalent to a unique $\wC_W$-reduced form with respect to $\uI$.
\end{theorem}

{\bf Proof.} Assume that $\I_j=(i_1, \dots, i_{\tau_j})$, $1\leq j \leq m$, let $s_j=\sum_{i=1}^jw_i$, $1\leq j \leq m$ 
and put $s_m=s$. Write

$$P(\I_m, :)=
 \left[
 \begin{smallarray}{ccc|ccc|c|c}
   P_{11}^{(1)}&\cdots &P_{1m}^{(1)}&P_{11}^{(2)}&\cdots &P_{1m-1}^{(2)}&\cdots &P_{11}^{(m)}\\
  P_{21}^{(1)}&\cdots &P_{2m}^{(1)}&P_{21}^{(2)}&\cdots &P_{2m-1}^{(2)}&\cdots &P_{21}^{(m)}\\
   \vdots  &\cdots &\vdots&   \vdots &\cdots &\vdots&   \cdots &\vdots\\
   P_{m1}^{(1)}&\cdots &P_{mm}^{(1)}&P_{m1}^{(2)}&\cdots &P_{mm-1}^{(2)}&\dots &P_{m1}^{(m)}\\
 \end{smallarray}\right],$$
 with 
 $P_{ik}^{(j)}\in \RR^{(\tau_i-\tau_{i-1})\times(\tau_k-\tau_{k-1})}$, $1\leq i, j\leq m$, $1\leq k \leq m-j+1$. 

If $\wh{P}=PY$ with $Y\in\wC_W$, then $\wh{P}(\I_m, :)=P(\I_m, :)Y$ and,
by Proposition \ref{propsameadmissiblerw}, $\uI$ is  an admissible sequence of indices for $\wh{P}$.

Since $\uI$ is an admissible sequence of indices for $P$, by \eqref{eqIrankr},
$\begin{bsmallmatrix}
P_{11}^{(1)}&\cdots &P_{1 j}^{(1)}\\
\vdots&\vdots &\vdots\\
P_{j1}^{(1)}&\cdots &P_{jj}^{(1)}\\
\end{bsmallmatrix}=P(\I_j, 1:\tau_j)\in \Gl(\tau_j)$, $1\leq j \leq m$.
We will prove by induction on $\ell$ that, for $1\leq \ell\leq m$, 
$P\sim P^{(\ell)}$
with
 \begin{equation}\label{eqfrlrw}
   \begin{array}{ll}
P_{ii}^{(\ell, 1)}=I_{(\tau_i-\tau_{i-1})}, & 1\leq i \leq \ell,\\
P_{ik}^{(\ell, 1)}=\nsc,& 1\leq i \leq \ell, \quad i<k \leq m_,\\
P_{ik}^{(\ell, j)}=\nsc, & 1\leq i \leq \ell, \; 2\leq j \leq m, \\
   &\max\{i-j+1, 1\}\leq k \leq m-j+1.
   \end{array}
 \end{equation}
 Taking into account \eqref{eqfrm}, this will prove that $P^{(m)}$ is a $\wC_W$-reduced form
 of $P$ with respect to $\uI$.
\begin{itemize}
 \item For $\ell=1$, we have $P_{11}^{(1)}\in \Gl(\tau_1)$. Let  $T_1={P_{11}^{(1)}}^{-1}$,
 $Y_{ii}=\diag(T_1, I_{(\tau_2-\tau_1)},$ $ \dots,  I_{(\tau_{m-i+1}-\tau_{m-i})})$, $1\leq i \leq m$ and 
$Y_I^{(a)}=\diag(Y_{11}, Y_{22}, \dots,  Y_{mm})$. Then  
$Y_I^{(a)}$ is an elementary matrix of type I and
$$
P(\I_m, :)Y_I^{(a)}=
 \left[
 \begin{smallarray}{cccc|ccc|c|c}
   I_{\tau_1}&P_{12}^{(a,1)}&\cdots &P_{1m}^{(a,1)}&P_{11}^{(a,2)}&\cdots &P_{1m-1}^{(a,2)}&
   \cdots &P_{11}^{(a,m)}\\
  P_{21}^{(a,1)}&P_{22}^{(a,1)}&\cdots &P_{2m}^{(a,1)}&P_{21}^{(a,2)}&\cdots &P_{2m-1}^{(a,2)}&
  \cdots &P_{21}^{(a,m)}\\
   \vdots  &\vdots &\vdots&   \vdots &\vdots &\vdots&   \vdots &\vdots&\vdots\\
   P_{m1}^{(a,1)}&P_{m2}^{(a,1)}&\cdots &P_{mm}^{(a,1)}&P_{m1}^{(a,2)}&\cdots &P_{mm-1}^{(a,2)}&
   \cdots &P_{m1}^{(a,m)}\\
 \end{smallarray}\right],
$$
where 
$P_{i1}^{(a,j)}=P_{i1}^{(j)}T_1$, $1\leq i,j\leq m$ and 
$P_{ik}^{(a,j)}$=$P_{ik}^{(j)}$, $1\leq i,j\leq m$, $2\leq k\leq m-j+1$.

Now, take the elementary matrix 
$Y_{II, 1, 2}^{(1)}$ with $D_{12}^{(1)}=-P_{12}^{(a, 1)}$. Then
$$
P(\I_m, :)Y_I^{(a)}Y_{II, 1, 2}^{(1)}=
\left[
 \begin{smallarray}{cccc|ccc|c|c}
   I_{\tau_1}&\nsc&\dots &P_{1m}^{(b,1)}&P_{1,1}^{(b,2)}&\cdots &P_{1m-1}^{(b,2)}&\cdots &P_{11}^{(b,m)}\\
  P_{21}^{(b,1)}&P_{22}^{(b,1)}&\cdots &P_{2m}^{(b,1)}&P_{21}^{(b,2)}&\cdots &P_{2m-1}^{(b,2)}&
  \cdots &P_{21}^{(b,m)}\\
   \vdots  &\vdots &\vdots&   \vdots &\vdots &\vdots&   \vdots &\vdots&\vdots\\
   P_{m1}^{(b,1)}&P_{m2}^{(b,1)}&\cdots &P_{mm}^{(b,1)}&P_{m1}^{(b,2)}&\cdots &P_{mm-1}^{(b,2)}&
   \cdots &P_{m1}^{(b,m)}\\
 \end{smallarray}\right],
$$
where 
$P_{i, 2}^{(b,j)}=P_{i, 2}^{(a,j)}-P_{i, 1}^{(a,j)}P_{1,2}^{(a, 1)}$, $1\leq i,j\leq m$ and 
$P_{i, k}^{(b,j)}$=$P_{i, k}^{(a,j)}$, $1\leq i,j\leq m$, $1\leq k\leq m-j+1$, $k\neq 2$.

Repeating the same transformation with the appropriate elementary matrices
$Y_{II,1,k}^{(1)}$, $3\leq k \leq m$ and $Y_{II,1,k}^{(j)}$, $2\leq j \leq m$, $1\leq k \leq m-j+1$, we get
$$P^{(1)}(\I_m, :)=
\left[
 \begin{smallarray}{cccc|ccc|c|c}
   I_{\tau_1}&\nsc&\cdots&\nsc&\nsc&\cdots&\nsc&\cdots&\nsc\\
  P_{21}^{(1,1)}&P_{22}^{(1,1)}&\cdots &P_{2m}^{(1,1)}&P_{21}^{(1,2)}&\cdots &P_{2m-1}^{(1,2)}&
  \cdots &P_{21}^{(1,m)}\\
   \vdots  &\vdots &\vdots&   \vdots &\vdots &\vdots&   \vdots &\vdots&\vdots\\
   P_{m1}^{(1,1)}&P_{m2}^{(1,1)}&\cdots &P_{mm}^{(1,1)}&P_{m1}^{(1,2)}&\cdots &P_{mm-1}^{(1,2)}&
   \cdots &P_{m,1}^{(1,m)}\\
 \end{smallarray}\right].
$$
Then, $P^{(1)}$ satisfies (\ref{eqfrlrw}) for $\ell=1$.

\item
Assume now that $\ell\in \{2, \ldots, m-1\}$ and  $P\sim P^{(\ell)}$ with $P^{(\ell)}$ satisfying (\ref{eqfrlrw}), i.e.
\[
P^{(\ell)}(\uI_m, :)=\begin{bmatrix}P_{11}^{(\ell)} &P_{12}^{(\ell)}&\cdots & P_{1m}^{(\ell)}\end{bmatrix}
\]
where
\[
P_{11}^{(\ell)}=\begin{bsmallmatrix} I_{\tau_1} & \nsc & \cdots & \nsc&\nsc& \cdots & \nsc\\
P_{21}^{(\ell,1)} & I_{(\tau_2-\tau_1)} & \cdots & \nsc&\nsc& \cdots & \nsc\\
\vdots &\vdots &\ddots & \vdots&\vdots&&\vdots\\
P_{\ell1}^{(\ell,1)}& P_{\ell2}^{(\ell,1)}&\cdots & I_{(\tau_\ell-\tau_{\ell-1})}& \nsc & \cdots & \nsc\\
P_{\ell+1\,1}^{(\ell,1)}& P_{\ell+1\,2}^{(\ell,1)}&\cdots &P_{\ell+1\,\ell}^{(\ell,1)}&P_{\ell+1\,\ell+1}^{(\ell,1)}&
\cdots &P_{\ell+1\,m}^{(\ell,1)}\\
\vdots &\vdots & & \vdots&\vdots&&\vdots\\
P_{m1}^{(\ell,1)}& P_{m2}^{(\ell,1)}&\cdots &P_{m\ell}^{(\ell,1)}&P_{m\,\ell+1}^{(\ell,1)}&
\cdots &P_{mm}^{(\ell)}
\end{bsmallmatrix},
\]
\[
P_{1j}^{(\ell)}=\begin{bsmallmatrix} \nsc& \nsc&\cdots & \nsc & \nsc &\cdots & \nsc\\
\vdots&\vdots& &\vdots&\vdots&&\vdots\\
\nsc& \nsc&\cdots & \nsc & \nsc &\cdots & \nsc\\
P_{j+1\, 1}^{(\ell,j)}& \nsc & \cdots & \nsc & \nsc & \cdots &\nsc\\
\vdots&\vdots& &\vdots&\vdots&&\vdots\\
P_{\ell\, 1}^{(\ell,j)} &P_{\ell\, 2}^{(\ell,j)} &\cdots&P_{\ell\, \ell-j}^{(\ell,j)}&\nsc &\cdots & \nsc\\
P_{\ell+1\, 1}^{(\ell,j)} &P_{\ell+1\, 2}^{(\ell,j)} &\cdots&P_{\ell+1\, \ell-j}^{(\ell,j)}&P_{\ell+1\, \ell-j+1}^{(\ell,j)}&
\cdots & P_{\ell+1\; m-j+1}^{(\ell,j)}\\
\vdots&\vdots& &\vdots&\vdots&&\vdots\\
P_{m1}^{(\ell,j)} &P_{m 2}^{(\ell,j)} &\cdots&P_{m\, \ell-j}^{(\ell,j)}&P_{m\, \ell-j+1}^{(\ell,j)}&
\cdots & P_{m\, m-j+1}^{(\ell,j)}
\end{bsmallmatrix},\quad 2\leq j\leq \ell-1,
\]

and
\[
P_{1j}^{(\ell)}=\begin{bsmallmatrix} 
\nsc & \nsc &\cdots & \nsc\\
\vdots&\vdots & &\vdots\\
\nsc&\nsc&\cdots & \nsc\\
P_{\ell+1\, 1}^{(\ell,j)}&P_{\ell+1\, 2}^{(\ell,j)} &\cdots & P_{\ell+1\, m-j+1}^{(\ell,j)}\\
\vdots&\vdots & &\vdots\\
P_{m1}^{(\ell,j)} &P_{m 2}^{(\ell,j)} &\cdots&P_{m\, m-j+1}^{(\ell,j)}
\end{bsmallmatrix},\quad \ell\leq j \leq m.
\]
By Proposition \ref{propsameadmissiblerw}, $\uI$ is an admissible sequence of indices for $P^{(\ell)}$. Thus
$$ P_{11}^{(\ell)}((1,\ldots, \tau_{\ell+1}),(1,\ldots, \tau_{\ell+1}))=\begin{bsmallmatrix}
   I_{\tau_1}&\dots &\nsc&\nsc\\
   \vdots  &\ddots &\vdots&  \vdots\\
   P_{\ell 1}^{(\ell,1)}&\cdots & I_{(\tau_{\ell}-\tau_{\ell-1})}&\nsc\\
   P_{\ell+1\,1}^{(\ell,1)}&\cdots & P_{\ell+1\,\ell}^{(\ell,1)}& P_{\ell+1\,\ell+1}^{(\ell,1)}
    \end{bsmallmatrix}\in \Gl(\tau_{\ell+1}),
 $$
 and, consequently, $P_{\ell+1\,\ell+1}^{(\ell,1)}\in \Gl(\tau_{\ell+1}-\tau_\ell)$.
Let $T_{\ell+1}={P_{\ell+1\,\ell+1}^{(\ell,1)}}^{-1}$,
$$
\begin{array}{rll}
Y'_{jj}=&\diag(I_{\tau_1}, \dots , T_{\ell+1} , \dots,  I_{(\tau_{m-j+1}-\tau_{m-j})}), &1\leq j \leq m- \ell, \\
Y'_{jj}=&\diag(I_{\tau_1}, \dots,   I_{(\tau_{m-j+1}-\tau_{m-j})}), & m-\ell+1\leq j\leq m, 
\end{array}
$$
and 
$Y_I^{(a')}=\diag(Y'_{11}, Y'_{22}, \dots,  Y'_{mm})$.
For $1\leq j \leq m-\ell, $ the matrix $P_{1j}^{(\ell)}Y'_{jj}$ is obtained form $P_{1j}^{(\ell)}$
by replacing the block $P_{i \, \ell+1}^{(\ell, j)}$ by  the block $P_{i\, \ell+1}^{(\ell, j)}T_{\ell+1}$  ($\ell+1\leq i\leq m$),
and for $m-\ell+1\leq j\leq m$, $P_{1j}^{(\ell)}Y'_{jj}=P_{1j}^{(\ell)}$.

Now, using successively appropriate matrices $Y_{II, \ell+1, k}^{(1)}$, $\ell+2\leq k \leq m$, we can annihilate the
blocks $P_{\ell+1\, \ell+2}^{(\ell,1)}, \ldots, P_{\ell+1\, m}^{(\ell,1)}$. Similarly, with appropriate
matrices $Y_{II, \ell+1, k}^{(j)}$, $2\leq j \leq m$, $\max\{\ell-j+2, 1\}\leq k \leq m-j+1$,
we can annihilate the blocks 
$P_{\ell+1\,\ell-j+2}^{(\ell,j)}$, \ldots, $P_{\ell+1\,m-j+1}^{(\ell,1)}$ for $j=2, \ldots \ell$
and $P_{\ell+1\, 1}^{(\ell,j)}$,\ldots, $P_{\ell+1\, m-j+1}^{(\ell,j)}$ for
$j=\ell+1,\ldots m$.
Therefore (\ref{eqfrlrw}) holds for $1\leq\ell\leq m$. 
\end{itemize}

Setting $P^{(\re)}=P^{(m)}$, $P^{(\re)}$ satisfies (\ref{eqfrm}) and so it is a $\wC_W$-reduced form of $P$ with
 respect to $\uI$.

Let us see now that the matrix $P^{(\re)}$ is unique, that is to say, that if $P\sim \wh{P}^{(\re)}$
with
\[
\begin{array}{ll}
\wh{P}_{ii}^{(\re,1)}=I_{(\tau_i-\tau_{i-1})}, &1\leq i \leq m,\\
\wh{P}_{ik}^{(\re,1)}=0,& 1\leq i \leq m-1,  \quad i<k \leq m,\\
\wh{P}_{ik}^{(\re,j)}=0, & 1\leq i \leq m,\quad 2\leq j \leq m, \\
&\max\{i-j+1, 1\}\leq k \leq m-j+1.
\end{array}
\]
then
$P^{(\re)}=\wh{P}^{(\re)}$.

As $P\sim \wh{P}^{(\re)}$,  there exists $Y\in \wC_{W}$ such that
$P^{(\re)}Y=\wh{P}^{(\re)}$; in particular,
\[
P^{(\re)}(\I_m, :)Y=\wh{P}^{(\re}(\I_m, :),
\]
where $Y$ is the matrix of \eqref{eq.QXQr} and so, it satisfies the properties \eqref{eq.QYQm}, \eqref{eq.qxqdijr}
and \eqref{eq.qxqyijr}. It is then enough to prove that $Y_{11}=I_{w_1}$ and $Y_{1j}=0$ for $j=2,\ldots, m$ because,
by \eqref{eq.qxqyijr}, this would imply that $Y=I_s$ (recall that $s=\sum_{i=1}^m w_i$).

We can split the columns of  $P^{(\re)}(\I_m, :)$ and $\wh{P}^{(\re}(\I_m, :)$ as follows
\[
P^{(\re)}(\I_m, :)=\begin{bmatrix} R_1 & R_2 &\cdots & R_m\end{bmatrix},\quad
\wh P^{(\re)}(\I_m, :)=\begin{bmatrix} \wh R_1 &\wh R_2 &\cdots &\wh R_m\end{bmatrix},
\]
with $R_j,\wh R_j\in\RR^{w_1\times w_j}$, $1\leq j\leq m$. Then $R_1$ and $\wh R_1$ are lower block-triangular
matrices with identity matrices as diagonal blocks (cf. \eqref{eq.defR11}) and $Y_{11}=R_1^{-1}\wh R_1$
is also a lower block-triangular
matrices with identity matrices as diagonal blocks. However, by definition (see \eqref{eq.QYQm}),
$Y_{11}$  is  an upper block-triangular
matrix whose blocks are of the same size as the blocks of $R_1$ and $\wh R_1$. Hence $Y_{11}=I_{w_1}$ and
by \eqref{eq.qxqyijr}, $Y_{jj}=I_{w_j}$ for $1\leq j\leq m$.

Let us prove now by induction that, for $j\in\{2,3,\ldots, m\}$, $Y_{1j}=0$.
In fact, $\wh R_2= R_1 Y_{12}+R_2 Y_{22}=R_1 Y_{12}+R_2$ because $Y_{22}=I_{w_2}$.  Thus,
$Y_{12}= R_1^{-1}(\wh R_2-R_2)$. Now, by \eqref{eq.QYQm} and  \eqref{eq.defR1j} $Y_{12}$ and
$\wh R_2-R_2$ are matrices with the form of \eqref{eq.QYQm} and \eqref{eq.defR1j} with $j=2$, respectively.
Since $R_1^{-1}$ is a lower block-triangular matrix, $R_1^{-1}(\wh R_2-R_2)$
has the same zero-nonzero block pattern as $\wh R_2-R_2$. Therefore $Y_{12}=\wh R_2-R_2=0$.

Assume that $Y_{1j}=0$ for $j=1,\ldots,\ell-1$ with $3\leq \ell\leq m$. By \eqref{eq.qxqyijr}, $Y_{i\ell}=0$
for $i=2,\ldots,\ell-1$ and so $\wh R_\ell= R_1 Y_{1\ell}+R_\ell Y_{\ell\ell}=R_1 Y_{1\ell}+R_\ell$. As above,
$Y_{1\ell}=R_1^{-1}(\wh R_\ell-R_\ell)$ and since $Y_{1\ell}$ and $R_1^{-1}(\wh R_\ell-R_\ell)$ have complementary
zero-nonzero block structures, $Y_{1\ell}=\wh R_\ell-R_\ell=0$. Therefore $Y=I_s$ and $P^{(\re)}=\wh{P}^{(\re)}$.
\hfill $\Box$

\begin{rem}\label{remnpredform}{\rm
The number of parameters of $P^{(\re)}$ is the number of parameters of $P$ minus the number of parameters
of $Y$. That is to say $\dim \PP_{(W;\rr)}- \dim C_{W}= rs-N_W$, where $s=\sum_{i=1}^m w_i$ and
$N_W=\sum_{i=1}^m w_i^2$ (see Remark \ref{rem.comweyr}). Note that it follows from Proposition \ref{prop.parnemp}
 that $r=r_1\geq w_i$, $1\leq i\leq m$ and so $rs\geq N_W$.
}\hfill $\Box$
\end{rem}

\subsection{Reduced form when $M$ has two conjugated complex eigenvalues}
\label{subsecredcomplex}

Let
$\whW= \whW (\lambda, \overline{\lambda}))$, $\lambda =a+bi \in \CC\setminus \RR$, with 
Weyr characteristic 
$(w_1, \dots, w_m)$  and $B=\begin{bsmallmatrix}a&b\\-b&a\end{bsmallmatrix}$.
As before, 
  $\tau_i=w_{m-i+1}$, $0\leq i \leq m$ ($w_{m+1}=0$), $s_j=\sum_{i=1}^jw_i$, $1\leq j \leq m$  ($s_m=s$) and assume that 
  $\PP_{(\whW;\rr)}\neq \emptyset$ (see Proposition \ref{prop.parnemp}); in particular, $r_1\geq w_1$.

\medskip
We will use additional notation. Let $B_0=\begin{bsmallmatrix}0&1\\-1&0\end{bsmallmatrix}$.
Given a row vector $z=\begin{bmatrix} z_1 & z_2\end{bmatrix}\in  \RR^{1\times 2}$,
$\dZ$ denotes the matrix $\dZ=\begin{bsmallmatrix} z\\zB_0\end{bsmallmatrix}=\begin{bsmallmatrix}
z_1&z_2\\-z_2&z_1\end{bsmallmatrix}\in \RR^{2\times 2}$. As noted in Section \ref{sec.centralizer}
(see the note after Lemma \ref{lemmazcomplex}),
$C_B=C_{B_0}=\{\dZ: {\color{blue}z}\in \RR^{1\times2}\}$. Since $\det \dZ=z_1^2+z_2^2$, $\dZ\in \wC_B$ if and only if 
$z\neq \nsc$.

If $Z=\begin{bsmallmatrix}
z_{11}&\cdots&z_{1n}\\\vdots & &\vdots \\
z_{m1}&\cdots&z_{mn}\\
\end{bsmallmatrix}
\in  \RR^{m\times 2n}$, with $z_{ij}\in \RR^{1\times 2}$, $1\leq i \leq m$, $1\leq j \leq n$, $\dZ$ is the 
matrix 
$\dZ=\begin{bsmallmatrix}
\dZ_{11}&\cdots&\dZ_{1n}\\\vdots & &\vdots \\
\dZ_{m1}&\cdots&\dZ_{mn}\\
\end{bsmallmatrix}
\in  \RR^{2m\times 2n}$.  Recall (see \eqref{eq.defX(n)}) that $I_n\otimes B=\diag(\overbrace{B, \dots, B}^n)$.
Since $C_B=C_{B_0}=\{\dZ: z\in \RR^{1\times2}\}$ it is easy to see that
\[
C_{B^{(n)}}=C_{B_0^{(n)}}=\{\dZ: Z\in \RR^{n\times2n}\}
\]

\medskip
Recall (Lemma \ref{lemmazcomplex})  that  $Y\in C_{\whW}$ if and only if $Y$ has the structure of 
\eqref{eq.QXQr} satisfying the properties \eqref{eq.QYQm},  \eqref{eq.qxqyijr} and
for $1\leq i, j\leq m$ and $\max\{i-j+1, 1\}\leq k \leq m-j+1$, (see \eqref{eq.qxqdijc})
$D^{(j)}_{i,k}=\begin{bmatrix}T_{\alpha\,\beta}^{(j)}\end{bmatrix}_{\begin{smallarray}{l}
 \tau_{i-1}+1\leq \alpha\leq \tau_i\\\tau_{k-1}+1\leq \beta\leq\tau_k\end{smallarray}}
 \in \RR^{2(\tau_i-\tau_{i-1})\times 2(\tau_k-\tau_{k-1})}$ and $T_{\alpha, \beta}^{(j)}=
\begin{bsmallmatrix}x_{\alpha\, \beta}^{(j)}&y_{\alpha\, \beta}^{(j)}\\-y_{\alpha\, \beta}^{(j)}&x_{\alpha\, \beta}^{(j)}
\end{bsmallmatrix}\in \RR^{2\times 2}$. Therefore, we can write
$
D^{(j)}_{i, k}={Z_{i,k}^{(j)}}^{\diamond} \in \RR^{2(\tau_i-\tau_{i-1})\times 2(\tau_k-\tau_{k-1})} 
\mbox{ with }Z_{i,k}^{(j)}\in \RR^{(\tau_i-\tau_{i-1})\times 2(\tau_k-\tau_{k-1})}, 
1 \leq i, j \leq m, \quad \max\{i-j+1, 1\}\leq k \leq m-j+1.
$
In addition $Y\in \wC_{\whW}$ if and only if 
  $$
  D_{i,i}^{(1)}\in \wC_{B^{(\tau_i-\tau_{i-1})}} 
  \quad 1\leq i \leq m.
 $$
\begin{definition}\label{defelementary2}\,
\begin{enumerate}
\item 
Let $T_i \in  \wC_{B^{(\tau_i-\tau_{i-1})}}$, 
$1\leq i \leq m$ and
$\widehat{Y_I}=\diag(\widehat{Y_{11}}, \dots,  \widehat{Y_{mm}})$ with 
$\widehat{Y_{ii}}=\diag(T_1, \dots,  T_{m-i+1})$, $1\leq i \leq m$. The matrices of this type will
be called {\em elementary matrices of type I} and they form a subgroup of $\wC_{\whW}$.
\item 
For $j=1$, $1\leq i <k\leq m$, and for
$2\leq j \leq m$, $1\leq k \leq m-j+1$ $1\leq i \leq k+j-1$, let 
$\widehat{Y_{II,i,k}^{(j)}}$ be a matrix of  \eqref{eq.QYQm} with, perhaps, $D_{ik}^{(j)}\neq 0$,
$$
D_{ii}^{(1)}=I_{2(\tau_i-\tau_{i-1})}, \quad 1\leq i \leq m,
$$
and all the other blocks zero. This type of matrices will be called 
{\em elementary matrices of type II} and they form a subgroup of $\wC_{\whW}$. 
\end{enumerate}
\end{definition}

\begin{lemma}\label{lemmaRdi}
Given an  integer $i\geq 2$ and  $Z \in \RR^{m\times 2n}$,  
  $$\rank  \begin{bsmallmatrix}Z\\ZB^{(n)} \\\vdots \\Z{B^{(n)}}^{i-1}\end{bsmallmatrix}=
   \rank \begin{bsmallmatrix}Z\\ZB^{(n)} 
 \end{bsmallmatrix}.$$
\end{lemma}
{\bf Proof.}
If $Z=\begin{bmatrix}Z_{1}&\dots&Z_{n}\end{bmatrix}\in  \RR^{m\times 2n}$, with 
$Z_{j}\in \RR^{m\times 2}$, $1\leq j \leq n$, then
$$\begin{bsmallmatrix}Z\\ZB^{(n)} \\\vdots \\Z{B^{(n)}}^{i-1}\end{bsmallmatrix}=
\begin{bsmallmatrix}
Z_{1}&\cdots&Z_{n}\\Z_{1}B&\cdots&Z_{n}B\\
\vdots & &\vdots \\
Z_{1}B^{i-1}&\cdots&Z_{n}B^{i-1}
\end{bsmallmatrix}.$$
It follows from  $B^2=2aB-(a^2+b^2)I_2$ that $Z_jB^k=2aZ_jB^{k-1}-(a^2+b^2)Z_jB^{k-2}$,
for $1\leq j \leq n$ and $k\geq 2$. Therefore there exists $S\in \Gl(mi)$ such that
$$S\begin{bsmallmatrix}Z\\ZB^{(n)} \\\vdots \\Z{B^{(n)}}^{i-1}\end{bsmallmatrix}=
\begin{bsmallmatrix}Z\\ZB^{(n)} \\
\nsc\\\vdots \\\nsc\end{bsmallmatrix},$$
and the lemma follows.
\hfill $\Box$

\begin{proposition}\label{propauxiliarrw22}
Let
\begin{equation}\label{eq.P}
P=\begin{bmatrix}P_1\\\vdots \\P_k\end{bmatrix}\in\PP_{(\whW,\rr)},\;
P_i=\begin{bmatrix}P_{i1}& P_{i2}&\cdots &P_{im}\end{bmatrix}
\end{equation}
with $P_{ij}\in\RR^{r_i\times 2w_j}$, $1\leq i\leq k$ and $1\leq j\leq m$.
Then $\rank P_{11}^{\diamond}=2w_1$.
\end{proposition}
{\bf Proof.}
Since $P\in\PP_{(\whW,\rr)}$, $\rank P=\sum_{j=1}^m 2w_j$ and so
$\rank \begin{bsmallmatrix}P_{1j}\\\vdots \\P_{kj}\end{bsmallmatrix}=2w_j$, $1 \leq j \leq m$.
On the other hand, it follows from $ P_{i+1}=I_{r_{1}, r_{i+1}}^TP_{1}\whW^{i}$,  $1\leq i \leq k-1$,
that $P_{i+1\,1}=I_{r_{1}, r_{i+1}}^TP_{11}{B^{(w_1)}}^i$,  $1\leq i \leq k-1$. Thus,
$$
\begin{bsmallmatrix}P_{11}\\\vdots \\P_{k1}\end{bsmallmatrix}=
\begin{bsmallmatrix}P_{11}\\I_{r_{1}, r_{2}}^T P_{11}B^{(w_1)}\\\vdots \\
I_{r_{1}, r_{k}}^T P_{11}{B^{(w_1)}}^{k-1}
\end{bsmallmatrix}=\diag(I_{r_1}, I_{r_{1}, r_{2}}^T\,\dots,  
I_{r_{1}, r_{k}}^T)\begin{bsmallmatrix}P_{11}\\P_{11}{B^{(w_1)}}\\\vdots \\
P_{11}{B^{(w_1)}}^{k-1}
\end{bsmallmatrix}.
$$
Hence
$$
2w_1=\rank  \begin{bsmallmatrix}P_{11}\\\vdots \\P_{k1}\end{bsmallmatrix}=
\rank  \begin{bsmallmatrix}P_{11}\\P_{11}{B^{(w_1)}}\\\vdots \\P_{11}{B^{(w_1)}}^{k-1}\end{bsmallmatrix}.
$$
By Lemma \ref{lemmaRdi},
$\rank  \begin{bsmallmatrix}P_{11}\\P_{11}{B^{(w_1)}}\end{bsmallmatrix}=2w_1$.

Put
\begin{equation}\label{eq.P11}
P_{11}=
\begin{bsmallmatrix}z_{11}&\cdots&z_{1w_1}\\
\vdots&&\vdots\\
z_{r1}&\cdots&z_{rw_1}\end{bsmallmatrix},
 \quad z_{ij}\in \RR^{1\times 2},
\quad 1\leq i \leq r, \, 1\leq j \leq w_1.
\end{equation}
Then
$$
2w_1=\rank  \begin{bsmallmatrix}P_{11}\\P_{11}{B^{(w_1)}}\end{bsmallmatrix}=
\rank\begin{bsmallmatrix}z_{11}&\cdots&z_{1w_1}\\
\vdots& &\vdots\\
z_{r1}&\cdots&z_{rw_1}\\\hline\\
z_{11}B&\cdots&z_{1w_1}B\\
\vdots& &\vdots\\
z_{r1}B&\cdots&z_{rw_1}B
\end{bsmallmatrix}=
\rank\begin{bsmallmatrix}z_{11}&\dots&z_{1w_1}\\
z_{11}B&\cdots&z_{1w_1}B\\\hline
\vdots& &\vdots\\\hline \\
z_{r1}&\cdots&z_{rw_1}\\
z_{r1}B&\cdots&z_{rw_1}B
\end{bsmallmatrix}.
$$
Let $T=\begin{bsmallmatrix}1&0\\ -\frac{a}{b}&\frac{1}{b}\end{bsmallmatrix}$. Then
$$
T\begin{bsmallmatrix}z_{ij}\\z_{ij}B \end{bsmallmatrix}=
\begin{bsmallmatrix}z_{ij}\\z_{ij}B_0 \end{bsmallmatrix}, \quad 1\leq i\leq r, \quad1\leq j \leq w_1.
$$
Therefore
$$
\diag(\overbrace{T, \dots, T}^r)
\begin{bsmallmatrix}z_{11}&\cdots&z_{1w_1}\\
z_{11}B&\cdots&z_{1w_1}B\\\hline
\vdots& &\vdots\\\hline\\
z_{r1}&\cdots&z_{r w_1}\\
z_{r1}B&\cdots&z_{rw_1}B
\end{bsmallmatrix}=
\begin{bsmallmatrix}z_{11}&\cdots&z_{1w_1}\\
z_{11}B_0&\cdots&z_{1,w_1}B_0\\\hline
\vdots&&\vdots\\\hline\\
z_{r1}&\cdots&z_{r w_1}\\
z_{r1}B_0&\cdots&z_{rw_1}B_0
\end{bsmallmatrix}=
\begin{bsmallmatrix}
\dZ_{11}&\cdots&\dZ_{1w_1}\\
\vdots & &\vdots \\
\dZ_{r1}&\cdots&\dZ_{rw_1}\\
\end{bsmallmatrix}=P_{11}^{\diamond}.
$$
As $\diag(\overbrace{T, \dots, T}^r)\in \Gl(2r)$, $\rank  P_{11}^{\diamond}=2w_1$.
 \hfill $\Box$

\begin{proposition}\label{corindices}
Let  $P\in \PP_{(W;\rr)}$. Then, for each $j=1,\ldots, m$, there is a sequence of $\tau_j$ indices
$\I_j\subseteq \{1, \dots, r\}$ satisfying \eqref{eqIsubsetr}, \eqref{eqIsetminus} and
\begin{equation}\label{eqIrank}
   P(\I_j, 1:2\tau_j)^\diamond\in \Gl(2\tau_j), \quad 1\leq j \leq m
\end{equation}
\end{proposition}
 
{\bf Proof.}
Let $P\in \PP_{(\whW;\rr)}$ be the matrix of \eqref{eq.P} and let $P_{11}\in\RR^{r\times 2w_1}$ be that of \eqref{eq.P11}.
Write $P_{11}=\begin{bmatrix} P_{11}^{(1)} & P_{11}^{(2)}&\cdots &P_{11}^{(m)}
\end{bmatrix}$ with $P_{11}^{(j)}\in\RR^{r_1\times 2(\tau_{j}-\tau_{j-1})}$, $1\leq j\leq m$ and
$X_j=\begin{bmatrix} P_{11}^{(1)} & P_{11}^{(2)}&\cdots &P_{11}^{(j)}\end{bmatrix}\in\RR^{r\times 2\tau_j}$,
$1\leq j\leq m$. We claim that $\rank X_j=\tau_j$, $1\leq j\leq m$. In fact, as in the proof of Proposition \ref{propauxiliarrw22},
$\rank P_{11}^\diamond(:,1:2\tau_j)=\rank \begin{bsmallmatrix}X_j\\X_jB_0^{(\tau_j)}\end{bsmallmatrix}$.
Since $\rank P_{11}^\diamond=2w_1$,  bearing in mind that $B_0$ is invertible,
\[
\begin{array}{rcl}
2\tau_j&=&\rank P_{11}^\diamond(:,1:2\tau_j)=\rank \begin{bmatrix}X_j\\X_jB_0^{(\tau_j)}\end{bmatrix}\\
&=&
\rank\left(\begin{bmatrix}X_j&0\\0&X_j\end{bmatrix}\begin{bmatrix}
I_{2\tau_j}&0\\0 &B_0^{(\tau_j)}\end{bmatrix}\right)=2\rank X_j.
\end{array},\quad 1\leq j\leq m.
\]
Since $\rank P_{11}^{(1)}=\tau_1$, in $P_{11}^{(1)}$  there must be $\tau_1$ linearly independent rows
$i_1<\dots< i_{\tau_1}$. Then $\I_1=(i_1,\ldots, i_{\tau_1})\in Q_{\tau_1, r}= Q_{\tau_1-\tau_0, r}$ and
$\rank P(\I_1,2\tau_1)^\diamond=\rank P_{11}^{(1)}(\I_1,:)^\diamond=
\rank \begin{bmatrix}P_{11}^{(1)}(\I_1,:)\\P_{11}^{(1)}(\I_1:)B_0^{(\tau_1)}\end{bmatrix}=2\tau_1$.

Now, $\rank\begin{bmatrix}P_{11}^{(1)} & P_{11}^{(2)}\end{bmatrix}=\rank X_2=\tau_2$. Thus, in $P_{11}^{(2)}$ 
there must be $\tau_2-\tau_1$ rows $i_{\tau_{1}+1}<i_{\tau_{1}+2}<\dots <i_{\tau_{2}}$ such that the rows
$i_1<\dots< i_{\tau_1}, \, i_{\tau_{1}+1}< \dots < i_{\tau_{2}}$ of $X_2$ are linearly independent. 
Put $\I_2=(i_1,\ldots, i_{\tau_2})$. Then $I_1\subseteq I_2$, $I_2\setminus I_1=
(i_{\tau_{1}+1}, \dots, i_{\tau_{2}})\in  Q_{\tau_2-\tau_{1},r}$, and
$\rank P(\I_2,2\tau_2)^\diamond=\rank X_2(\I_2,:)^\diamond= 
\rank \begin{bmatrix}X_2(\I_2,:)\\X_2(\I_2:)B_0^{(\tau_2)}\end{bmatrix}=2\tau_2$.
Continuing the process, we can obtain $m$ sequences $\I_1, \dots, \I_m$ satisfying \eqref{eqIsubsetr}, 
\eqref{eqIsetminus} and \eqref{eqIrank}. \hfill $\Box$

\begin{definition}\label{def.admissible}
Given  $P\in \PP_{(\whW;\rr)}$, let
$\I_i$, $1\leq i \leq m$, be sequences of indices satisfying \eqref{eqIsubsetr}, \eqref{eqIsetminus} and \eqref{eqIrank}. 
Then we  say that $\uI=(\I_1, \dots, \I_{m})$ is an {\em admissible sequence of indices} for $P$.
\end{definition}

\begin{proposition}\label{propsameadmissible}
Let $P, \wh P\in \PP_{(\whW;\rr}$ be matrices  such that $\wh{P} {\sim} P$ and let $\uI=(\I_1, \dots, \I_{m})$
be an admissible sequence of indices for $P$.
Then $\uI$ is also an admissible sequence of indices for $\wh{P}$.
\end{proposition}

\noindent
{\bf Proof.} The proof is analogous to that of Proposition \ref{propsameadmissiblerw}.
\hfill $\Box$

\bigskip
Let
\begin{equation}\label{eq.defAshW}
\A_{\whW}=\{\uI=(\I_1, \dots, \I_{m})\; : \; \I_j  \mbox{ satisfies \eqref{eqIsubsetr} and \eqref{eqIsetminus}, } 1\leq j \leq m\}.
\end{equation}
For $\uI \in \A_{\whW}$, $\U_{\uI}$ denotes the open subset of 
 $\PP_{(\whW; \rr)}$ formed by the matrices  of $\PP_{(\whW;\rr)}$
with $\uI$ as an admissible sequence of indices.

\medskip
As in the case of only one real eigenvalue (Section \ref{subsecredreal}) we introduce now the notion of reduced  form
of a matrix in $\PP_{(\whW; \rr)}$.

Let $P=\begin{bsmallmatrix} P_1\\P_2\\\vdots\\P_k\end{bsmallmatrix}\in\PP_{(\whW; \rr)}$ with $P_1=\begin{bmatrix}
P_{11}& P_{12}&\cdots &P_{1m}\end{bmatrix}$ and $P_{1j}\in\RR^{r\times 2w_j}$, $1\leq j\leq m$. Let
$\uI=(\I_1, \dots, \I_{m})$ be an admissible sequence of indices for $P$, with $\I_j=(i_1, \dots, i_{\tau_j})$, $1\leq j\leq m$.
A matrix
$R=\begin{bsmallmatrix} R_1\\R_2\\\vdots\\R_k\end{bsmallmatrix}\in\PP_{(\whW; \rr)}$ with $R_1=\begin{bmatrix}
R_{11}& R_{12}&\cdots &R_{1m}\end{bmatrix}$ and $R_{1j}\in\RR^{r\times 2w_j}$, $1\leq j\leq m$, is said to be
a \textit{$\wC_{\whW}$-reduced form of $P$ with respect to $\uI$} if

\[
R_{11}(\I_m,:)^\diamond=\begin{bmatrix} I_{2\tau_1} & 0 & 0 & \cdots & 0\\R_{21}^{(1)\diamond}& 
I_{2(\tau_2-\tau_1)} & 0 &\cdots & 0\\
R_{31}^{(1)\diamond }& R_{32}^{(1)\diamond}& I_{2(\tau_3-\tau_2)} &\cdots & 0\\\vdots&\vdots& \vdots &\ddots& \vdots\\
R_{m1}^{(1)\diamond}& R_{m2}^{(1)\diamond}& R_{m3}^{(1)\diamond} & \cdots & I_{2(\tau_m-\tau_{m-1})}\end{bmatrix}
\]
and for $j=2,3,\ldots, m$,
\[
\begin{array}{l}
R_{1j}(\I_m,:)^\diamond=\\
\begin{blockarray}{ccccccc}
 \mbox{\scriptsize $2\tau_1$}&
 \mbox{\scriptsize $2(\tau_2-\tau_1)$}& \mbox{\scriptsize $2(\tau_3-\tau_2)$}&\cdots &
 \mbox{\scriptsize $2(\tau_{m-j}-\tau_{m-j-1})$}&
\mbox{\scriptsize $2(\tau_{m-j+1}-\tau_{m-j})$}\\
\begin{block}{[cccccc]c}
0 & 0 & 0&\cdots &0&0&\mbox{\scriptsize $2\tau_1$}\\
\vdots &\vdots&\vdots& & \vdots & \vdots &\vdots\\
0& 0 & 0&\cdots &0& 0&\mbox{\scriptsize $2(\tau_j-\tau_{j-1})$}\\
R_{j+1\, 1}^{(j)\diamond}&0&0& \cdots &0&0&\mbox{\scriptsize $2(\tau_{j+1}-\tau_j)$}\\
R_{j+2\, 1}^{(j)\diamond} & R_{j+2\, 2}^{(j)\diamond} &0& \cdots & 0 & 0&\mbox{\scriptsize $2(\tau_{j+2}-\tau_{j+1})$}\\
\vdots&\vdots&\vdots&&\vdots &\vdots&\vdots\\ 
R_{m1}^{(j)\diamond} & R_{m2}^{(j)\diamond}&R_{m3}^{(j)\diamond}&\cdots &R_{m\,m-j}^{(j)\diamond} &0&\mbox{\scriptsize $2(\tau_{m}-\tau_{m-j})$}\\
\end{block}
\end{blockarray}
\end{array}
\]
\begin{theorem}\label{theofrcomplex}
Let $P\in \PP_{(\whW; \rr)}$ and let $\uI=(\I_1, \dots, \I_{m})$ be an admissible sequence of indices for $P$. Then $P$ is
$\wC_{\wh W}$-equivalent to a unique $\wC_{\wh W}$-reduced form $P^{(\re)}\in \PP_{(\whW; \rr)}$ with respect to $\uI$.
\end{theorem}
{\bf Proof.}
The proof is analogous to that of Theorem \ref{theofrreal}, replacing $P_{ik}^{(j)}$ by $P_{ik}^{(j)\diamond}$ in
\eqref{eqfrlrw}
and using the elementary matrices $\widehat{Y_{I}}$ and $\widehat{Y_{II,i,k}^{(j)}}$ instead of $Y_{I}$ and $Y_{II,i,k}^{(j)}$. 
\hfill $\Box$


\begin{rem}\label{rem.npredform}
The number of parameters of $P^{(\re)}$ is $2sr-N_{\whW}=\dim \PP_{(\wh W,\rr)}-\dim \wC_{\wh W}$.
  \end{rem}

\subsection{Local parameterization and local system of coordinates of   $\PP_{(A;\rr)}/\wC_A$}
\label{subsecredgeneral}

We consider now the general case: let $\ualpha:\alpha_1(s)\mid\cdots\mid \alpha_n(s)$ be monic polynomials such that
$\sum_{i=1}^n\deg(\alpha_i(s))=n$ and assume that \eqref{eq.primefacalfa} is the prime factorization of
$\alpha_{n-i+1}(s)$. Let $A$ be the associated real Weyr canonical form of  \eqref{eq.AWR}.
Assume also that condition \eqref{eqnecss} holds true. It follows from Remark \ref{rem.deqn} that
$\PP_{(A; \rr)}\neq \emptyset$.


Given $P,  \wh P\in \PP_{(A; \rr)}$, we can partition $P$ and $\wh P$ in the form
$$
P=\begin{bmatrix}P_1 &\dots &P_p&P_{p+1}&\dots &P_{p+q}\end{bmatrix},
\;
 \wh P=\begin{bmatrix}\wh P_1 &\dots &\wh P_p&\wh P_{p+1}&\dots &\wh P_{p+q}\end{bmatrix},
$$
with $P_i, \wh P_i\in \PP_{(W_i; \rr)}$, $1\leq i \leq p$ and
$P_i, \wh P_i\in \PP_{(\whW_i; \rr)} $, $p+1\leq i \leq p+q$.

Then
$P\stackrel{\wC_A}{\sim}\wh P$ if and only if
$P_i\stackrel{\wC_{W_i}}{\sim}\wh P_i$, $1\leq i \leq p$, and
$P_i\stackrel{\wC_{\whW_i}}{\sim} \wh P_i$, $p+1\leq i \leq p+q$.

\begin{definition}\label{defmultiindex}
With the above notation, let $P=\begin{bmatrix}P_1 &\dots &P_p&P_{p+1}&\dots &P_{p+q}\end{bmatrix}\in \PP_{(A; \rr)}$
and let $\uI^{(i)}$ be an admissible sequence of indices for $P_i$, $1\leq i \leq p+q$.
Then we say that $\uI=(\uI^{(1)}, \dots, \uI^{(p+q)})$ is a {\em multi-index} for $P$.
\end{definition}

\begin{definition}\label{deffrgeneral}
Let $P=\begin{bmatrix}P_1 &\dots &P_p&P_{p+1}&\dots &P_{p+q}\end{bmatrix}\in \PP_{(A; \rr)}$ and let
 $\uI=(\uI^{(1)}, \dots, \uI^{(p+q)})$ be a {\em multi-index} for $P$. Let
 $ P_i^{(\re)}$ be the $\wC_{W_i}$-reduced form of $P_i$ with respect to $\uI^{(i)}$, $1\leq i \leq p$,  and 
let $ P_i^{(\re)}$ be the $\wC_{\wh W_i}$reduced form of $P_i$ with respect to $\uI^{(i)}$ , $p+1\leq i \leq p+q$.
Then
\[
P^{(\re)}=\begin{bmatrix}P^{(\re)}_1 &\dots &P^{(\re)}_p&P^{(\re)}_{p+1}&\dots &P^{(\re)}_{p+q}\end{bmatrix}
\]
is said to be the $\wC_{A}$-{\em reduced form of $P$}.
\end{definition}

Let $s^{(i)}=\sum_{j=1}^{m_{i,1}}w_{i, j}$, $1\leq i \leq p+q$,
 $N_{i}=\dim \wC_{W_i}$, $1\leq i \leq p$, and
 $N_{i}=\dim \wC_{\whW_i}$, $p+1\leq i \leq p+q$,
then $n=\sum_{i=1}^{p}s^{(i)}+2\sum_{i=p+1}^{p+q}s^{(i)}$ and
$N=\dim\wC_A=\sum_{i=1}^{p+q}N_i$.
The number of parameters in $P^{(\re)}$ is
$\sum_{i=1}^{p}(s^{(i)}r-N_i)+\sum_{i=p+1}^{p+q}(2s^{(i)}r-N_i)= nr-N$.

\bigskip
Recalling the definitions of $\A_{W}$ and $\A_{\whW}$ in \eqref{eq.defAsW} and \eqref{eq.defAshW}, respectively, let
$$
\A_A=\A_{W_1}\times \dots \times \A_{W_p}\times\A_{\whW_{p+1}}\times  \dots \times \A_{\whW_{p+q}}.
$$
Given $\uI=(\uI^{(1)}, \dots, \uI^{(p+q)})\in \A_A$, let
$\U_{\uI}=\U_{\uI^{(1)}}\times \dots\times \U_{\uI^{(p+q)}}$. Note that the matrices in $\U_{\uI^{(j)}}$ are,
for $1\leq j\leq p+q$, full column rank matrices; however, there may be matrices in $\U_{\uI}$ which do not have full
column rank. Thus, we must define $\V_{\uI}=\PP_{(A; \rr)}\cap \U_{\uI}$. This is an open subset of $\PP_{(A; \rr)}$.

Let $\R_{\uI^{(i)}}$ denote the subset of $\U_{\uI^{(i)}}$ formed by the matrices $P_i^{(\re)}\in \U_{\uI^{(i)}}$ in
$\wC_{W_i}$-reduced form or $\wC_{\wh W_i}$-reduced form with respect to $\uI^{(i)}$ according as $1\leq i \leq p$
or $p+1\leq i \leq p+q$. Now, for $1\leq i\leq p$, let 
$\nu_i:\; \RR^{s^{(i)}r-N_i}\; \longrightarrow \; \R_{\uI^{(i)}}$ be defined as follows: Let $x\in\RR^{s^{(i)}r-N_i}$. Use the
$(\tau_2-\tau_1)\tau_1$ first elements of $x$ to successively construct, row by row, a $(\tau_2-\tau_1)\times\tau_1$
matrix  and call it $R_{2,1}^{(1)}$ as in \eqref{eq.defR11}. Then use the following $(\tau_3-\tau_2)\tau_1$ elements of $x$ 
to successively construct, row by row, a $(\tau_3-\tau_2)\times\tau_1$ matrix  and call it $R_{31}^{(1)}$ as in \eqref{eq.defR11}.
Use the following $(\tau_3-\tau_2)(\tau_2-\tau_1)$ to construct the matrix $R_{32}^{(1)}$ of \eqref{eq.defR11} .
Use the same rules to successively construct the remaining blocks of the lower block-triangular matrix $R_{11}(\I_{m_i},:)$
of \eqref{eq.defR11}. Then use the remaining elements of $x$ to construct the matrices $R_{1j}(\I_{m_i},:)$ of 
\eqref{eq.defR1j}, $2\leq j\leq m_i-1$ (note that $R_{1m_i}(\I_{m_i},:)=0$). Next, for $p+1\leq i\leq p+q$ define 
$\nu_i:\; \RR^{2s^{(i)}r-N_i}\; \longrightarrow \; \R_{\uI^{(i)}}$  to map $x\in\RR^{2s^{(i)}r-N_i}$ into 
$(R_{11}(\I_{m_i},:)^\diamond, R_{12}(\I_{m_i},:)^\diamond,\ldots, R_{1m_i-1}(\I_{m_i},:)^\diamond)$ where these
matrices are constructed as in the previous case.

Define now $\R_{\uI}=\R_{\uI^{(1)}}\times \dots\times \R_{\uI^{(p+q)}}$ and 
$\nu_{\uI}:\; \RR^{nr-N}\; \longrightarrow \; \R_{\uI}$ as $\nu(x)=(\nu_1(x_1),\ldots, \nu_p(x_p),\nu_{p+1}(x_{p+1}),
\ldots, \nu_{p+q}(x_{p+q}))$ where $x=(x_1,\ldots, x_p,x_{p+1},\ldots,$ $x_{p+q})$ and $x_i\in \RR^{s^{(i)}r-N_i}$ if
$1\leq i\leq p$ and $x_i\in \RR^{2s^{(i)}r-N_i}$ if $p+1\leq i\leq q$ . It is plain that $\nu_{\uI}$ is a diffeomorphism and we
can identify $\R_{\uI}$ with  $\RR^{nr-N}$.

Let
$\W_{\uI}=\nu_{\uI}^{-1}(\R_{\uI}\cap \PP_{(A; \rr)})$. Then $\W_{\uI}$ is an open subset of $\RR^{nr-N}$. 
If $x \in \W_{\uI}$, we will will denote
$\nu_{\uI}(x)$ by $P_x^{(\re)}$
\medskip

Let $\pi$ be the submersion $\pi:\;\PP_{(A; \rr)}\; :\; \longrightarrow\; \PP_{(A; \rr)}/\wC_A$. Then $\pi$ is an open map
(see, for example, \cite[Theorem 7.16]{Lee03}) and so
$\wt{\V_{\uI}}=\pi(\V_{\uI})$ is an open subset of $\PP_{(A; \rr)}/\wC_A$.

\begin{theorem}\label{theopar}
The map $\psi_{\uI}:\;\W_{\uI}\, \longrightarrow \, \wt\V_{\uI}$, defined by
$\psi_{\uI}(x)=\wt{P_x^{(\re)}}$, where $\wt{P_x^{(\re)}}$ is the orbit of $P_x^{(\re)}$ under the action of $\wC_A$, 
is a diffeomorphism.
\end{theorem}
{\bf Proof.} It is clear that $\psi_{\uI}$ is well defined and  bijective.

On one hand,  the map $\varphi_{\uI}:\;\W_{\uI}\, \longrightarrow \, \V_{\uI}$ defined by $\varphi_{\uI}(x)=P_x^{(\re)}$ is 
differentiable. Hence, the map $\psi_{\uI}=\pi\mid_{\V_{\uI}}\circ \varphi_{\uI}$ is also  differentiable.

On the other  hand, the map  $\alpha_{\uI}:\;\V_{\uI}\, \longrightarrow \, \R_{\uI}\cap \PP_{(A; \rr)}$ defined by
$\alpha_{\uI}(P)=P^{(\re)}$, where $P^{(\re)}$
is the $\wC_{A}$-reduced form of $P$ with respect to $\uI$, is differentiable. Hence, the map
$\eta_{\uI}= \nu_{\uI}^{-1}\circ \alpha_{\uI}:\; \V_{\uI}\; \longrightarrow \; \W_{\uI}$ is also differenciable.
Since $\eta_{\uI}=\psi_{\uI}^{-1}\circ \pi\mid_{\V_{\uI}}$, by Proposition 7.17 of \cite{Lee03}, we conclude that
$\psi_{\uI}^{-1}$ is differentiable.
\hfill $\Box$

\begin{rem}
The map $\psi_{\uI}$ defined in Theorem \ref{theopar} is a local parameterization and
$\psi_{\uI}^{-1}$ is a local system of coordinates for $\PP_{(A; \rr)}/\wC_{A}$.
\end{rem}

\section{Local parameterization and local system of coordinates of $\HH_{(F,G)}$}
\label{secparameterization}

Assume that we are given  a sequence of monic polynomials $\ualpha:\alpha_1(s)\mid\cdots\mid\alpha_n(s)$ such that 
$\sum_{i=1}^n\deg(\alpha_i)=n$ and  a controllable pair $(F, G) \in  \RR^{n\times n}\times\RR^{n\times m}$ in
$p$-Brunovsky canonical form (cf. \eqref{eq.FpGp}) with $\rr: r_1\dots \geq r_{k}>0=r_{k+1}=\dots =r_n$ as Brunovsky indices.
Let $\kk: k_1\geq \dots \geq k_{r}>0=k_{r+1}=\dots =k_m$ be its 
controllability indices and assume that condition \eqref{eqnecss} holds.
Then, by Proposition \ref{prop.necss},
$\HH_{(F, G)}\neq \emptyset$. Assume also (Remark \ref{remrankr}) that $G=\begin{bmatrix}G_1&0\end{bmatrix}$,
$ G_1\in \RR^{n\times r}$ and $\rank G_1 =r$. Also, let $\rr: r_1\dots \geq r_{k}>0=r_{k+1}=\dots =r_n$ be the
Brunovsky  indices of $(F, G)$ and let $A\in\OO(\ualpha)$.
Then the map $ \phi: \, \PP_{(A;\rr)}/\wC_A   \, \longrightarrow \,    \HH_{(F,G_1)}$ 
defined in Theorem \ref{teodifeo}, is a diffeomorphism.


\medskip
 
 Assume that $A$ is in real Weyr canonical form (cf. \eqref{eq.AWR}) and with
the notation of Section \ref{subsecredgeneral}, let $\uI\in \A_A$ and $\wh {\V}_I=\phi(\wt{\V}_{\uI})$.
Then $\wh {\V}_I$ is an open subset of $\HH_{(F,G_1)}$. Hence, if 
 $\psi_{\uI}:\;\W_{\uI}\, \longrightarrow \, \wt{\V}_{\uI}$
 is the diffeomorphism defined in Theorem \ref{theopar}, then
 $\alpha_{\uI}= \phi\circ \psi_{\uI}:\;\W_{\uI}\, \longrightarrow \, {\wh\V}_{\uI}$ is also a diffeomorphism.

 Recall that $\HH_{(F, G)}=\HH_{(F, G_1)}\times \RR^{(m-r)\times n}$. Then
 $\wh {\V'}_I=\wh {\V}_I\times \RR^{(m-r)\times n}$ is an open subset of $\HH_{(F,G)}$ and
 $$
\begin{array}{rccc}
\alpha'_{\uI}:&\W_{\uI}\times \RR^{(m-r)\times n}& \longrightarrow &{\wh\V'_{\uI}}\\
&(x, K_2)& \mapsto & \begin{bsmallmatrix}\alpha_{\uI}(x)\\\\K_2\end{bsmallmatrix}
\end{array}
$$
is a diffeomorphism.
Therefore, $\alpha'_{\uI}$
 is a local parameterization and
$\psi_{\uI}^{-1}$ is a local system of coordinates of $\HH_{(F,G)}$.

\begin{rem}
If $m=n$ and  $r_1=n$, then  $(F, G)$ is feedback equivalent to $(\nsc, I_n)$ and 
$
\HH_{(F, G)}=\OO(\ualpha).
$
In this case we obtain a parameterization of $\OO(\ualpha)$.
\end{rem}

We finish with an example illustrating the whole procedure to obtain a parametrization of $\HH_{(F,G)}$ when
$(F,G)$ is in $p$-Brunovsky canonical form.

\begin{example}\label{exfinalw}{\rm
Let $n=5$, $\alpha_1(s)=\alpha_2(s)=\alpha_3(s)=1$, $\alpha_4(s)=s$, $\alpha_5(s)=s^2(s^2+1)$, 
and let  $\rr=(r_1, r_2, r_3)=(2, 2, 1)$.
Then $\ualpha=1\mid 1\mid 1\mid s\mid s^2(s^2+1)$ and the Segre characteristic of any $A\in\OO(\ualpha)$ for the
eigenvalue $0$ is $(2,1)$ and that for the eigenvalues $i$ and $-i$ is $(1)$. Thus their Weyr characteristics are
the conjugate partitions of $(2,1)$ and $(1)$, respectivley:
$$
w(0)=(2,1), \quad 
w(i)=w(-i)=(1).
$$
Also $\tau_1=2$, $\tau_2=1$ for the eigenvalue $0$ and $\tau_1=1$ for the eigenvalues $i$ and $-i$.

On the other hand, the controllability indices of $(F,G)$ are $\kk=(3,2)$. Therefore, $(3,2)\prec (4,1)$ and,
by Proposition \ref{prop.necss}, $\HH_{(F,G)}\neq \emptyset$.

Let
\[
A=\left[\begin{array}{cc|c||cc}
    0&0&1&0&0\\0&0&0&0&0\\\hline 0&0&0&0&0\\
\hline \hline
0&0&0&0&1\\0&0&0&-1&0\\
  \end{array}\right]
=\begin{bmatrix}W(0)&\nsc\\\nsc&\whW(i, -i)\end{bmatrix},
\]
and
\[
\begin{bmatrix}F&G\end{bmatrix}
  =\left[\begin{array}{cc|cc|c||cc}
      0&0&1&0&0&0&0\\0&0&0&1&0&0&0\\\hline
      0&0&0&0&1&0&0\\ 
      0&0&0&0&0&0&1\\ \hline
      0&0&0&0&0&1&0\\    
       \end{array}\right].
\]
Note that
$$w(0)\cup w(i)\cup w(-i)=(2,1,1,1)\prec \rr=(2, 2, 1).$$
It follows from  \eqref{eq.necssweyr} and Proposition \ref{prop.necss} that $\HH_{(F,G)}\neq\emptyset$ and
from Remark \ref{rem.deqn} that $\PP_{(A; \rr)}\neq \emptyset$.
 $$
  \PP_{(A; \rr})=\left\{\begin{bsmallmatrix}p_1\\p_2\\\hline\\p_1A\\p_2A\\\hline p_1A^2\end{bsmallmatrix}=
\left[\begin{smallarray}{cc|c||cc}
    p_{11}^{(1)}& p_{12}^{(1)}& p_{11}^{(2)}& p_{14}& p_{15}\\
    p_{21}^{(1)}& p_{22}^{(1)}& p_{21}^{(2)}& p_{24}& p_{25}\\\hline
    0&0&p_{11}^{(1)}& -p_{15}& p_{14}\\
    0& 0&p_{21}^{(1)}& -p_{25}& p_{24}\\\hline 
    0&0& 0& -p_{14}& -p_{15}\\  
       \end{smallarray}\right]\in \Gl(5)
\right\}.
$$
The map 
$$
\begin{array}{rccc}
\phi: &
\PP_{(A;\rr)}/\wC_A   &\longrightarrow &   \HH_{(F,G)}\\
&\wt P& \mapsto & \begin{bsmallmatrix}p_1A^{3}\\p_2A^{2}\\
\end{bsmallmatrix}P^{-1},
\end{array}
$$
is a diffeomorphism. The possible multi-indices for the matrices in  $\PP_{(A; \rr)}$ are $\uI=(\uI^{(1)}, \uI^{(2)})$ with
$\uI^{(1)}=((1), (1,2))$ or $\uI^{(1)}=((2), (2,1))$, and $\uI^{(2)}=(1)$ or $\uI^{(2)}=(2)$.

Assume that $\uI=(\uI^{(1)}, \uI^{(2)})$ with
   $\uI^{(1)}=((2), (2,1))$ and
$\uI^{(2)}=(1)$; i.e.,
$p_{21}^{(1)}\neq 0$, $\begin{bsmallmatrix}p_{11}^{(1)}&p_{12}^{(1)}\\p_{21}^{(1)}&p_{22}^{(1)}\end{bsmallmatrix}\in \Gl(2)$, 
and $\begin{bmatrix}p_{14}&p_{15}\end{bmatrix}\neq \begin{bmatrix}0&0\end{bmatrix}$.
Let $\V_{\uI}$ be the open subset of matrices in $\PP_{(A; \rr)}$ which admit
$\uI$ as a multi-index, and $\wh V_{\uI}=\phi\circ \pi(\V_{\uI})\subseteq\HH_{(F,G)}$.

In order to obtain  the $\wC_A$-reduced form with respect to $\uI$ of the matrices in $\V_{\uI}$ we proceed as follows.
Let $$P_1=\left[\begin{array}{cc|c}
     p_{11}^{(1)}& p_{12}^{(1)}& p_{11}^{(2)}\\
    p_{21}^{(1)}& p_{22}^{(1)}& p_{21}^{(2)}\\\hline
    0&0&p_{11}^{(1)}\\
    0& 0&p_{21}^{(1)}\\\hline 
    0&0& 0 \\  
       \end{array}\right], \quad P_2=\left[\begin{array}{cc}
     p_{14}& p_{15}\\
    p_{24}& p_{25}\\ \hline
    -p_{15}& p_{14}\\
      -p_{25}& p_{24}\\ \hline
     -p_{14}& -p_{15}\\
       \end{array}\right].
$$
According to Theorem \ref{theofrreal}, $P_1$ is equivalent to a unique matrix $P_1^{(\re)}$ such that
$$P_1^{(\re)}((2,1), :)=
\left[\begin{array}{cc|c}
   1&0&0\\
   p_{2,1}^{(\re)}&1&0
       \end{array}\right];\mbox{ 
i.e., }
P_1^{(\re)}
=\left[\begin{array}{cc|c}
    p_{21}^{(\re)}&1&0\\
    1&0& 0\\\hline
    0&0&p_{21}^{(\re)}\\
    0&0& 1\\ \hline
    0&0& 0\\
       \end{array}\right].   
$$
Therefore there exists a matrix $Y_1=\begin{bmatrix}d^{(1)}_{11}&d^{(1)}_{12}&d^{(2)}_{11}\\ 
0&d^{(1)}_{22}&d^{(2)}_{21}\\0&0&d^{(1)}_{11}
\end{bmatrix}\in \wC_{W(0)}$ such that $P_1Y_1=P_1^{(\re)}$. It follows from this that 
$p_{21}^{(\re)}=\frac{p_{11}}{p_{21}}$.

Analogously, there exists a matrix $Y_2=\begin{bmatrix}z_{11}&z_{12}\\-z_{12}&z_{11}
\end{bmatrix}\in \wC_{\whW(i,-i)}$ such that $P_2Y_2=
\left[\begin{smallarray}{cc}
     1& 0\\
      p_{24}^{(\re)}& p_{25}^{(\re)}\\ \hline
      0&1\\
    -p_{25}^{(\re)}& p_{24}^{(\re)}\\\hline
    -1&0
       \end{smallarray}\right]=P_2^{(\re)}$.
In fact, $Y_2=\begin{bsmallmatrix} p_{14} & p_{15}\\-p_{15} & p_{14}\end{bsmallmatrix}^{-1}$. So,
$z_{11}=\frac{p_{14}}{p_{14}^2+p_{15}^2}$, $z_{12}=-\frac{p_{15}}{p_{14}^2+p_{15}^2}$,
$p_{24}^{(\re)}=\frac{1}{p_{14}^2+p_{15}^2}(p_{14}p_{24}+p_{15}p_{25})$ and
$p_{25}^{(\re)}=\frac{1}{p_{14}^2+p_{15}^2}(p_{14}p_{25}-p_{15}p_{24})$.

Summarizing,
$$
P^{(\re)}=
\begin{bmatrix}p^{(\re)}_1\\p^{(\re)}_2\\\hline
p^{(\re)}_1A\\p^{(\re)}_2A\\\hline
p^{(\re)}_1A^2\end{bmatrix}=
\left[\begin{array}{cc|c||cc}
    p^{(\re)}_{21}& 1&0&1&0\\
    1&0&0& p^{(\re)}_{24}& p^{(\re)}_{25}\\\hline
    0&0&p^{(\re)}_{21}& 0&1\\
    0& 0&1& -p^{(\re)}_{25}& p^{(\re)}_{24}\\\hline 0&0&0& -1&0
       \end{array}\right].
$$
The free parameters of $P^{(\re)}$ are $p^{(\re)}_{21}$, $p^{(\re)}_{24}$ and $p^{(\re)}_{25}$ (recall that,
by Theorem \ref{theovd}, $\dim \HH(F,G)=nr-N$ where $n=\sum_{i=1}^n \deg(\alpha_i(s))$, $r=r_1$ and $N$
is given by \eqref{eqqN}; in this case $\dim \HH(F,G)=5\cdot 2-7=3$.
Since $P^{(\re)}$ must be invertible, the free parameters must satisfy $p^{(\re)}_{21}p^{(\re)}_{24}\neq 1$.
Then, recalling the definition of $\phi$ in  \eqref{eq.defphi},
$$
\begin{bmatrix}p^{(\re)}_1A^3\\ p^{(\re)}_2A^2\end{bmatrix}(P^{(\re)})^{-1}=
 \begin{bmatrix}
    0&0&\frac{1}{p^{(\re)}_{21}p^{(\re)}_{24}-1}&
    -\frac{p^{(\re)}_{21}}{p^{(\re)}_{21}p^{(\re)}_{24}-1}&
    \frac{p^{(\re)}_{21}p^{(\re)}_{25}}{p^{(\re)}_{21}p^{(\re)}_{24}-1}\\
     0&0&\frac{p^{(\re)}_{25}}{p^{(\re)}_{21}p^{(\re)}_{24}-1}&
    -\frac{p^{(\re)}_{21}p^{(\re)}_{25}}{p^{(\re)}_{21}p^{(\re)}_{24}-1}&
    p^{(\re)}_{24}+\frac{p^{(\re}_{21}{p^{(\re)}_{25}}^2}{p^{(\re)}_{21}p^{(\re)}_{24}-1}
  \end{bmatrix}.
  $$
Taking 
$\W_{\uI}=\left\{\begin{bsmallmatrix}x\\y\\z\end{bsmallmatrix}\in \RR^3\; : \; xy\neq 1\right\}$, then
$\W_{\uI}$ is an open set of $\RR^3$ and
$$\begin{array}{rccc}
  \alpha_{\uI}:& \W_{\uI} &
  \longrightarrow &  \wh V_{\uI} \\\
  &\begin{bsmallmatrix}x\\y\\z\end{bsmallmatrix}& \mapsto & \begin{bsmallmatrix}
    0&0&\frac{1}{xy-1}&-\frac{x}{xy-1}&\frac{xz}{xy-1}\\
     0&0&\frac{z}{xy-1}&-\frac{xz}{xy-1}&y+\frac{xz^2}{xy-1}   
   \end{bsmallmatrix},
\end{array}
$$
is a parameterization of $\wh V_{\uI}$.
}\hfill$\Box$
\end{example}
\section{Conclusions}
\label{secconclusions}
Given a sequence $\alpha_1(s)\mid\cdots\mid\alpha_n(s)$ of monic polynomials
with $\sum_{i=1}^n\deg(\alpha_i(s))$ $=n$ and a
controllable linear control system $(F, G)$, the geometry of the set $\HH_{(F,G)}$ of feedback matrices $K$
such that the state matrix of the closed loop system $F + GK$ has $\alpha_1(s)\mid\cdots\mid\alpha_n(s)$
as invariant polynomials, has been studied. It is proved that $\HH_{(F,G)}$ is a differentiable manifold diffeomorphic
to an orbit space by the action of a Lee group. Namely, the orbit space is an orbit space of
truncated observability matrices whose state matrix is fixed and has the given sequence of polynomials as invariant
polynomials; and the Lee group is the centralizer of that matrix. Then the dimension,  a local parametrization
and a local system of coordinates of $\HH_{(F,G)}$ are provided.

\appendix
\section{Proof of Proposition \ref{prop.parnemp}}
The first step of the proof is to prove that  $\PP_{(A;\rr)}\neq\emptyset$ if and only if there exist nonnegative 
integers $k'_1\leq k_1, \ldots, k'_r\leq k_r$ (recall that $(k_1,\ldots k_r)$ is the conjugate partition of
$(r_1,\ldots, r_k)$) such that
\begin{equation}\label{eq.necss3}
(k'_{\sigma(1)},\ldots, k'_{\sigma(r)})\prec(\deg(\alpha_n(s)),\ldots, \deg(\alpha_d(s)))
\end{equation}
where $(k'_{\sigma(1)},\ldots, k'_{\sigma(r)})$ is a permutation of $(k'_1,\ldots, k'_r)$ rearranged in nonincreasing
order; i. e., $k'_{\sigma(1)}\geq\cdots\geq k'_{\sigma(r)}$, and $\alpha_1(s)\mid \cdots\mid \alpha_d(s)$
are the invariant polynomials of $A$. In fact,
bearing in mind the relationship of the Antoulas' truncated and permuted observability matrices 
and the matrices of $\PP_{(A;\rr)}$, when $n=\sum_{i=1}^d r_i=d$,  $P\in\PP_{(A;\rr)}$ if and only if 
(using Antoulas' notation of \cite[Section 2.2]{Antou83}) $\left\{p_1, p_1A,\ldots, p_1A^{k_1-1},\right.$
 $\left.\ldots, p_r, p_rA,\ldots, p_rA^{k_r-1}\right\}$ form a \textit{nice basis} of $\RR^{d}$. In that case, 
 $k_1\geq\dots\geq k_r$ are said to be \textit{nice indices} of $(P_1,A)$. Hence $\PP_{(A;\rr)}\neq \emptyset$
 if and only if there exists a matrix $P_1\in\RR^{r\times d}$ such that $k_1\geq\dots\geq k_r$ are nice indices
 of $(P_1,A)$. Thus, when $n=d$, it follows from \cite[Corollary 2.7]{Za97} that 
 $\PP_{(A;\rr)}\neq \emptyset$ if and only if 
 \[
 (k_1,\ldots, k_r)\prec(\deg(\alpha_n(s)),\ldots, \deg(\alpha_d(s)))
\]
or, equivalenty, by Proposition \ref{prop.uncmajc},
\[
(\deg(\alpha_n(s)),\ldots, \deg(\alpha_d(s)))^\ast\prec (r_1,\ldots, r_k).
\]
 If $n=\sum_{i=1}^d r_i>d$ then some rows of $P\in\PP_{(A;\rr)}$ must be linear dependent on other rows.
 If we take the first $d$ linearly independent rows of $P$ then they form a nice basis of $\RR^d$ because if,
 for some indices $1\leq i\leq r$ and $0\leq q\leq k_i-1$,
 $p_iA^q$ linearly depends on the rows preceding it in $P$ then $p_iA^{q+1}$ also depend on the rows
 preceding it in $P$. Hence there are nonnegative integers $k'_1\leq k_1, \ldots, k'_r\leq k_r$ such that
 they are nice indices of $(P_1,A)$ and so they satisfy \eqref{eq.necss3}.
 And conversely, if there are indices $k'_1\leq k_1$,
 \ldots, $k'_r\leq k_r$ satisfying \eqref{eq.necss3} then there is $P'\in\PP_{(A;\rr')}$
 where $\rr'=(r'_1,\ldots, r'_k)$ is the conjugate partition of $\kk'=(k'_1,\ldots, k'_r)$. Then we can add the rows
 $p_iA^{k'_i},\ldots p_iA^{k_i-1}$, $1\leq i\leq r$, in the appropriate positions to obtain a matrix 
 $P\in\PP_{(A;\rr)}$.

 The second part of the proof is to show that \eqref{eq.kalpha} and \eqref{eq.wr} are equivalent.
Put $x=n-d=\sum_{j=1}^r k_j-\sum_{j=1}^d\deg(\alpha_j(s))$. Then, for $i\geq 1$,
\[
\begin{array}{l}
\sum_{j=i+1}^r k_j\geq \sum_{j=1}^{d-i}\deg(\alpha_j(s)) \Leftrightarrow 
n-\sum_{j=1}^i k_j\geq d-\sum_{j=d-i+1}^{d}\deg(\alpha_j(s))\\
 \Leftrightarrow 
x+\sum_{j=d-i+1}^{d}\deg(\alpha_j(s))\geq \sum_{j=1}^i k_j\\
 \Leftrightarrow 
x+\deg(\alpha_d(s))+\cdots+\deg(\alpha_{d-i+1}(s))\geq k_1+\cdots+k_i\\
\Leftrightarrow 
(k_1,\ldots,k_r)\prec(x+\deg(\alpha_d(s)),\deg(\alpha_{d-1}(s),\ldots, \deg(\alpha_1(s)))
\end{array}
\]
where we have used that $x+\sum_{i=1} ^d\deg(\alpha_i(s))=\sum_{i=1}^r k_i$.
Taking into account that $(r_1,\ldots, r_k)=(k_1,\ldots, k_r)^\ast$, $(w_1,\ldots, w_d)=
(\deg(\alpha_d(s)),$ $\ldots, \deg(\alpha_1(s)))^\ast$ and using item (ii) of Proposition \ref{prop.uncmajc},
we get
\[
(w_1, \dots, w_d)\cup(x)^\ast\prec (r_1, \dots, r_d).
\]
In conclusion, \eqref{eq.kalpha} and \eqref{eq.wr} are equivalent conditions.

Finally, we are to prove that $\PP_{(A;\rr)}\neq\emptyset$ if and only if \eqref{eq.wr} holds.

We have seen in the first step of the proof that $\PP_{(A; \rr)}\neq \emptyset$ if and only if there exist indices 
$k'_1\leq k_1, \ldots, k'_r\leq k_r$ satisfying \eqref{eq.necss3}
where $k'_{\sigma(1)}\geq\cdots\geq k'_{\sigma(r)}$.
Let us see that $k'_{\sigma(j)}\leq k_j$, $1\leq j \leq r$. 
Let $j\in \{1, \dots, r\}$ and assume that $k'_{\sigma(j)}> k_j$. Then
$k'_{\sigma(1)}\geq\cdots\geq k'_{\sigma(j)}>k_j\geq k'_j$. This means that
$j\not \in \{\sigma(1), \dots, \sigma(j)\}$ and so, there exists $\ell>j$ such that 
$\ell \in \{\sigma(1), \dots, \sigma(j)\}$. Thus,
 $k_j\geq k_\ell\geq k'_\ell\geq k'_{\sigma(j)}>k_j$, which is a contradiction.

Taking $\wt k_j=k'_{\sigma(j)}$,  $1\leq j \leq r$ we can conclude that 
$\PP_{(A; \rr)}\neq \emptyset$ if and only if there exist indices
$\wt k_1\geq \dots\geq  \wt k_r$ such that
 $\wt k_j\leq k_j$, $1\leq j \leq r$, and
 \begin{equation}\label{eq.wtkalpha}
 (\wt k_1, \dots, \wt k_r)\prec (\deg(\alpha_d(s)), \dots, \deg(\alpha_1(s))).
 \end{equation}
Equivalently, there exist nonnegative integers  $\wt r_1\geq \dots\geq  \wt r_d$ (the conjugate
partition of $(\wt k_1,\ldots, \wt k_r)$) such that
\begin{equation}\label{eq.wtrr}
\wt r_j\leq r_j, \quad 1\leq j \leq d,
\end{equation}
\begin{equation}\label{eq.wwtr}
(w_1, \dots, w_d)\prec (\wt r_1, \dots, \wt r_d).
\end{equation}
 Note that \eqref{eq.wtrr}  is equivalent to $\wt k_j\leq k_j$ because $\wt r_i=\#\{j: \wt k_j\geq i\}\leq
 \#\{j: k_j\geq i\}=r_i$ and by item (ii) of Proposition \ref{prop.uncmajc}, \eqref{eq.wwtr} is equivalent to
 \eqref{eq.wtkalpha}.

The last step of the proof is to show that condition \eqref{eq.wr} holds if and only if  the exist indices
$\wt r_1\geq \dots\geq  \wt r_d$  satisfying  \eqref{eq.wtrr} and \eqref{eq.wwtr}.

The ``if'' part es immediate: it follows from \eqref{eq.wtrr} and \eqref{eq.wwtr}
that $\sum_{j=1}^i w_j\leq \sum_{j=1}^i\wt r_j\leq \sum_{j=1}^ir_j,\quad 1\leq i\leq d$.

Conversely, assume that \eqref{eq.wr} holds and let $h=\min\{i : d\leq \sum_{j=1}^i r_j\}$. Then
 $\sum_{j=1}^{h-1} r_j<d\leq\sum_{j=1}^h r_j$. Define
\[
\begin{array}{ll}
\wt r_j=r_j, & 1\leq j \leq h-1,\\
\wt r_h=d-\sum_{i=1}^{h-1}r_j,\\
 \wt r_j=0, &h+1\leq j\leq d. 
\end{array}
 \]
Then $\wt r_h=d-\sum_{i=1}^{h-1}r_j\leq r_h\leq r_{h-1}=\wt r_{h-1}$. Therefore $\wt r_1\geq \dots\geq  \wt r_d$ 
 and  \eqref{eq.wtrr} holds. Since $\sum_{i=1}^{d}\wt r_j=d=\sum_{i=1}^{d}\deg(\alpha_j)$,
 \eqref{eq.wwtr}  follows from \eqref{eq.wr}.   \hfill$\Box$

\section*{Disclosure statement}

The authors report there are no competing interests to declare.

\bibliographystyle{acm}

\end{document}